\RequirePackage{ifpdf}
\ifpdf 
\documentclass[pdftex]{sigma}
\else
\documentclass{sigma}
\fi

\numberwithin{equation}{section}

\begin{document}

\newcommand{\sh}{{\rm sh}}
\newcommand{\ch}{{\rm ch}}

\newcommand{\De}{\Delta}
\newcommand{\de}{\delta}
\newcommand{\Z}{{\mathbb Z}}
\newcommand{\N}{{\mathbb N}}
\newcommand{\C}{{\mathbb C}}
\newcommand{\Cs}{{\mathbb C}^{*}}
\newcommand{\R}{{\mathbb R}}
\newcommand{\Q}{{\mathbb Q}}
\newcommand{\T}{{\mathbb T}}
\newcommand{\re}{{\rm Re}\, }
\newcommand{\im}{{\rm Im}\, }
\newcommand{\cW}{{\cal W}}
\newcommand{\cJ}{{\cal J}}
\newcommand{\cE}{{\cal E}}
\newcommand{\cA}{{\cal A}}
\newcommand{\cR}{{\cal R}}
\newcommand{\cP}{{\cal P}}
\newcommand{\cM}{{\cal M}}
\newcommand{\cN}{{\cal N}}
\newcommand{\cI}{{\cal I}}
\newcommand{\cMs}{{\cal M}^{*}}
\newcommand{\cB}{{\cal B}}
\newcommand{\cD}{{\cal D}}
\newcommand{\cC}{{\cal C}}
\newcommand{\cL}{{\cal L}}
\newcommand{\cF}{{\cal F}}
\newcommand{\cH}{{\cal H}}
\newcommand{\cS}{{\cal S}}
\newcommand{\cT}{{\cal T}}
\newcommand{\cU}{{\cal U}}
\newcommand{\cQ}{{\cal Q}}
\newcommand{\cV}{{\cal V}}
\newcommand{\cK}{{\cal K}}
\newcommand{\intR}{\int_{-\infty}^{\infty}}
\newcommand{\intI}{\int_{0}^{\pi/2r}}
\newcommand{\limp}{\lim_{\re x \to \infty}}
\newcommand{\limn}{\lim_{\re x \to -\infty}}
\newcommand{\limpn}{\lim_{|\re x| \to \infty}}
\newcommand{\diag}{{\rm diag}}
\newcommand{\Ln}{{\rm Ln}}
\newcommand{\Arg}{{\rm Arg}}
\newcommand{\pri}{^{\prime}}
\newcommand{\Ccs}{C_0^{\infty}(\R)}
\newcommand{\rE}{{\mathrm E}}
\newcommand{\rF}{{\mathrm F}}
\newcommand{\vh}{\hat{v}}
\newcommand{\aaa}{a_{+},a_{-}}
\newcommand{\iR}{\int_{\R}}

\newtheorem{lem}{Lemma}[section]
\newtheorem{theor}[lem]{Theorem}
\newtheorem{cor}[lem]{Corollary}
\newtheorem{prop}[lem]{Proposition}

\allowdisplaybreaks

\renewcommand{\thefootnote}{$\star$}

\renewcommand{\PaperNumber}{101}

\FirstPageHeading

\ShortArticleName{A Relativistic Conical Function and its Whittaker Limits}

\ArticleName{A Relativistic Conical Function\\ and its Whittaker Limits\footnote{This paper is a
contribution to the Special Issue ``Relationship of Orthogonal Polynomials and Special Functions with Quantum Groups and Integrable Systems''. The
full collection is available at
\href{http://www.emis.de/journals/SIGMA/OPSF.html}{http://www.emis.de/journals/SIGMA/OPSF.html}}}

\Author{Simon RUIJSENAARS}

\AuthorNameForHeading{S. Ruijsenaars}

\Address{School of Mathematics, University of Leeds, Leeds LS2 9JT, UK}
\Email{\href{mailto:siru@maths.leeds.ac.uk}{siru@maths.leeds.ac.uk}}
\URLaddress{\url{http://www.maths.leeds.ac.uk/~siru/}}

\ArticleDates{Received April 30, 2011, in f\/inal form October 23, 2011;  Published online November 01, 2011}

\Abstract{In previous work we introduced and studied a function $R(a_{+},a_{-},{\bf c};v,\hat{v})$ that is a generalization of the hypergeometric function ${}_2F_1$ and the Askey--Wilson polynomials. When the coupling vector ${\bf c}\in\C^4$ is specialized to $(b,0,0,0)$, $b\in\C$, we obtain a function $\cR (a_{+},a_{-},b;v,2\hat{v})$ that generalizes the conical function specialization of ${}_2F_1$ and the $q$-Gegenbauer polynomials. The function $\cR$ is the joint eigenfunction of four analytic dif\/fe\-ren\-ce operators associated with the relativistic Calogero--Moser system of $A_1$ type, whereas the function $R$ corresponds to $BC_1$, and is the joint eigenfunction of four hyperbolic Askey--Wilson type dif\/ference operators. We show that the $\cR$-function admits f\/ive novel integral representations that involve only four hyperbolic gamma functions and plane waves. Taking their nonrelativistic limit, we arrive at four representations of the conical function. We also show that a limit procedure leads to  two commuting relativistic Toda Hamiltonians and two commuting dual Toda Hamiltonians, and that a similarity transform of the function $\cR$  converges to a joint eigenfunction of the latter four dif\/ference operators.}


\Classification{33C05; 33E30; 39A10; 81Q05; 81Q80}

\tableofcontents

\section{Introduction}
This article may be viewed as a continuation of our previous work on a `relativistic' generalization~$R$ of the Gauss hypergeometric function $_2F_1$, introduced in~\cite{I}. The latter paper and two later parts in a series~\cite{II,III} will be referred to as I, II and III in the sequel. The def\/inition of the $R$-function in I is in terms of a contour integral, whose integrand involves eight hyperbolic gamma functions. (We review this in Section~\ref{section2}, cf.~\eqref{R}--\eqref{defsj}.)

In recent years, van de Bult~\cite{bult} tied in the $R$-function with the notion of modular double of the quantum group ${\mathcal U}_q(sl(2,\C))$, as def\/ined by Faddeev~\cite{fadd1}. As a spin-of\/f, he obtained a new representation of the $R$-function. Also,
van de Bult, Rains and Stokman~\cite{bust} have shown (among other things) that the 8-variable $R$-function
\begin{gather}\label{Rdef}
R(a_{+},a_{-},{\bf c};v,\hat{v}),\qquad a_{+}, a_{-},v,\hat{v}\in\C,\qquad a_{+}/a_{-}\notin (-\infty,0],\qquad {\bf c}\in\C^4,
\end{gather}
can be obtained as a limit of Spiridonov's 9-variable hyperbolic hypergeometric function~\cite{spir}. Their novel viewpoint leads to a third representation for the $R$-function. (See Proposition~4.20 and Theorem~4.21 in~\cite{bust} for the latter two representations.)

In this paper we are concerned with a 5-variable specialization of the $R$-function, def\/ined by
\begin{gather}\label{cRdef}
\cR(a_{+},a_{-},b;x,y)\equiv R(a_{+},a_{-},(b,0,0,0);x,y/2).
\end{gather}
Suitable discretizations of this function give rise to the $q$-Gegenbauer polynomials, whereas discretizations of the $R$-function yield the Askey--Wilson polynomials, cf.~I; moreover, the nonrelativistic limit of the $\cR$-function yields the conical function specialization of $_2F_1$. Hence it may be viewed as corresponding to the Lie algebra $A_1$, whereas the $R$-function can be tied in with $BC_1$.

The key new result of this paper concerning $\cR$ consists of the integral representation
\begin{gather}\label{key}
\cR(b;x,y)=\sqrt{\frac{\alpha}{2\pi}}\frac{G(2ib-ia)}{G(ib-ia)^2}\int_{\R}dz \frac{G(z+ (x-y)/2-ib/2)G(z- (x-y)/2-ib/2)}{G(z+ (x+y)/2+ib/2)G(z- (x+y)/2+ib/2)}.
\!\!\!
\end{gather}
Here and throughout the paper we use parameters
\begin{gather}\label{aconv}
\alpha\equiv 2\pi/a_+a_-,\qquad a\equiv (a_++a_-)/2,
\end{gather}
$G(\aaa;z)$ is the hyperbolic gamma function (cf.~Appendix~\ref{appendixA}), and  the dependence on~$a_{+}$,~$a_{-}$ is suppressed. (We shall often do this when no confusion can arise.) Furthermore, in~\eqref{key} we choose at f\/irst
\begin{gather}\label{first}
(\aaa,b,x,y)\in(0,\infty)^2\times(0,2a)\times\R^2.
\end{gather}

By contrast to the previous three integral representations following from I, \cite{bult} and~\cite{bust}, the integrand in~\eqref{key}
involves only four hyperbolic gamma functions. We also obtain several closely related representations that involve in addition plane waves, cf.~\eqref{Ri}--\eqref{Rv}. As will transpire in Section~\ref{section3}, upon using the f\/irst one~\eqref{key} of these novel representations (which we dub `minimal' representations) to introduce the $\cR$-function, it is possible to rederive in a more transparent and self-contained way a great many features that also follow upon specialization of the $R$-function theory, developed not only in~I, II and~III, but also in our later papers~\cite{anne} and~\cite{quad}. Moreover, special cases and limits of the $\cR$-function are far more easily obtained from the minimal representations than from the original integral representation of I or from the alternative representations following from~\cite{bult} and~\cite{bust}. (The integrands of these earlier representations involve at least eight hyperbolic gamma function factors.)

A survey of the results of I--III and~\cite{anne} can be found in~\cite{Rsurv}, but the def\/inition~\eqref{cRdef} of the $A_1$-analog of the ($BC_1$) $R$-function dates back to the more recent paper~\cite{quad}. In Section~\ref{section2}, we review in particular the pertinent results from~\cite{quad}. However, we have occasion to  add a lot more information that follows by specializing previous f\/indings concerning the $R$-function and related functions to their $A_1$ counterparts. This includes the asymptotic behavior and Hilbert space properties obtained in~II and~III, resp., which are adapted to the $A_1$ setting in Subsection~\ref{section2.2}, and the parameter shifts obtained in~\cite{anne}, which we focus on in~Subsection~\ref{section2.4}.  Moreover, in~\eqref{cRM} we detail the connection of the renormalized function
\begin{gather}\label{cRr}
\cR_r(\aaa,b;x,y)\equiv  \frac{G(ib-ia)}{G(2ib-ia)}\cR(\aaa,b;x,y),
\end{gather}
to the $A_1$ type functions $M(ma_{+}+na_{-};x,y)$, $m,n\in\Z$, which featured in our previous papers~\cite{GLFII} and~\cite{Hilb}. We present the proof of~\eqref{cRM} in~Subsection~\ref{section2.3}, together with various corollaries.

Altogether, Section~\ref{section2} invokes a considerable amount of information from our previous work. We have attempted to sketch this in such a way that the reader need only consult the pertinent papers for quite technical aspects (in case of doubt and/or inclination, of course). Even so, it is probably advisable to skim through Section~\ref{section2} at f\/irst reading, referring back to it when the need arises.

By contrast, Section~\ref{section3} (combined with Appendices~A and~C) is largely self-contained. Its starting point is a hyperbolic functional identity that f\/irst arose as a specialization of elliptic functional identities expressing the relation of certain Hilbert--Schmidt integral kernels to the elliptic $BC_1$ relativistic Calogero--Moser dif\/ference operators introduced by van Diejen~\cite{diej}. We need not invoke these identities (which can be found in~\cite{HSI}, cf.~also~\cite{noum}), since the relevant hyperbolic version is quite easily proved directly. The key point is that the hyperbolic identities can be rewritten in terms of two pairs of hyperbolic $A_1$-type relativistic Calogero--Moser dif\/fe\-rence operators $A_{\pm}(b_j;x)$, $j=1,2$, with distinct couplings $b_1$, $b_2$. The dif\/ference operators are given by
\begin{gather}\label{Apm}
A_{\de}(b;x)\equiv \frac{s_{\de}(x-ib)}{s_{\de}(x)}T^x_{ia_{-\de}}+\frac{s_{\de}(x+ib)}{s_{\de}(x)}T^x_{-ia_{-\de}},\qquad \de=+,-.
\end{gather}
Here, the translations are def\/ined on analytic functions by
\begin{gather}\label{Tdef}
(T^z_cf)(z)\equiv f(z-c),\qquad c\in\C^{*}.
\end{gather}
Also, throughout this paper we use the abbreviations
\begin{gather}\label{denot}
s_{\de}(z)=\sinh(\pi z/a_{\de}),\! \qquad c_{\de}(z)=\cosh(\pi z/a_{\de}),\! \qquad e_{\de}(z)=\exp(\pi z/a_{\de}),\! \qquad \de=+,-.\!\!
\end{gather}
Choosing $b_1=b$, $b_2=0$, an auxiliary function $B(b;x,y)$ can be def\/ined that satisf\/ies the eigenvalue equations (at f\/irst under certain restrictions on the $B$-arguments)
\begin{gather}\label{AdeB}
A_{\de}(b;x)B(b;x,y)=2c_{\de}(y)B(b;x,y),\qquad \de=+,-.
\end{gather}
 More specif\/ically, the function $B(b;x,y)$ is the Fourier transform of the hyperbolic kernel function, which is a product of four hyperbolic gamma functions.

When we write the integrand of the integral def\/ining $B$ as a product of two factors that involve only two hyperbolic gamma functions, we can use the Plancherel relation and the explicit Fourier transform formula for factors of this type (derived in Appendix~\ref{appendixC}) to obtain two new integral representations for $B$.  In particular, this leads to a function $C(b;x,y)$ given by
\begin{gather}\label{defC}
C(b;x,y)\equiv \sqrt{\frac{\alpha}{2\pi}}\int_{\R}dz \frac{G(z+ (x-y)/2-ib/2)G(z- (x-y)/2-ib/2)}{G(z+ (x+y)/2+ib/2)G(z- (x+y)/2+ib/2)}.
\end{gather}

Comparing~\eqref{defC}, \eqref{key} and~\eqref{cRr}, we read of\/f
\begin{gather}\label{cRC}
\cR_r(b;x,y)=G(ib-ia)C(b;x,y).
\end{gather}
However, we need a further study of the $C$-function~\eqref{defC} to arrive at a proof of this relation to the function $\cR_r$, as def\/ined originally by~\eqref{cRr} and~\eqref{cRdef}. Indeed, as already alluded to below~\eqref{first}, we can use~\eqref{defC} as a starting point to derive many features that $C$ and $\cR_r$ have in common.

In particular, the general analysis in~Appendix~B of~I can be applied to the integral on the r.h.s.\ of~\eqref{defC}, which yields a complete elucidation of the behavior of $C(b;x,y)$ under meromorphic continuation.
Moreover, via the A$\De$Es (analytic dif\/ference equations)~\eqref{AdeB} and the manifest invariance of $C$ under interchanging $x$ and $y$, it follows that
$C(b;x,y)$ is a joint eigenfunction of the four A$\De$Os (analytic dif\/ference operators)
\begin{gather}\label{4A}
A_{+}(b;x),\qquad A_{-}(b;x),\qquad A_{+}(b;y),\qquad
A_{-}(b;y),
\end{gather}
with eigenvalues
\begin{gather}\label{4ev}
2c_{+}(y),\qquad 2c_{-}(y),\qquad
2c_{+}(x),\qquad 2c_{-}(x),
\end{gather}
resp. This is also the case for $\cR_r(b;x,y)$ and, moreover, the equality~\eqref{cRC} can be shown for the special case $y=ib$ by a further application of Appendix~\ref{appendixC}. The general case then follows by a uniqueness argument already used in Subsection~\ref{section2.3}.

We reconsider the special $b$-values
\begin{gather}\label{bmn}
b_{mn}\equiv ma_{+}+na_{-},\qquad m,n\in\Z,
\end{gather}
in Subsection~\ref{section4.1}, inasmuch as they satisfy $b_{mn}\in(0,2a)$.  Indeed, the new Fourier transform representations in Section~\ref{section3} are only well def\/ined for $b\in(0,2a)$, but they can be explicitly evaluated by a residue calculation when $b$ is of this form. The key point is that the $G$-ratios in the integrand can then be written in terms of the hyperbolic cosines $c_{\pm}(z)$ by using the $G$-A$\De$Es~\eqref{Gades}. In principle, this yields again the functions $M(b_{mn};x,y)$ from~\cite{GLFII}, but we have not tried to push through a direct equality proof (as opposed to appealing to uniqueness).

Subsection~\ref{section4.2} deals with the nonrelativistic limit. Specializing the results of~I yields the hypergeometric function in terms of which the conical function can be expressed (cf.~Chapter~14 in~\cite{dlmf}). The f\/ive minimal representations~\eqref{Ri}--\eqref{Rv} of the $\cR$-function lead to four representations~\eqref{Rlimi}--\eqref{Rlimv} of the limit function. Rewriting them in terms of the conical function, three of these can be found in the literature (by looking rather hard). This is reassuring, since just as in I we were not able to get rigorous control on the nonrelativistic limits.

In order to describe the results of Section~\ref{section5}, we begin by recalling that in our paper~\cite{toda} we arrived at relativistic nonperiodic Toda $N$-particle systems by taking a limit of the relativistic hyperbolic Calogero--Moser $N$-particle systems. In this limit the self-duality of the latter is not preserved, inasmuch as the dual commuting Hamiltonians have a very dif\/ferent character from the def\/ining Hamiltonian and its commuting family. Specialized to the present context, this limit can be used to obtain a joint eigenfunction of two Toda Hamiltonians $H^{T}_{\pm}(\eta;x)$ and two dual Toda Hamiltonians $\hat{H}^{T}_{\pm}(\eta;y)$, with the real parameter $\eta$ playing the role of a coupling constant.

The limit transition proceeds in two stages. The f\/irst step is to set
\begin{gather}\label{gam}
b=a-i\gamma,\qquad \gamma\in\R.
\end{gather}
At the classical level the analogous $b$-choice still yields real-valued Hamiltonians with a well-def\/ined self-dual action-angle map and scattering theory~\cite{aa1}.  Correspondingly, the four reduced $N=2$ quantum Hamiltonians at issue here are still formally self-adjoint for this $b$-choice. (They are similarity transforms of the A$\De$Os~\eqref{4A} with a weight function factor.) Moreover, restricting attention to
\begin{gather}\label{real}
(a_{+},a_{-},x,y)\in(0,\infty)^2\times\R^2,
\end{gather}
their joint eigenfunction remains real-valued, although this reality property is no longer manifest: It hinges on a symmetry property under taking $b$ to $2a-b$, which translates into evenness in the parameter $\gamma$.

The next step is to substitute
\begin{gather}\label{Ts}
x\to x+\Lambda,\qquad \gamma\to \eta +\Lambda,
\end{gather}
and take $\Lambda$ to $\infty$. In this limit the Hamiltonians and their joint eigenfunction converge, whereas the dual Hamiltonians must be multiplied by a factor $e_{\de}(-\Lambda)$ to obtain a f\/inite limit. This can be understood from their $\Lambda$-dependent eigenvalues $2c_{\de}(x+\Lambda)$ following from the $x$-shift~\eqref{Ts}, cf.~\eqref{4ev}. Indeed, after multiplication by this renormalizing factor the eigenvalues have the f\/inite limits $e_{\de}(x)$, $\de=+,-$.

The f\/ive representations of the $\cR$-function give rise to four representations of the relativistic Toda eigenfunction F$^T(\eta;x,y)$, namely~\eqref{T1eq}, \eqref{T2eq}, \eqref{T3} and \eqref{T4}. Suitably paired of\/f, however, these dif\/ferent formulas express real-valuedness with~\eqref{real} in ef\/fect. Taking this into account, we wind up with two essentially dif\/ferent representations that are intertwined via the Plancherel formula for the Fourier transform. The key formula involved here is derived in~Corollary~\ref{corollaryC2}.

The results just delineated can be found in Subsection~\ref{section5.1}.
In Subsection~\ref{section5.2} we f\/irst study the asymptotic behavior of F$^T(\eta;x,y)$ for $x\to \pm \infty$ and $y\to\infty$. We then clarify the analyticity properties of F$^T(\eta;x,y)$ by introducing a similarity transform $\cH(x-\eta,y)$. Using the four representations~\eqref{cH1}--\eqref{cH4} of the latter, we  show that the function~$\cH(x,y)$ is holomorphic for $(x,y)\in\C^2$.

Subsection~\ref{section5.3} deals with the joint eigenfunction properties of F$^T(\eta;x,y)$ and its similarity transforms. Formally, these follow from those of the $\cR$-function. However, the Toda limit is not easy to control analytically, and the direct derivation of the eigenvalue equations is not too hard and quite illuminating.

Our results in Section~\ref{section5} have some overlap with earlier results by Kharchev, Lebedev and Semenov-Tian-Shansky~\cite{lebe},  who obtained functions closely related to
F$^T(\eta;x,y)$ from the viewpoint of harmonic analysis for Faddeev's modular double of a quantum group~\cite{fadd1}. The non\-relativistic nonperiodic Toda eigenfunctions are widely known as Whittaker functions, and meanwhile it has become customary to call eigenfunctions for $q$-Toda Hamiltonians Whittaker functions as well. In particular, $q$-Whittaker functions were introduced by Olshanetsky and Rogov for rank~1 (their work can be traced from~\cite{olro}) and by Etingof for arbitrary rank~\cite{etin}, and these functions have been further studied in various later papers (see e.g.~\cite{cher} and refe\-ren\-ces given there). We would like to stress that these functions are quite dif\/ferent from the ones at issue here and in~\cite{lebe}. The crux is that the former are only well def\/ined for $q$ not on the unit circle, whereas here and in~\cite{lebe} the eigenfunctions have a symmetric dependence on two generically distinct $q$'s, given by
\begin{gather}\label{2q}
q_{+}=\exp(i\pi a_{+}/a_{-}),\qquad q_{-}=\exp(i\pi a_{-}/a_{+}).
\end{gather}
This state of af\/fairs is closely related to the dif\/ferent character of the trigonometric gamma function (more widely known as the $q$-gamma function, with the restriction $|q|\ne 1$ being indispensable) and the hyperbolic gamma function (which depends on parameters $a_+$ and $a_-$ in the right half plane).

In Section~\ref{section6} we study the nonrelativistic limit of the representations of the relativistic eigenfunction, arriving at two distinct representations for the nonperiodic Toda eigenfunction that have been known for a long time. Just as for the relativistic case, its property of being also an eigenfunction for a dual Hamiltonian seems not to have been observed before. (These duality features are the quantum counterparts of duality features of the pertinent action-angle maps, f\/irst pointed out in~\cite{toda}.) To control one of the two pertinent limits, a novel limit transition for the hyperbolic gamma function is needed, whose proof is relegated to~Appendix~\ref{appendixB}.

\section[The $\cR$-function as a special case of the $R$-function]{The $\boldsymbol{\cR}$-function as a special case of the $\boldsymbol{R}$-function}\label{section2}

\subsection[The functions $\cR$ and $\cR_r$]{The functions $\boldsymbol{\cR}$ and $\boldsymbol{\cR_r}$}\label{section2.1}

The $R$-function \eqref{Rdef} is def\/ined as a contour integral over a variable $z$, with the $z$-dependence of the integrand encoded in a product of eight hyperbolic gamma function factors.
(See Appendix~\ref{appendixA} for a review of the relevant features of the hyperbolic gamma function.)
Specif\/ically,
with suitable restrictions on the eight variables, the $R$-function is given by
\begin{gather}\label{R}
R({\bf c};v,\hat{v})=\sqrt{\frac{\alpha}{2\pi}}\int_{\cal
C}F(c_0;v,z)K({\bf c};z)F(\hat{c}_0;\hat{v},z)dz.
\end{gather}
Here we have
\begin{gather}
\hat{c}_0\equiv (c_0+c_1+c_2+c_3)/2,
\\
\label{F}
F(d;y,z)\equiv
\frac{G(z\pm y+id-ia)}{G(\pm y+id-ia)},
\end{gather}
(with $f(w\pm y)$ denoting $f(w+y)f(w-y)$), and $K$ is given by
\begin{gather}\label{K}
K({\bf c};z)\equiv
\frac{1}{G(z+ia)}\prod_{j=1}^3\frac{G(is_j)}{G(z+is_j)},
\end{gather}
with new parameters
\begin{gather}\label{defsj}
s_1\equiv c_0+c_1-a_{-}/2,\qquad s_2\equiv c_0+c_2-a_{+}/2,\qquad s_3\equiv
c_0+c_3.
\end{gather}
Also, recall $a$ and $\alpha$ are def\/ined by~\eqref{aconv}.

We do not need the def\/inition of the contour $\cC$ for general variable choices (this is discussed in I and Section~4 of the survey~\cite{Rsurv}); instead we presently def\/ine $\cC$ for the cases at issue. For the special ${\bf c}$-choice in \eqref{cRdef} we can use the duplication formula~\eqref{dupl} to obtain
\begin{gather}
K((b,0,0,0);z)=\cK(b;z),
\end{gather}
where
\begin{gather}
\cK(b;z)\equiv \frac{1}{G(z+ia)}\frac{G(2ib-ia)}{G(ib-ia)}\frac{G(z+ib-ia)}{G(2z+2ib-ia)}.
\end{gather}
Using also the ref\/lection equation~\eqref{refl}
we deduce that $\cR$ is given by
\begin{gather}
\cR(b;x,y)     =     \left(\frac{\alpha}{2\pi}\right)^{1/2}\frac{G(2ib-ia)}{G(ib-ia)}G(\pm x-ib+ia)G(\pm y/2-ib/2+ia)
\nonumber \\
\phantom{\cR(b;x,y)     =}{}\times \int_{\cal
C}\frac{G(z\pm x+ib-ia)G(z+ib-ia)G(z\pm y/2+ib/2-ia)}{G(z+ia)G(2z+2ib-ia)}dz.\label{cRcont}
\end{gather}

For the variable choice that is most relevant for Hilbert space purposes, namely,
\begin{gather}\label{Hpar}
\aaa,b,x,y>0,
\end{gather}
the contour $\cC$ may be chosen equal to the real line in the $z$-plane, indented downwards near 0 so as to avoid a pole of $\cK(b;z)$. From~\eqref{Gpo}--\eqref{Gze} it follows that the poles of $\cK(b;z)$ are located on the imaginary axis at
\begin{gather}\label{Kpoles}
z-z_{kl}=0,\qquad z-z_{kl}+ib=ia_{+}/2, \, ia_{-}/2, \, 0,\qquad k,l\in\N.
\end{gather}
Thus they are above the contour, whereas the remaining $z$-poles of the integrand at
\begin{gather}
z+z_{kl}=\pm x-ib,\  \pm y/2-ib/2,\qquad  k,l\in\N,
\end{gather}
are below $\cC$.

From the above representation it is immediate that $\cR$ is symmetric under the interchange of the parameters $a_{+}$ and $a_{-}$:
\begin{gather}\label{Rmodinv}
\cR(a_{+},a_{-},b;x,y)=\cR(a_{-},a_{+},b;x,y).
\end{gather}
It is not at all clear, however, that $\cR$ is also symmetric under the interchange of the positions~$x$ and~$y$:
\begin{gather}\label{sd}
\cR(a_{+},a_{-},b;x,y)=\cR(a_{+},a_{-},b;y,x).
\end{gather}
 This self-duality feature follows in particular from a second relation between $\cR$ and $R$, namely,
\begin{gather}\label{cR2}
\cR(b;x,y)=R((b,b,b,b)/2;x/2,y).
\end{gather}
(This is equation~(4.8) in~\cite{quad}.) Indeed, this second ${\bf c}$-choice yields the same function $\cK(b;z)$ as the f\/irst one, so that substitution of~\eqref{R} (with the same contour $\cC$) now yields \eqref{cRcont} with $x$ and $y$ interchanged on the r.h.s.

There are two more ${\bf c}$-choices that lead from $R$ to $\cR$, namely, $(b,0,b,0)$ and $(b,b,0,0)$. Specif\/ically, from equations~(4.6) and (4.7) in~\cite{quad} we have
\begin{gather}\label{cR3}
\cR(\aaa,b;x,y)=R(a_{+},2a_{-},(b,0,b,0);x,y),
\\
\label{cR4}
\cR(\aaa,b;x,y)=R(2a_{-},a_{+},(b,b,0,0);x,y).
\end{gather}
From the def\/inition of the $R$-function we then obtain alternative integral representations for the $\cR$-function from which the self-duality property \eqref{sd} is manifest. (Indeed, since we have $c_0=\hat{c}_0=b$ for these two choices, the integrand is invariant under the interchange of $x$ and $y$.) On the other hand, the modular invariance property \eqref{Rmodinv} is not at all clear, since the integral representations involve the hyperbolic gamma function
 with $a_{-}$ replaced by $2a_{-}$. Using~\eqref{G2},
 they can be re-expressed in terms of the modular invariant function $G(\aaa;z)$. However, the resulting integrand is then still not modular invariant. Since it seems not to simplify and does not look
illuminating, we do not detail it any further.

The analyticity properties of the $R$-function are known in great detail from~Theorem~2.2 in I, cf.~also Section~4 in the survey~\cite{Rsurv}. Combining this theorem with the def\/inition~\eqref{cRdef} of $\cR$ and its self-duality property~\eqref{sd}, we deduce in particular that $\cR$ extends from the intervals \eqref{Hpar} to a meromorphic function in~$b$,~$x$ and~$y$, whose poles in~$x$ and~$y$ can only occur at the locations
\begin{gather}\label{cRpoles}
\pm z=2ia-ib+z_{kl},\qquad  z=x,y,\qquad k,l\in \N.
\end{gather}

We proceed to list further consequences of the $R$-function theory for $\cR$. Two features that are clear from each of the above integral representations are evenness and scale invariance (given scale invariance of $G$):
\begin{gather}\label{even}
\cR (b;x,y)=\cR(b;\de x, \de'y),\qquad \de,\de'=+,-,
\\
\label{scale}
\cR(\aaa,b;x,y)=\cR(\lambda  a_{+},\lambda a_{-}, \lambda b;\lambda x,\lambda y),\qquad \lambda>0.
\end{gather}
A less obvious feature is the explicit evaluation
\begin{gather}\label{cRnorm}
\cR(b;x,ib)=1.
\end{gather}
It follows from the formula
\begin{gather}
R({\bf c};v,i\hat{c}_0)=1,
\end{gather}
(cf.~equation~(3.26) in I or Section~6 in~\cite{Rsurv}), by using any of the four relations \eqref{cRdef}, \eqref{cR2}, \eqref{cR3}, \eqref{cR4}.
Def\/ining next
\begin{gather}
\cR_n(x)\equiv \cR(b;x,y_n),\qquad y_n\equiv ib+ina_{-},\qquad n\in\N,
\end{gather}
the eigenvalue A$\De$E (analytic dif\/ference equation) for $A_{+}(b;y)$ entails
\begin{gather}
 \frac{s_{+}(y_n-ib)}{s_{+}(y_n)}\cR_{n-1}(x)+\frac{s_{+}(y_n+ib)}{s_{+}(y_n)}\cR_{n+1}(x)=2c_{+}(x)\cR_n(x).
\end{gather}
In view of \eqref{cRnorm}, it follows from this that $\cR_n$ is of the form
\begin{gather}\label{qGeg}
\cR_n(x)=P_n(c_{+}(x)),
\end{gather}
where $P_n(z)$ is a polynomial in $z$ of degree $n$ and parity $(-)^n$. The relation of these polynomials to the $q$-Gegenbauer polynomials and to the  Askey--Wilson polynomials associated with the four relevant ${\bf c}$-choices is detailed at the end of Section~4 of~\cite{quad}.

The renormalized $\cR$-function~$\cR_r$ given by~\eqref{cRr}
is the counterpart of the renormalized $R$-function $R_r$ obtained from \eqref{R} by omitting the product $\prod_jG(is_j)$ in $K$, cf.~\eqref{K}. (To see this, use \eqref{refl} and \eqref{dupl}.) Clearly, it shares the features \eqref{Rmodinv}, \eqref{sd}, \eqref{even} and \eqref{scale} of the $\cR$-function, whereas \eqref{cRnorm} is replaced by
\begin{gather}\label{cRrnorm}
\cR_r(b;x,ib)=\frac{G(ib-ia)}{G(2ib-ia)}.
\end{gather}
The renormalizing factor in the function~$\cR_r$ ensures that it has no poles that are independent of $x$ and $y$, cf.~Theorem~2.2 in~I. More precisely, $\cR_r(\aaa,b;x,y)$ extends to a function that is meromorphic in the domain{\samepage
\begin{gather}\label{Dp}
D_+\equiv \big\{ (\aaa,b,x,y)\in \C^5 \mid \re a_+>0,\ \re a_->0\big\},
\end{gather}
and whose poles can only occur at the locations~\eqref{cRpoles}.}

It is not obvious, but true that for the special $b$-choices $b_{mn}$~\eqref{bmn}
we have an equality
\begin{gather}\label{cRM}
\cR_r(b_{mn};x,y)=M(b_{mn};x,y)+M(b_{mn};-x,y),
\end{gather}
where $M(b_{mn};x,y)$ is the function def\/ined at the end of Section~III in our paper~\cite{GLFII}. Therefore, $\cR_r(b;x,y)$ is the continuous (indeed, real-analytic) interpolation to arbitrary $b\in\R$ of  the function given by equation~(3.74) in~\cite{GLFII}, which is def\/ined only for the $b$-values~$b_{mn}$. (Note the latter are dense in $\R$ when the ratio $a_{+}/a_{-}$ is irrational.)

It will not cause surprise that in the free case we have
\begin{gather}\label{M00}
M(b_{00};x,y)=\exp(i\alpha xy/2).
\end{gather}
It would take us too far af\/ield, however, to detail all of the functions $M(b_{mn};x,y)$, $m,n\in\Z$. For our purposes it is enough to specify their general structure: They are `elementary' in the sense that they can be written
\begin{gather}\label{Mform}
M(b_{mn};x,y)=\exp(i\alpha xy/2)R_{mn}(e_{+}(x),e_{-}(x),e_{+}(y),e_{-}(y)),
\end{gather}
where $R_{mn}$ is a rational function of its four arguments, cf.~Section~III in~\cite{GLFII}. In Subsection~\ref{section2.4} we deduce this structure in another way (namely, by exploiting parameter shifts). Moreover, for the case where $m$ and $n$ are not both positive or both non-positive, this structure can be understood from the novel Fourier transform representations~\eqref{Rii}--\eqref{Rv}, cf.~Subsection~\ref{section4.1}.

We postpone the proof of the equality assertion \eqref{cRM} to Subsection~\ref{section2.3}. An ingredient of this proof is the asymptotic behavior of $\cR_r(b;x,y)$ as $x$ goes to $\infty$, and this is most easily obtained as a corollary of the asymptotics of a closely related function~$\rE(b;x,y)$, def\/ined by~\eqref{EcE}.

\subsection[The functions $\rE$ and $\rF$]{The functions $\boldsymbol{\rE}$ and $\boldsymbol{\rF}$}\label{section2.2}

The function $\rE(b;x,y)$ can be viewed as a specialization of the function denoted $\cE(\gamma;v,\hat{v})$ in II and~\cite{Rsurv}. The relation between $\gamma$ and ${\bf c}$ reads
\begin{gather}\label{gamc}
\gamma(\aaa,{\bf c})=(c_0-a,c_1-a_{-}/2,c_2-a_{+}/2,c_3).
\end{gather}
In particular, the `free' case ${\bf c}={\bf 0}$ yields
\begin{gather}\label{gamf}
\gamma_{\rm f}\equiv (-a,-a_{-}/2,-a_{+}/2,0).
\end{gather}
The switch from ${\bf c}$ to $\gamma$ is crucial for uncovering further symmetries: The function $\cE(\gamma;v,\hat{v})$ is invariant under $D_4$ transformations on $\gamma$ (i.e., permutations and even sign changes), cf.~II. (In~\cite{bust} this $D_4$ symmetry has been reobtained in a quite dif\/ferent way.) It is def\/ined by
\begin{gather}
\cE(\gamma;v,\hat{v})=\frac{\chi(\gamma)}{c(\gamma;v)c(\hat{\gamma};\hat{v})}R_r(\gamma;v,\hat{v}).
\end{gather}
Here, the generalized ($BC_1$) Harish-Chandra $c$-function is given by
\begin{gather}\label{c}
c(\gamma;v)\equiv \frac{1}{G(2v+ia)}\prod_{\mu =0}^3G(v-i\gamma_{\mu}),
\end{gather}
the dual of $\gamma$ by
\begin{gather}\label{dualgam}
\hat{\gamma}\equiv J\gamma,\qquad
J\equiv \frac{1}{2} \left( \begin{array}{rrrr}
1  &  1  &  1  &  1  \\
1  &  1  &  -1  &  -1  \\
1  &  -1  &  1  &  -1  \\
1  &  -1  &  -1  &  1
\end{array} \right) ,
\end{gather}
and the constant by
\begin{gather}\label{chigam}
\chi(\gamma) \equiv \exp\big(i\alpha [\gamma\cdot \gamma /4-\big(a_{+}^2+a_{-}^2+a_{+}a_{-}\big)/8]\big),\qquad
\alpha =2\pi/a_{+}a_{-}.
\end{gather}

Denoting the $\gamma$-vectors corresponding to the two one-parameter families
\begin{gather}
{\bf c}=(b,0,0,0), (b,b,b,b)/2,
\end{gather}
by $\gamma^{(1)}$,  $\gamma^{(2)}$, it is easy to check that
\begin{gather}\label{J12}
J\gamma^{(1)}=\gamma^{(2)}.
\end{gather}
Def\/ining
\begin{gather}\label{EcE}
\rE(b;x,y)\equiv \cE\big(\gamma^{(1)};x,y/2\big),
\end{gather}
a straightforward calculation (using the duplication formula \eqref{dupl}) yields
\begin{gather}\label{ER}
\rE(b;x,y)=\frac{\phi(b)}{c(b;x)c(b;y)}\cR_r(b;x,y),
\end{gather}
where we have introduced the constant
\begin{gather}\label{phib}
\phi(b)\equiv \exp(i\alpha b(b-2a)/4),
\end{gather}
and generalized ($A_1$) Harish-Chandra $c$-function
\begin{gather}\label{cb}
c(b;z)\equiv \frac{G(z+ia-ib)}{G(z+ia)}.
\end{gather}
Recalling \eqref{cR2} and using \eqref{J12}, it readily follows that we also have
\begin{gather}
\rE(b;x,y)=\cE\big(\gamma^{(2)};x/2,y\big).
\end{gather}

It involves more work to obtain the relations between $\rE$ and $\cE$ corresponding to \eqref{cR3} and~\eqref{cR4}. Setting
\begin{gather}
\gamma^{(3)}\equiv \gamma(a_{+},2a_{-},(b,0,b,0)),\qquad \gamma^{(4)}\equiv \gamma(2a_{-},a_{+},(b,b,0,0)),
\end{gather}
these are given by
\begin{gather}
\rE(\aaa,b;x,y)=\cE\big(a_{+},2a_{-},\gamma^{(3)};x,y\big)=\cE\big(2a_{-},a_{+},\gamma^{(4)};x,y\big).
\end{gather}
(These formulas amount to special cases of the doubling identity for the $\cE$-function obtained in Section~6 of~\cite{quad}.)

The relation \eqref{ER} between $\rE$ and $\cR_r$ yields a similarity transformation turning the A$\De$Os \eqref{4A} into $\cA_{\pm}(b;x)$, $\cA_{\pm}(b;y)$,
where
\begin{gather}\label{cAdef}
\cA_{\de}(b;z)\equiv c(b;z)^{-1}A_{\de}(b;z)c(b;z)=T^z_{ia_{-\de}}+V_{\de}(b;z)T^z_{-ia_{-\de}},\qquad \de=+,-,
\\
V_{\de}(b;z)\equiv \frac{s_{\de}(z+ib)s_{\de}(z-ib +ia_{-\de})}{s_{\de}(z)s_{\de}(z+ia_{-\de})}.
\end{gather}
From this it is easy to verify that these A$\De$Os are formally self-adjoint operators on $L^2(\R)$ (by contrast to the A$\De$Os \eqref{4A}), and that they are invariant under the transformation
\begin{gather}\label{bsym}
b \mapsto a_{+}+a_{-}-b.
\end{gather}
It is not obvious, but true that we also have
\begin{gather}\label{Esym}
\rE(b;x,y)=\rE(a_{+}+a_{-}-b;x,y).
\end{gather}
This symmetry property can be derived from \eqref{EcE} and the $D_4$ invariance of the $\cE$-function: We have
\begin{gather}\label{gam1}
\gamma^{(1)}=(b-a,-a_{-}/2,-a_{+}/2,0),
\end{gather}
so the map $b\mapsto 2a-b$ yields a $\gamma$-vector related to $\gamma^{(1)}$ by a sign f\/lip of the f\/irst and last component.

Combining~\eqref{ER} and~\eqref{cb} with the analyticity features of $\cR$ (cf.~the paragraph con\-taining~\eqref{cRpoles}), we deduce that $\rE(b;x,y)$ is meromorphic in $b$, $x$ and $y$, with $b$-independent pole locations
\begin{gather}
z=-2ia-z_{kl},\qquad z=x,y,\qquad k,l\in\N,
\end{gather}
corresponding to the factor $G(x+ia)G(y+ia)$, and $b$-dependent poles at
\begin{gather}
z=ib+z_{kl},\qquad z=-ib+2ia+z_{kl},\qquad z=x,y,\qquad k,l\in\N.
\end{gather}

The main disadvantage of the function $\rE(b;x,y)$ compared to the $\cR_r$-function is that it is not even in $x$ and $y$, since the $c$-functions in \eqref{ER} are not even. Instead, it satisf\/ies
\begin{gather}\label{Erefl}
\rE(b;-x,y)=-u(b;x)\rE(b;x,y),
\end{gather}
where
\begin{gather}\label{u}
u(b;z)\equiv -c(b;z)/c(b;-z)=-G(z\pm (ia-ib))/G(z\pm ia).
\end{gather}
On the other hand, $\rE$ inherits all other important properties of $\cR_r$, and is the simplest function to use for Hilbert space purposes. In particular, it has the `unitary asymptotics'
\begin{gather}\label{Eas}
\rE(b;x,y)\sim \exp(i\alpha xy/2)- u(b;-y)\exp(-i\alpha xy/2),\\
 b\in\R,\qquad y\in(0,\infty), \qquad x\to \infty,\nonumber
\end{gather}
cf.~Theorem~1.2 in~II. Here, the $u$-function
encodes the scattering associated with the A$\De$Os $\cA_{\pm}(b;x)$, reinterpreted as commuting self-adjoint operators on the Hilbert space $L^2((0,\infty),dx)$. More precisely, using corresponding results on the $\cE$-function from~II and~III, it follows that the generalized Fourier
transform
\begin{gather}\label{cF}
\cF
  : \ \cC\equiv C_0^{\infty}((0,\infty))\subset \hat{\cH}\equiv L^2((0,\infty),dy) \to\cH
\equiv L^2((0,\infty),dx),
\end{gather}
def\/ined by
\begin{gather}\label{cFE}
(\cF\psi)(x)\equiv\left(\frac{\alpha}{4\pi}\right)^{1/2}\int_0^{\infty}\rE(b;x,y)\psi(y)dy,\qquad
\psi\in\cC,
\end{gather}
extends to a unitary operator, provided the coupling $b$ is suitably restricted. Specif\/ically, it suf\/f\/ices to require
\begin{gather}\label{bres}
b\in [0,a_{+}+a_{-}].
\end{gather}
The self-adjointness of the operators $\hat{\cA}_{\pm}(b)$ on $\cH$ associated to the A$\De$Os $\cA_{\pm}(b;x)$ for $b$ in this interval can then be easily understood from the unitarity of $\cF$: they are the pullbacks to $\cH$ under $\cF$ of the self-adjoint operators of multiplication by $2c_{\pm}(y)$ on $\hat{\cH}$.

As already mentioned, these statements follow from II and III by specialization, but it may help to look f\/irst at Section~9 in the survey~\cite{Rsurv}. Starting from the representation~\eqref{EcE}, the
vector $\gamma^{(1)}$ belongs to the polytope $P$ given by equation~(9.2) in~\cite{Rsurv}, provided $b\in (0,2a)$. Therefore, the transform associated with $\cE(\gamma^{(1)};v,\hat{v})$ (def\/ined by equation~(9.4)) is an isometry. As a consequence, the transform~\eqref{cFE} is an isometry. (The normalization factor in~\eqref{cFE} dif\/fers from that in~equation~(9.4) in~\cite{Rsurv} to accommodate the scale factor 1/2 in the $y$-dependence of~$\cE$ in~\eqref{EcE}.) Isometry of the inverse transform is then clear from the self-duality of $\rE(b;x,y)$. Next, for $b=0$ we have the identity
\begin{gather}\label{ER0}
\rE(0;x,y)=\cR_r(0;x,y)=2\cos(\alpha xy/2),
\end{gather}
cf.~\eqref{ER}--\eqref{cb} and~equation~(7.33) in~\cite{Rsurv} (also, note that for $b=0$ we have $\gamma^{(1)}=\gamma_{\rm f}$, cf.~\eqref{gam1} and~\eqref{gamf}). In view of the symmetry~\eqref{Esym},
it follows that $\cF$ amounts to the cosine transform for $b=0$ and $b=2a$, so these transforms are unitary as well. More generally, we obtain a family of unitary operators
\begin{gather}
\cF(\aaa,b),\qquad (\aaa,b)\in \Pi_u\equiv(0,\infty)^2\times[0,a_{+}+a_{-}],
\end{gather}
which is strongly continuous on the parameter set $\Pi_u$ and satisf\/ies
\begin{gather}\label{cFsym}
\cF(\aaa, b)=\cF(\aaa, a_{+}+a_{-}-b),
\end{gather}
cf.~Theorem~3.3 in III.

The $G$-function asymptotics~\eqref{Gas}
entails that the $c$-function~\eqref{cb} has asymptotics
\begin{gather}\label{cas}
c(b;x)\sim \phi(b)^{\pm 1}\exp(\mp \alpha bx/2),\qquad \re (x)\to\pm\infty,
\end{gather}
with $\phi(b)$ given by~\eqref{phib}. Hence the
$u$-function~\eqref{u} has asymptotics
\begin{gather}\label{uas}
u(b;x)\sim -\phi(b)^{\pm 2},\qquad \re ( x)\to  \pm\infty.
\end{gather}
Also, the ref\/lection equation~\eqref{refl} and the complex conjugation relation~\eqref{Gcon}
entail
\begin{gather}
u(b;-x)u(b;x)=1,\qquad |u(b;x)|=1,\qquad b,x\in\R.
\end{gather}
Thus, if we set
\begin{gather}\label{Fdef}
\rF(b;x,y)\equiv \phi(b)^{-1}(-u(b;x))^{1/2}(-u(b;y))^{1/2}\rE(b;x,y),\qquad b,x,y>0,
\end{gather}
(with the square root phase factors reducing to 1 for $b=0$), then $\rF$ has asymptotics
\begin{gather}\label{Fas}
\rF(b;x,y)\sim [-u(b;y)]^{1/2}\exp(i\alpha xy/2) +[-u(b;y]^{-1/2}\exp(-i\alpha xy/2),\qquad x\to \infty,
\end{gather}
and if we replace $\rE$ by $\rF$ in the above unitary transform~\eqref{cFE}, we retain unitarity.

Introducing the weight function
\begin{gather}\label{w}
w(b;x)\equiv 1/c(b;\pm x)=G(\pm x +ia)/G(\pm x +ia-ib),
\end{gather}
we have
\begin{gather}
w(b;x)>0, \qquad b,x>0,
\end{gather}
and we can also write $\rF$ in terms of $\cR_r$ as
\begin{gather}\label{FR}
\rF(b;x,y) =w(b;x)^{1/2}w(b;y)^{1/2}\cR_r(b;x,y),\qquad b,x,y>0,
\end{gather}
with the positive square roots understood.
Hence, $\rF$ is a joint eigenfunction of the A$\De$Os $H_{\pm}(b;x)$ and
$H_{\pm}(b;y)$
with eigenvalues~\eqref{4ev}, where
\begin{gather} 
H_{\de}(b;z)     \equiv    w(b;z)^{1/2} A_{\de}(b;z)w(b;z)^{-1/2}
\nonumber \\
\label{H} \phantom{H_{\de}(b;z)}{}    =
\sum_{\tau=+,-}\left( \frac{s_{\de}(z-\tau ib)}{s_{\de}(z)}\right)^{1/2}T^z_{\tau ia_{-\de}}\left(\frac{s_{\de}(z+\tau ib)}{s_{\de}(z)}\right)^{1/2}.
\end{gather}

From the $G$-A$\De$Es~\eqref{Gades}
we deduce
\begin{gather}\label{wwr}
w(b;x) =4s_{+}(x)s_{-}(x)w_r(b;x), \qquad w_r(b;x)\equiv G(\pm x -ia+ib).
\end{gather}
For $b\in (0,2a)$ the reduced weight function $w_r(b;x)$ is positive for all real $x$, and since it is also even, its positive square root for $x>0$ has a real-analytic extension to an even positive function on all of~$\R$. By contrast, it is clear from~\eqref{wwr} that $w(b;x)^{1/2}$, $x>0$, extends to an odd real-analytic function on~$\R$. As a consequence, one can also view the transform associated with $\rF(b;x,y)$, $b\in (0,2a)$, as a unitary transform from the odd subspace of $L^2(\R,dy)$ onto the odd subspace of $L^2(\R,dx)$. This is the viewpoint taken in~\cite{Hilb}, where we studied this transform (among other ones) for the special $b$-values $Na_{+}$ with $N\in\N^{*}$. As shown there, for $b>2a$ unitarity and self-adjointness generically break down in a way that can be understood in great detail.

To be sure, the precise connection between the above functions~$\rF$ and~$\cR_r$ and the  functions~$F_r$ and~$E_r$ from~\cite{Hilb} is not clear at face value. But the latter are derived from the functions $M((N+1)a_{+};x,y)$ of~\cite{GLFII}, as specif\/ied below equation~(1.42) in~\cite{Hilb}, so this connection is encoded in the identities~\eqref{cRM} for the special cases $(m,n)=(N+1,0)$, $N\in\N$.

\subsection{The identities (\ref{cRM}) and their consequences}\label{section2.3}

We proceed to prove the general identities~\eqref{cRM}.
Our reasoning involves in particular a comparison of the behavior for $x\to\infty$ of the functions on the l.h.s.\ and r.h.s. For $\cR_r(b;x,y)$ this asymptotics easily follows upon combining~\eqref{ER}, \eqref{Eas} and \eqref{cas}:
\begin{gather}\label{Rras}
\cR_r(b;x,y)\sim \exp(-\alpha bx/2)\sum_{\tau=+,-}c(b;\tau y)\exp(\tau i\alpha xy/2),\\
 b\in\R,\qquad y>0,\qquad x\to\infty.\nonumber
\end{gather}

{\sloppy Next, we consider the functions $M(b_{mn};\pm x,y)$. To begin with, they are eigenfunctions of the four A$\De$Os~\eqref{4A} (where $b=b_{mn}$) with eigenvalues~\eqref{4ev}, cf.~Theorem~II.3 in~\cite{GLFII}. Their `ele\-men\-ta\-ry' form~\eqref{Mform} follows from equations~(3.65)--(3.68) in~\cite{GLFII}. The function $K_{N_{+},N_{-}}(\aaa;x,y)$ occurring in these formulas is specif\/ied in equation~(3.2), with~$S_{N_{\de}}$ given by equation~(2.21). In turn, the coef\/f\/icients in equation~(2.21) are def\/ined via equations~(2.2)--(2.5) in~\cite{GLFII}. (See also Subsection~\ref{section4.1} for more information on these special cases.) It is straightforward to obtain the asymptotics for $\re(x)\to \infty$ from these explicit formulas. Specif\/ically, this yields
\begin{gather}\label{Mas}
M(b_{mn};\pm x,y)=\exp (- \alpha b_{mn}x/2)c(b_{mn};\pm y)e^{ \pm i\alpha xy/2}\big[1+O\big(e^{-\rho \re(x)}\big)\big],\\
 y>0,\qquad \re(x)\to \infty.\nonumber
\end{gather}
The decay rate $\rho$ is the minimum of the two numbers $2\pi /a_{\pm}$, and the implied constant can be chosen uniform for $\im (x)$ varying over $\R$.

}

Comparing \eqref{Rras} and \eqref{Mas}, it follows that the functions on the l.h.s.\ and r.h.s.\ of~\eqref{cRM} have the same asymptotics for $x\to \infty$. It therefore suf\/f\/ices to prove that for f\/ixed $\aaa,y>0$ and $m,n\in\Z$, they must be proportional as functions of $x$. Moreover, we may as well assume $a_{+}/a_{-}$ is irrational, since equality for this case entails equality for all $\aaa>0$. (Indeed, the functions $M(b_{mn};\pm x,y)$ are manifestly real-analytic in $a_{+}$ and $a_{-}$ for $\aaa >0$, and this real-analyticity property is also valid for $\cR_r(b;x,y)$, cf.~I.)

The key consequence of the irrationality assumption is that the vector space of meromorphic joint solutions $f(x)$ to the A$\De$Es
\begin{gather}\label{joint}
A_{\pm}( ma_{+}+na_{-};x)f(x)=2c_{\pm}(y)f(x),\qquad  m,n\in\Z,\qquad y>0,\qquad  a_{+}/a_{-}\notin \Q,
\end{gather}
is two-dimensional.
To explain why this is so, we f\/irst note that the functions $M(b_{mn};\pm x,y)$ are independent solutions to~\eqref{joint}, their independence already being clear from their general form~\eqref{Mform}. Moreover, it follows from their uniform asymptotics~\eqref{Mas} that there exists a~po\-si\-ti\-ve number $\Lambda$, depending on the f\/ixed variables $a_+$, $a_-$, $m$, $n$ and $y$, but not on $x$, such that in the half plane $\re(x)>\Lambda$ both functions are zero-free, and satisfy
\begin{gather}\label{qas}
\lim_{\im(x)\to\infty}M(b_{mn};x,y)/M(b_{mn};-x,y)=0,\qquad \re(x)\in (c_{-},c_{+})\subset [\Lambda, \infty).
\end{gather}

We are now in the position to invoke a result from Section~1 in~\cite{side}, to the ef\/fect that the above  suf\/f\/ices for any joint meromorphic solution $f(x)$ of~\eqref{joint} to be a linear combination of the two functions $M(b_{mn};\pm x,y)$.  Since $\cR_r(b_{mn};x,y)$ is an even meromorphic joint solution, the functions on the l.h.s.\ and r.h.s.\ of~\eqref{cRM} are proportional, so their equality now follows. In particular, for the free case $m=n=0$ we recover the identity \eqref{ER0} from~\eqref{cRM}--\eqref{M00}. Moreover, taking $y=ib_{mn}$ in~\eqref{cRM}, we can invoke~\eqref{cRrnorm} to deduce the corollary
\begin{gather}\label{Mid}
M(b_{mn};x,ib_{mn})+M(b_{mn};-x,ib_{mn})=\frac{G(i(m-1/2)a_{+}+i(n-1/2)a_{-})}{G(i(2m-1/2)a_{+}+i(2n-1/2)a_{-})}.
\end{gather}
Using the $G$-A$\De$Es~\eqref{Gades}, the r.h.s.\ can be rewritten in terms of sine-functions. For the special case $m=N+1$, $n=0$, the resulting identity amounts to equation~(2.78) in~\cite{GLFII}, cf.~also Subsection~\ref{section4.1}.

We would like to add in passing that it is very plausible that~\eqref{qas} is not necessary for two-dimensionality. Indeed, denoting by $\cP_c$ the f\/ield of meromorphic functions with period $c\in\C^{*}$, any third independent joint meromorphic solution would have to be both of the form
\begin{gather}\label{form1}
f(x)=p_1(x)M(b_{mn};x,y)+p_2(x)M(b_{mn};-x,y),\qquad p_1,p_2\in\cP_{ia_{+}},
\end{gather}
and of the form
\begin{gather}\label{form2}
f(x)=q_1(x)M(b_{mn};x,y)+q_2(x)M(b_{mn};-x,y),\qquad  q_1,q_2\in\cP_{ia_{-}}.
\end{gather}
Since the intersection of the f\/ields $\cP_{ia_{+}}$ and $\cP_{ia_{-}}$ reduces
 to the constants when $a_{+}/a_{-}$ is irrational, we expect (but are unable to prove) that this simultaneous representation should lead to a contradiction without appealing to~\eqref{qas}.

Now that we have proved~\eqref{cRM},  it follows that the function $\cR_r(\aaa,b;x,y)$, which is real-analytic on the parameter set
\begin{gather}\label{Pi}
\Pi\equiv \big\{(\aaa,b)\in (0,\infty)^2\times\R\big\},
\end{gather}
 is the continuous interpolation of the functions on the r.h.s.\ of~\eqref{cRM}, which are only def\/ined for the dense subset of `elementary' parameters
\begin{gather}\label{elem}
\Pi_{\rm el}\equiv \{ (\aaa,b)\in\Pi \mid b=ma_{+}+na_{-},\ m,n\in\Z\}.
\end{gather}
The natural question whether another linear combination of $M(b_{mn};\pm x,y)$ that is independent from the even one admits a continuous interpolation as well remains open. In this connection we should point out  that our reasoning at the end of~Section~3 of~\cite{koor} renders this extremely unlikely, but is not conclusive. Indeed, we cannot rule out that the sequence of functions $Q_{-}$ given by equation~(3.15) in~\cite{koor}, with $N_{+}\in\N$, gives rise to an inf\/inity of distinct limits $L_{-}$, corresponding to distinct subsequences. (This oversight is of no consequence for the later sections in~\cite{koor}.) For the same reason, the analogous assertion about the $R$-function, made at the end of~\cite{Rsurv}, has not been completely proved.

Before turning to parameter shifts, we derive  a non-obvious reality feature of $\cR_r$ from the  relations~\eqref{cRM}, namely,
\begin{gather}\label{Rreal}
\overline{\cR_r}(\aaa,b;x,y)=\cR_r(\aaa,b;x,y),\qquad \forall\, (\aaa,b,x,y)\in\Pi\times\R^2.
\end{gather}
The point is that from the explicit formulas for the functions $M$ it is apparent that we have
\begin{gather}
\overline{M}(b_{mn};x,y)=M(b_{mn};-x,y), \qquad x,y\in\R,
\end{gather}
cf.~equation~(3.73) in~\cite{GLFII}. Therefore, \eqref{Rreal} is clear from \eqref{cRM} and interpolation. As a corollary, this yields reality of $\cR$ and $F$ for real parameters and variables, cf.~\eqref{cRr} and~~\eqref{FR}. (Alternatively, this reality property of the $\cR$-function follows from that of the $R$-function proved in Lemma~2.1 of~III.)

\subsection{Parameter shifts}\label{section2.4}

We continue to summarize results concerning parameter shifts from~\cite{anne}, inasmuch as they apply to the present $A_1$ context. In Section~1 of~\cite{anne} we introduced the up-shifts
\begin{gather}\label{Su}
S_{\de}^{(u)}(x)\equiv \frac{-i}{2s_{\de}(x)}\big(T^x_{ia_{-\de}}-T^x_{-ia_{-\de}}\big),
\end{gather}
satisfying
\begin{gather}\label{ush}
S_{\de}^{(u)}(x)A_{\de'}(b;x)=A_{\de'}(b+a_{\de};x)S_{\de}^{(u)}(x),
\end{gather}
and the down-shifts
\begin{gather}\label{Sd}
S_{\de}^{(d)}(b;x)\equiv \frac{2i}{s_{\de}(x)}\big[s_{\de}(x-ib)s_{\de}(x+ia_{-\de}-ib)
T^x_{ia_{-\de}}-(i\to -i)\big],
\end{gather}
satisfying
\begin{gather}\label{dsh}
S_{\de}^{(d)}(b;x)A_{\de'}(b;x)=A_{\de'}(b-a_{\de};x)S_{\de}^{(d)}(b;x),
\end{gather}
where $\de,\de'=+,-$. Clearly, the up-shifts $S_{+}^{(u)}(x)$ and $S_{-}^{(u)}(x)$ commute, and the down-shifts $S_{+}^{(d)}(b_1;x)$ and $S_{-}^{(d)}(b_2;x)$ commute as well. The shifts are also related by
\begin{gather}\label{SA1}
S_{\de}^{(u)}(x)S_{\de}^{(d)}(b;x)=A_{\de}(b;x)^2-4\cos^2(\pi(b-a_{-\de})/a_{\de}),
\\
\label{SA2}
S_{\de}^{(d)}(b+a_{\de};x)S_{\de}^{(u)}(x)=A_{\de}(b;x)^2-4\cos^2(\pi b/a_{\de}),
\end{gather}
where $\de=+,-$. It is a matter of straightforward calculations to verify the formulas \eqref{ush} and~\eqref{dsh}--\eqref{SA2}.

Starting from the joint eigenfunctions $\exp(\pm i\alpha xy/2)$ of $A_{\pm}(0;x)$ with eigenvalues $2c_{\pm}(y)$, one can now obtain joint eigenfunctions with the same eigenvalues for $A_{\pm}(b_{mn};x)$ by acting with the shifts on the plane waves. By construction, these joint eigenfunctions are of the elementary form~\eqref{Mform}. Choosing $a_{+}/a_{-}$ irrational, it follows from two-dimensionality of the joint eigenspace that these eigenfunctions are (generally $y$-dependent) multiples of $M(b_{mn};\pm x;y)$.

A more telling action of the shifts is encoded in
\begin{gather}\label{cRu}
S_{\de}^{(u)}(x)\cR_r(b;x,y)=4s_{\de}(y+ ib)s_{\de}(y-ib)\cR_r(b+a_{\de};x,y),
\\
\label{cRd}
S_{\de}^{(d)}(b;x)\cR_r(b;x,y)=\cR_r(b-a_{\de};x,y).
\end{gather}
Indeed, these relations hold for arbitrary $b$. For $b=b_{mn}$, it then follows by using~\eqref{cRM} that they also hold for the summands $M(b_{mn};\pm x,y)$. (This is because their plane wave factors are independent, cf.~\eqref{Mform}.)

The equations~\eqref{cRd} and~\eqref{cRu} follow from a suitable specialization of equations~(3.11) and~(3.13) in~\cite{anne}. But in the present $A_1$ case we can also derive them quite easily by using the elementary joint eigenfunctions $M(b_{mn};x,y)$ with $a_{+}/a_{-}\notin\Q$. Indeed, once we have shown that~\eqref{cRd},~\eqref{cRu} hold for $y>0$, $b=b_{mn}$,  $m,n\in\Z$, and with $\cR_r$ replaced by $M$, it is easy to deduce~\eqref{cRd},~\eqref{cRu} from~\eqref{cRM} and interpolation. (Note that the four shifts commute with parity.) Their validity for these special cases can be readily verif\/ied: One need only show that the functions on the l.h.s.\ and r.h.s.\ have the same $x\to \infty$ asymptotics, and using~\eqref{Mas} this causes little dif\/f\/iculty.

Next, we obtain the counterparts of the $A_{\de}(b;x)$- and $\cR_r$-shifts for the A$\De$Os $\cA_{\de}(b;x)$ and their joint eigenfunction $\rE(b;x,y)$. (For the $BC_1$ setting we did this in~Section~8 of~\cite{quad}; as in previous cases, it is in fact simpler and more illuminating to obtain the relevant formulas by direct means, instead of by specialization.) They are given by
\begin{gather}\label{sSu}
\cS_{\de}^{(u)}(b;x)\equiv \frac{1}{c(b+a_{\de};x)}S_{\de}^{(u)}(x)c(b;x),
\\
\label{sSd}
\cS_{\de}^{(d)}(b;x)\equiv \frac{1}{c(b-a_{\de};x)}S_{\de}^{(d)}(b;x)c(b;x).
\end{gather}
A moment's thought shows that this entails the validity of \eqref{dsh}--\eqref{SA2} with~$S$, $A$ replaced by~$\cS$,~$\cA$. Also, using the def\/inition~\eqref{cb} of the $c$-function and the $G$-A$\De$Es~\eqref{Gades}, we obtain the explicit formulas
\begin{gather}\label{sSuex}
\cS_{\de}^{(u)}(b;x)=T^x_{ia_{-\de}}-\frac{s_{\de}(x-ib)s_{\de}(x-ib+ia_{-\de})}{s_{\de}(x)s_{\de}(x+ia_{-\de})}T^x_{-ia_{-\de}},
\\
\label{sSdex}
\cS_{\de}^{(d)}(b;x)=T^x_{ia_{-\de}}-\frac{s_{\de}(x+ib)s_{\de}(x+ib-ia_{-\de})}{s_{\de}(x)s_{\de}(x+ia_{-\de})}T^x_{-ia_{-\de}}.
\end{gather}
Notice that they imply
\begin{gather}\label{udrel}
\cS_{\de}^{(u)}(2a-b;x)=\cS_{\de}^{(d)}(b;x).
\end{gather}
Finally,  a straightforward calculation yields the counterparts of~\eqref{cRu} and~\eqref{cRd}:
\begin{gather}
\cS_{\de}^{(u)}(b;x)\rE(b;x,y)=2e_{\de}(-i\pi b)s_{\de}(y+ib)
\rE(b+a_{\de};x,y),
\\
\cS_{\de}^{(d)}(b;x)\rE(b;x,y)=2e_{\de}(i\pi (b-a_{-\de}))s_{\de}(y-ib+ia_{-\de})\rE(b-a_{\de};x,y).
\end{gather}

To conclude this section, we point out that the eight shifts acting on $x$ have duals acting on $y$ given by the formulas~\eqref{Su}, \eqref{Sd} and~\eqref{sSu}--\eqref{sSdex} with $x\to y$. By self-duality, their respective actions on  $\cR_r(b;x,y)$ and $\rE(b;x,y)$ follow from the above by interchanging~$x$ and~$y$.

\section[Five minimal representations of the $\cR$-function]{Five minimal representations of the $\boldsymbol{\cR}$-function}\label{section3}

We begin this section by focusing on the kernel function
\begin{gather}
K(b;x,v):=G((\pm x \pm v -ib)/2)\equiv\prod_{\de_1,\de_2=+,-}G((\de_1 x+\de_2 v-ib)/2).
\end{gather}
We have established that this function satisf\/ies three independent kernel identities. We expect that these might be useful in other contexts than the present one. Indeed, here we only need the special case~\eqref{dv} of the f\/irst of the identities. We  collect the three identities in the following proposition.

\begin{prop}\label{proposition3.1}
Letting $b,d\in\C$, we have
\begin{gather}
\frac{s_{\de}(x-ib+id)}{s_{\de}(x)}K(b;x-ia_{-\de},v)+(i\to -i)\nonumber\\
\qquad{} =\frac{s_{\de}(v-id)}{s_{\de}(v)}K(b;x,v-ia_{-\de})+(i\to -i),\label{idd}
\\
\frac{s_{\de}(x-ib)}{s_{\de}(x)}K(b;x-2ia_{-\de},v)+(i\to -i)\nonumber\\
\qquad{}=
\frac{s_{\de}(v-ib)}{s_{\de}(v)}K(b;x,v-2ia_{-\de})+(i\to -i),\label{id2}
\\
\frac{s_{\de}((x-ib)/2)}{s_{\de}(x/2)}K(b;x-ia_{-\de},v)+(i\to -i)\nonumber\\
\qquad{} =
\frac{s_{\de}((v-ib)/2)}{s_{\de}(v/2)}K(b;x,v-ia_{-\de})+(i\to -i).\label{id3}
\end{gather}
\end{prop}

\begin{proof}
To prove \eqref{idd}, we divide l.h.s.\ and r.h.s.\ by $K(b;x-ia_{-\de},v)$ and use the A$\De$Es~\eqref{Gades} to write the result as
\begin{gather}
  \frac{s_{\de}(x-ib+id)}{s_{\de}(x)}+\frac{s_{\de}(x+ib-id)}{s_{\de}(x)}\frac{c_{\de}((x+ v-ib)/2)}{c_{\de}((x+ v+ib)/2)}\frac{c_{\de}((x- v-ib)/2)}{c_{\de}((x- v+ib)/2)}  \nonumber \\
\qquad{} = \frac{s_{\de}(v-id)}{s_{\de}(v)} \frac{c_{\de}((x-v-ib)/2)}{c_{\de}((x-v+ib)/2)} +   \frac{s_{\de}(v+id)}{s_{\de}(v)} \frac{c_{\de}((x+v-ib)/2)}{c_{\de}((x+v+ib)/2)}.\label{idq}
\end{gather}
Both sides are $2ia_{\de}$-periodic functions of $x$ with equal limits
\begin{gather}
e_{\de}(\pm (-ib+id))+e_{\de}(\pm (-ib-id)),\qquad \re x \to\pm\infty.
\end{gather}
The residues at the (generically simple) poles $x=0$, $x=ia_{\de}$ in the period strip clearly cancel. By Liouville's theorem, it remains to check that the residues at the poles $x=\pm v -ib/2\pm ia_{\de}$ cancel as well, and this is a routine calculation.

Next, we divide~\eqref{id2} by $K(b;x,v)$ and use~\eqref{Gades} to obtain
\begin{gather}
  \frac{s_{\de}(x-ib)}{s_{\de}(x)}\frac{c_{\de}((x-ia_{-\de}/2\pm v +ib)/2)}{c_{\de}((x-ia_{-\de}/2\pm v -ib)/2)}+ (x\to -x)
\nonumber \\
\qquad{} = \frac{s_{\de}(v-ib)}{s_{\de}(v)} \frac{c_{\de}((v-ia_{-\de}/2\pm x +ib)/2)}{c_{\de}((v-ia_{-\de}/2\pm x -ib)/2)}+ (v\to -v).
\end{gather}
Both sides are $2ia_{\de}$-periodic functions of $x$ with equal limits
\begin{gather}
e_{\de}(\pm ib)+e_{\de}(\mp ib)=\frac{s_{\de}(v-ib)}{s_{\de}(v)} +\frac{s_{\de}(v+ib)}{s_{\de}(v)},\qquad \re x \to\pm\infty.
\end{gather}
The residues at $x=0$ and $x=ia_{\de}$ manifestly cancel. It is a straightforward calculation to verify that the residues at the remaining poles $x=\pm ia_{\de}\pm ia_{-\de}/2\pm v +ib$ cancel, too. Hence~\eqref{id2} follows.

Finally, to prove~\eqref{id3}, we divide both sides by
$K(b;x-ia_{-\de},v)$ and use~\eqref{Gades} to get as the counterpart of~\eqref{idq}:
\begin{gather}
   \frac{s_{\de}((x-ib)/2)}{s_{\de}(x/2)}+\frac{s_{\de}((x+ib)/2)}{s_{\de}(x/2)}\frac{c_{\de}((x+ v-ib)/2)}{c_{\de}((x+ v+ib)/2)}\frac{c_{\de}((x- v-ib)/2)}{c_{\de}((x- v+ib)/2)}  \nonumber \\
\qquad{} = \frac{s_{\de}((v-ib)/2)}{s_{\de}(v/2)} \frac{c_{\de}((x-v-ib)/2)}{c_{\de}((x-v+ib)/2)} +   \frac{s_{\de}((v+ib)/2)}{s_{\de}(v/2)} \frac{c_{\de}((x+v-ib)/2)}{c_{\de}((x+v+ib)/2)}.
\end{gather}
Both sides are $2ia_{\de}$-periodic functions of $x$ with equal limits
\begin{gather}
e_{\de}(\pm (-ib/2))+e_{\de}(\pm (-3ib/2)),\qquad \re x \to\pm\infty.
\end{gather}
As before, residue cancellation at $x=0$ and $x=ia_{\de}$ is immediate, whereas the verif\/ication that the residues at the remaining poles $x=\pm ia_{\de} \pm v  -ib/2$ cancel as well involves a bit more work.
\end{proof}

From~\eqref{Apm} we see that the identity~\eqref{idd} can be rewritten as
\begin{gather}
A_{\de}(b-d;x)K(b;x,v)=A_{\de}(d;v)K(b;x,v),\qquad \de=+,-.
\end{gather}
At f\/irst sight, one might think that the identities~\eqref{id2} and~\eqref{id3} can also be rewritten by using a rescaled version of the two commuting $A_1$ dif\/ference operators $A_{\pm}(b;z)$. The two dif\/ference operators
\begin{gather}
\frac{s_{\de}((z-ib)/2)}{s_{\de}(z/2)}T^z_{ia_{-\de}} +
\frac{s_{\de}((z+ib)/2)}{s_{\de}(z/2)}T^z_{-ia_{-\de}},\qquad \de=+,-,
\end{gather}
featuring in~\eqref{id3}, do not even commute, however. Thus no such rescaling is possible for~\eqref{id3}. The dif\/ference operators
\begin{gather}
\frac{s_{\de}(z-ib)}{s_{\de}(z)}T^z_{2ia_{-\de}} +
\frac{s_{\de}(z+ib)}{s_{\de}(z)}T^z_{-2ia_{-\de}},\qquad \de=+,-,
\end{gather}
corresponding to~\eqref{id2} do commute. Even so, one can only rescale one of the operators such that it takes the~$A_1$ form~\eqref{Apm}, but not both at once.

For the purpose of studying the $A_1$ operators, then, we can only make use of~\eqref{idd}. More specif\/ically, our starting point is the special case $d=0$:
\begin{gather}\label{dv}
A_{\de}(b;x)K(b;x,v)= \big(T^v_{ia_{-\de}} +T^v_{-ia_{-\de}} \big)K(b;x,v).
\end{gather}
In order to exploit this identity, we introduce  the Fourier transform
\begin{gather}\label{defB}
B(b;x,y)\equiv \frac{1}{2}\int_{\R}dvK(b;x,v)\exp(i\alpha vy/2),\qquad b\in(0,2a),\qquad x,y\in \R.
\end{gather}
The integral is well def\/ined, since the $b$-restriction ensures that the $v$-poles of the integrand at
\begin{gather}\label{Kpo}
\pm v=x+2ia-ib+z_{kl}, \qquad \pm v=-x+2ia-ib+z_{kl},\qquad k,l\in\N,
\end{gather}
(cf.~\eqref{Gpo}--\eqref{zkl}), stay away from the real axis, and since the $G$-asymptotics~\eqref{Gas} entails an exponential decay
\begin{gather}
K(b;x,v)\sim \exp(\mp \alpha bv/2),\qquad b>0,\qquad x\in\C,\qquad \re v\to\pm\infty.
\end{gather}
These features also imply that $B$ extends from the real $x$-axis to a function that is holomorphic in the strip $\im x\in (-2a+b,2a-b)$.

Next, we temporarily  assume
\begin{gather}
b\in(0,a_s/2),\qquad a_s\equiv\min(\aaa).
\end{gather}
Then the action of the shifts in the A$\De$Os $A_{\pm}(b;x)$ given by~\eqref{Apm} is well def\/ined on $B(b;x,y)$, provided we restrict~$x$ to a strip $|\im x|<a_s/2$. Moreover, we may take the shifts under the integral sign in~\eqref{defB} and use the kernel identity~\eqref{dv} to obtain
\begin{gather}
A_{\de}(b;x)B(b;x,y)=\frac{1}{2}\int_{\R}dv\sum_{\tau=+,-}K(b;x,v+\tau ia_{-\de})\exp(i\pi vy/a_+a_-),\\
 |\im x|<a_s/2.\nonumber
\end{gather}
Upon shifting contours $\R\to \R\pm ia_{-\de}$, no poles are met, and so we obtain
\begin{gather}
\frac12\sum_{\tau=+,-}e_{\de}(\tau y)\int_{\R-i\tau a_{-\de}}dvK(b;x,v)\exp(i\pi vy/a_+a_-).
\end{gather}
The integrands of both terms are now equal, and the contours can be shifted back to $\R$ without changing the value of the integrals. Hence we deduce the eigenvalue equations
\begin{gather}\label{eigen}
A_{\pm}(b;x)B(b;x,y)=2c_{\pm}(y)B(b;x,y),\qquad |\im x|<a_s/2.
\end{gather}

Reverting to our previous assumption $b\in(0,2a)$, we proceed to obtain two dif\/ferent representations of $B(b;x,y)$. To this end we use the Plancherel relation
\begin{gather}\label{Plan}
\int_{\R}dpf(p)g(p)=\frac{\alpha}{2\pi}\int_{\R}dq\hat{f}(q)\hat{g}(-q),\qquad f,g\in L^2(\R)\cap L^1(\R),
\end{gather}
with the Fourier transform def\/ined by
\begin{gather}\label{Four}
\hat{h}(q)=\int_{\R}dp\exp(i\alpha pq)h(p),\qquad h=f,g.
\end{gather}
Rewriting $B$ as
\begin{gather}\label{Bnew}
B(b;x,y)=\int_{\R}dp\frac{G(p+(x-ib)/2)G(p-(x+ib)/2)\exp(i\alpha py)}{G(p-(x-ib)/2)G(p+(x+ib)/2)},\\
 b\in(0,2a),\qquad x,y\in\R,\nonumber
\end{gather}
we now def\/ine
\begin{gather}
f_1(p)\equiv\frac{G(p+(x-ib)/2)}{G(p-(x-ib)/2)},\qquad f_2(p)\equiv\frac{G(p+(x-ib)/2)}{G(p+(x+ib)/2)},
\\
g_1(p)\equiv\frac{G(p-(x+ib)/2)\exp(i\alpha py)}{G(p+(x+ib)/2)},\qquad
g_2(p)\equiv\frac{G(p-(x+ib)/2)\exp(i\alpha py)}{G(p-(x-ib)/2)}.
\end{gather}
We can calculate the Fourier transforms of these four functions by using Proposition~\ref{propositionC.1}. Doing so, we use the Plancherel relation~\eqref{Plan} and then replace $q$ by $z+y/2$ to obtain the two representations announced above:
\begin{gather}\label{rep2}
B(b;x,y)=G(\pm x +ia-ib)\int_{\R}dz \frac{G(z\pm (x-y)/2-ia+ib/2)}{G(z\pm (x+y)/2+ia-ib/2)},
\\
\label{rep3}
B(b;x,y)=G(ia-ib)^2\int_{\R}dz G(\pm z\pm y/2-ia+ib/2)\exp(i\alpha zx).
\end{gather}

Next, we compare~\eqref{rep3} to~\eqref{Bnew}, deducing that it can be rewritten as
\begin{gather}\label{BB}
B(b;x,y)=G(ia-ib)^2B(2a-b;y,x).
\end{gather}
Also, def\/ining a new function $C(b;x,y)$ by~\eqref{defC},
we see that~\eqref{rep2} amounts to
\begin{gather}\label{BC}
B(b;x,y)=\sqrt{\frac{2\pi}{\alpha}}G(\pm x +ia-ib)C(2a-b;x,y).
\end{gather}

The function $C(\aaa,b;x,y)$ is of central importance for what follows. Writing it as
\begin{gather}
C(\aaa,b;x,y)= \sqrt{\frac{\alpha}{2\pi}}\int_{\R}dz \prod_{\de=+,-}\frac{G(\aaa;z-u_{\de})}{G(\aaa;z-d_{\de})},
\end{gather}
where we have introduced
\begin{gather}
u_{\pm}\equiv \pm (x-y)/2+ib/2,\qquad d_{\pm}\equiv \pm (x+y)/2-ib/2,
\end{gather}
we infer that its behavior under analytic continuation in its 5 variables is immediate from the general analysis in Appendix~B of~I (with $N$ specialized to 2). We proceed to summarize the salient information.

To this end, we need the function $E(\aaa;z)$ discussed in~Appendix~\ref{appendixA}, cf.~\eqref{GE}. Introducing
\begin{gather}\label{P}
P(\aaa,b;x,y)\equiv C(\aaa,b;x,y)E(\aaa;\pm x+ib-ia)E(\aaa;\pm y+ib-ia),\!\!\!\!
\end{gather}
the product function $P(\aaa,b;x,y)$ extends from $(0,\infty)^2\times (0,a_{+}+a_-)\times \R^2$ to a function that is holomorphic in the domain
\begin{gather}
D(\aaa,b)\equiv \big\{ (\aaa,b,x,y)\in \C^5 \mid \re a_+>0,\ \re a_->0, \ \re (b/a_+a_-)>0\big\}.
\end{gather}
 Hence $C$ is meromorphic in $D(\aaa,b)$, with poles occurring solely at the zeros
\begin{gather}
\pm x=2ia-ib+z_{kl},\qquad \pm y=2ia-ib+z_{kl},\qquad k,l\in\N,
\end{gather}
of the $E$-product, cf.~\eqref{Eze}; moreover, the maximal multiplicity of a pole at $z=z_0$, with $z=x,y$, is given by the zero multiplicity at $z=z_0$ of the pertinent $E$-factor.

The corresponding meromorphy properties of $B$ are now clear from its relation~\eqref{BC} to $C$: It continues meromorphically to the domain $D(\aaa,2a-b)$. From~\eqref{BB} we then deduce that $B(\aaa,b;x,y)$ has a meromorphic extension to the larger domain $D_+$~\eqref{Dp}.
Using~\eqref{BC} again, it now follows that $C$ has a meromorphic extension to $D_+$ as well.

We proceed to obtain further information on the function $C$. First, we list features that are immediate from its def\/inition~\eqref{defC} and properties of the $G$-function, cf.~Appendix~\ref{appendixA}:
\begin{gather}\label{Cmodinv}
C(a_{+},a_{-},b;x,y)=C(a_{-},a_{+},b;x,y),\qquad  ({\rm modular\  invariance}),
\\
\label{Csd}
C(a_{+},a_{-},b;x,y)=C(a_{+},a_{-},b;y,x),\qquad  (\mbox{self-duality}),
\\
\label{Ceven}
C(b;x,y)=C(b;\de x, \de'y),\qquad \de,\de'=+,-,\qquad ({\rm evenness}),
\\
\label{Cscale}
C(\aaa,b;x,y)=C(\lambda  a_{+},\lambda a_{-}, \lambda b;\lambda x,\lambda y),\qquad \lambda>0,\qquad ({\rm scale\ invariance}),
\\
\label{Creal}
C(\aaa,b;x,y)\in\R,\!\!\!\!\qquad  \aaa>0,\!\!\!\!\qquad b\in(0,2a),\!\!\!\!\qquad x,y\in\R,\!\!\!\!\qquad \mbox{(real-valuedness)}.\!\!\!\!
\end{gather}
Clearly, the relations \eqref{Cmodinv}--\eqref{Cscale} are well def\/ined and hold true on $D_+$~\eqref{Dp}.  Also, the property~\eqref{Creal} can be rendered manifest by substituting the integral representation~\eqref{ghyp} in the four $G$-functions and combining factors to obtain
\begin{gather}
C(\aaa,b;x,y)=  \sqrt{\frac{\alpha}{2\pi}}  \int_{\R}dz\cos\left( \int_{\R}\frac{dw}{w}\frac{\sin(xw)\sin(yw)\cosh(bw)\sin(2zw)}{\sinh(a_{+}w)\sinh(a_-w)}\right)  \\
\phantom{C(\aaa,b;x,y)=}{}  \times \exp\left( \int_{\R}\frac{dw}{w}\left(\frac{\cos(xw)\cos(yw)\sinh(bw)\cos(2zw)}{\sinh(a_{+}w)\sinh(a_-w)}-\frac{b}{a_+a_-w}\right)\right),  \nonumber
\end{gather}
  where $(\aaa,b,x,y)\in(0,\infty)^2\times(0,2a)\times\R^2$.

Second, combining~\eqref{Csd} with~\eqref{BB} and~\eqref{BC}, we deduce
\begin{gather}\label{Cbsym}
\frac{G(ia-ib)C(b;x,y)}{G(x+ia-ib)G(y+ia-ib)}=
\frac{G(ib-ia)C(2a-b;x,y)}{G(x-ia+ib)G(y-ia+ib)},\qquad (b\mbox{-symmetry}).
\end{gather}
Third, we can use Proposition~\ref{propositionC.1} once more to obtain from~\eqref{defC} the explicit result
\begin{gather}\label{Cnorm}
C(b;x,ib)=G(ia-2ib)G(ib-ia)^2,\qquad ({\rm normalization}),
\end{gather}
where $(\aaa,b,x)\in(0,\infty)^2\times (0,a)\times \R$.
Last but not least, we claim that we have the joint eigenvalue equations
\begin{gather}\label{Cades}
A_{\pm}(b;x)C(b;x,y)=2c_{\pm}(y)C(b;x,y),\qquad A_{\pm}(b;y)C(b;x,y)=2c_{\pm}( x)C(b;x,y).
\end{gather}

 To prove this claim, we f\/irst note that the eigenvalue equations~\eqref{eigen} continue meromorphically to $D_+$. Next, we observe that the $G$-A$\De$Es~\eqref{Gades} imply
\begin{gather}
G(\pm x +ia-ib)^{-1}A_{\de}(b;x)G(\pm x +ia-ib)=A_{\de}(2a-b;x),\qquad \de=+,-.
\end{gather}
(Here, we view the l.h.s.\ as the product of three operators acting on meromorphic functions.)
Hence we obtain via~\eqref{BC}
\begin{gather}
A_{\pm}(2a-b;x)C(2a-b;x,y)=2c_{\pm}(y)C(2a-b;x,y),
\end{gather}
which is equivalent to the f\/irst two A$\De$Es in~\eqref{Cades}. The last two are then clear from the self-duality relation~\eqref{Csd}.

All of the properties of $C$ just derived also hold true for the function $G(ib-ia)\cR_r$ def\/ined by~\eqref{cRr} and~\eqref{cRdef}, cf.~Section~\ref{section2}. By using solely the eigenvalue properties~\eqref{Cades}, the evenness properties~\eqref{Ceven}, and the normalization~\eqref{Cnorm}, we can now show that these two functions coincide, as announced in the Introduction, cf.~\eqref{cRC}. Specif\/ically, applying the uniqueness argument in Subsection~\ref{section2.3} to $C$ in its dependence on $x$, we obtain the equality~\eqref{cRC} up to a~proportionality factor $p(\aaa,b,y)$. Repeating this argument for the $y$-dependence, we see that the proportionality factor can only depend on the parameters $a_+$, $a_-$, $b$. From the normalization relation~\eqref{Cnorm} it then follows that $p=1$, thus proving~\eqref{cRC}.

Let us now collect the resulting minimal representations of the $\cR$-function. From~\eqref{cRC} and~\eqref{cRr} we obtain
\begin{gather}\label{Ri}
\cR(b;x,y)=\frac{1}{\sqrt{a_+a_-}}\frac{G(2ib-ia)}{G(ib-ia)^2}
\int_{\R}dz \frac{G(z\pm (x-y)/2-ib/2)}{G(z\pm (x+y)/2+ib/2)}.
\end{gather}
Next, combining~\eqref{Bnew} and~\eqref{BC}, we deduce
\begin{gather}
\cR(b;x,y)=\frac{1}{\sqrt{a_+a_-}}\frac{G(2ib-ia)}{G(ib-ia)^2}
G(\pm x +ia-ib)\nonumber\\
\phantom{\cR(b;x,y)=}{} \times
\int_{\R}dz\frac{G(z\pm x/2-ia+ib/2)}{G(z\pm x/2+ia-ib/2)}\exp(i\alpha zy),\label{Rii}
\end{gather}
and using~\eqref{BB} we infer
\begin{gather}\label{Riii}
\cR(b;x,y)=\frac{1}{\sqrt{a_+a_-}}G(2ib-ia)
G(\pm x +ia-ib)
\int_{\R}dz\frac{G(z\pm y/2-ib/2)}{G(z\pm y/2+ib/2)}\exp(i\alpha zx).\!\!\!
\end{gather}
Finally, using the self-duality property of $\cR$, we obtain from~\eqref{Rii} and~\eqref{Riii} the representations
\begin{gather}
\cR(b;x,y)=\frac{1}{\sqrt{a_+a_-}}\frac{G(2ib-ia)}{G(ib-ia)^2}
G(\pm y +ia-ib)\nonumber\\
\phantom{\cR(b;x,y)=}{}\times
\int_{\R}dz\frac{G(z\pm y/2-ia+ib/2)}{G(z\pm y/2+ia-ib/2)}\exp(i\alpha zx),\label{Riv}
\\
\cR(b;x,y)=\frac{1}{\sqrt{a_+a_-}}G(2ib-ia)
G(\pm y +ia-ib)\nonumber\\
\phantom{\cR(b;x,y)=}{}\times
\int_{\R}dz\frac{G(z\pm x/2-ib/2)}{G(z\pm x/2+ib/2)}\exp(i\alpha zy).\label{Rv}
\end{gather}
The f\/ive representations~\eqref{Ri}--\eqref{Rv} are well def\/ined and hold true for $(\aaa,b,x,y)\in(0,\infty)^2\times(0,2a)\times\R^2$. (Combining~\eqref{Ri} with~\eqref{Cbsym}, we get a further representation that we do not consider.)

Taking stock of the above developments, we note that we might have started from the f\/irst minimal representation~\eqref{Ri} to {\em define} the $\cR$-function. Then many of its properties follow quite easily. On the other hand, it seems not feasible to give a direct proof of its crucial joint eigenfunction property. With hindsight, however, this can be shown by f\/irst obtaining the second representation~\eqref{Rii} (say) via Proposition~\ref{propositionC.1}, and then using the identity~\eqref{dv} to arrive at~\eqref{eigen}. From this the joint eigenfunction property~\eqref{Cades} follows as before.

Another important property of $\cR$ is its asymptotic behavior for $x\to \infty$. Like other features addressed in this section, this can already be gleaned from Section~\ref{section2}, via the specialization of the more general asymptotics  of the function $\cE(\gamma;v,\hat{v})$ obtained in~Theorem~1.2 of~II. However, provided we restrict $b$ to the interval $(0,2a)$, it is quite easy to obtain the $x\to \infty$ asymptotics directly from the new representations of $\cR$ in terms of a Fourier transform.

To detail this, let us f\/irst note that we need only consider the function E$(b;x,y)$, which we can now view as being def\/ined via~\eqref{ER}--\eqref{cb}. (Indeed, there is no dif\/f\/iculty in obtai\-ning the asymptotics of the $c$-function; in this connection, compare~\eqref{cas}, \eqref{Eas} and~\eqref{Rras}.) Using~\eqref{Riii}, we deduce that for $b\in(0,2a)$ the E-function has the representation
\begin{gather}\label{Enew}
\rE(b;x,y)=\frac{\phi(b)}{\sqrt{a_+a_-}}\frac{G(ib-ia)}{c(b;y)}\frac{G(x+ia)}{G(x-ia+ib)}
\int_{\R}dz\frac{G(z\pm y/2-ib/2)}{G(z\pm y/2+ib/2)}\exp(i\alpha zx).
\end{gather}
Letting $y\in(0,\infty)$, we can shift the contour up by $a-b/2 +\epsilon$, where $\epsilon>0$ is small enough so that only the simple poles at
\begin{gather}
z=\pm y/2 -ib/2+ia,
\end{gather}
are encountered. The residues at these poles easily follow from~\eqref{Gres}, yielding a contribution
\begin{gather}\label{ressum}
M(b;x)(\exp(i\alpha xy/2)+c(b;-y)\exp(-i\alpha xy/2)/c(b;y)),
\end{gather}
with the multiplier given by
\begin{gather}
M(b;x)\equiv \phi(b) G(x+ia)G(-x +ia-ib)\exp(-\alpha x (a-b/2)).
\end{gather}
Now from~\eqref{Gas} we see that $M(b;x)$ converges to 1 for $x\to \infty$. To recover the asymptotics~\eqref{Eas}, therefore, it is enough to show that the r.h.s.\ of~\eqref{Enew} with $z$ replaced by $z+ia-ib/2+i\epsilon$ vanishes for $x\to \infty$.

To prove this, we write the shifted contour integral as
\begin{gather}\label{sci}
\frac{1}{\sqrt{a_+a_-}}\frac{G(ib-ia)}{c(b;y)}M(b;x)\exp(-\epsilon\alpha x)
\int_{\R}dz\frac{G(z+i\epsilon \pm y/2-ib+ia)}{G(z+i\epsilon \pm y/2+ia)}\exp(i\alpha zx).
\end{gather}
Now one need only use~\eqref{Gas} to verify that the integrand is bounded by a multiple of \linebreak $\exp(-\alpha b|z|/2)$, which implies the function~\eqref{sci} does converge to 0 for $x\to\infty$.

We stress that this short argument only yields~\eqref{Eas} under the restriction $b\in(0,2a)$.  In particular, by contrast to the previous contour integral representation used in~II, one must cope with an inevitable contour pinching when one tries to use~\eqref{Enew} to go beyond this $b$-interval.

Another issue is that stronger asymptotic estimates than just obtained are necessary to recover the Hilbert space transform features for the E-function sketched in Subsection~\ref{section2.3}, cf.~\eqref{cF}--\eqref{cFsym}. It is beyond our scope to study this further, but we would like to repeat that the $b$-interval $[0,2a]$ cannot be enlarged without losing the critical unitarity and self-adjointness properties~\cite{Hilb}.

At face value, the new representations~\eqref{Ri}--\eqref{Rv} seem to hold promise for a direct proof of the shift properties of $\cR_r$, cf.~\eqref{cRu}--\eqref{cRd}. Even so, we were unable to push this through.  To date, therefore, the only reasoning yielding the properties for general $b$ is to f\/irst derive them for special $b$-values, as sketched in~Subsection~\ref{section2.4}. We now turn to a study of the $\cR_r$-function for these special values.

\section{Specializations and nonrelativistic limit}\label{section4}

\subsection{Elementary special cases}\label{section4.1}

As already mentioned, the functions
\begin{gather}\label{RN}
R_N(\aaa;x,y)\equiv \cR_r(\aaa,(N+1)a_+;x,y),\qquad N\in\N,
\end{gather}
have been extensively studied before. They were f\/irst obtained more than twenty years ago~\cite{fdss}, and then reconsidered from an  algebraic and function-theoretic viewpoint in~\cite{GLFII} and from a~representation-theoretic viewpoint in a paper by van Diejen and Kirillov~\cite{diki}. The corresponding Hilbert space transforms were studied in great detail in~\cite{Hilb}.

Our f\/irst goal in this section is to demonstrate how the elementary character of these functions can be directly understood from the Fourier transform representations~\eqref{Rii}--\eqref{Rv}. Indeed, thus far the relation encoded in~\eqref{cRM} has only been shown by appealing to a uniqueness argument. The crux is that for the choices $b=(N+1)a_{\pm}$  one can use the $G$-A$\De$Es~\eqref{Gades} to obtain integrals that can be explicitly evaluated by a residue calculation.

Specif\/ically, let us start from~\eqref{Riii} to obtain f\/irst
\begin{gather}
R_N(x,y)    =   \frac{4^{-N-1}}{\sqrt{a_+a_-}}G(i(N+1)a_+-ia)
G(\pm x +ia-i(N+1)a_+)
\nonumber \\
\phantom{R_N(x,y)    =}{} \times
\int_{\R}dze^{i\alpha zx}/\prod_{j=0}^Nc_-(z\pm y/2-i(N-2j)a_+/2).\label{RNrep}
\end{gather}
Now we recall that~\eqref{Riii} is valid for $b\in(0,2a)$, which implies we have $Na_+<a_-$ in~\eqref{RNrep}. Taking $y>0$ from now on, it follows that the integrand has $2N+2$ simple poles in the strip $\im z\in(0,a_-)$. The product yields a function that is $ia_-$-periodic in $z$. Thus, denoting the integral by~$I_N$, we have
\begin{gather}
I_N-\exp(-2\pi x/a_+)I_N=2\pi i\sum_{\de=+,-}\sum_{j=0}^N {\rm Res}\, \left(z=\frac{1}{2}\big(\de y +ia_-+i(N-2j)a_+\big)\right).
\end{gather}
The residues are easily calculated, and hence we obtain
\begin{gather}\label{IN}
I_N=\frac{(-i)^{N+1}a_-}{s_+(x)}\!\left( \frac{e^{i\alpha yx/2}}{s_-(y)}\!
\sum_{j=0}^N\frac{e_-(-(N{-}2j)x)}{\prod\limits_{k\ne j}s_-(y{+}i(k{-}j)a_+)\sin(\pi(j-k)a_+/a_-)}{+}(y{\to} -y)\!\right).\!\!\!
\end{gather}
Also, the prefactor can be calculated by using~\eqref{Gades} once more, combined with~\eqref{Geval}. Introducing
\begin{gather}\label{PN}
P_N(z)\equiv \prod_{j=-N}^N 2s_-(z+ija_+),\qquad z=x,y,
\end{gather}
we get
\begin{gather}\label{pref}
\frac{2^{-2N-1}}{a_-}\frac{s_+(x)}{P_N(x)}\prod_{j=1}^N 2\sin(\pi ja_+/a_-).
\end{gather}

For $N=0$ the product of~\eqref{IN} and~\eqref{pref} yields
\begin{gather}\label{Rzero}
R_0(x,y)=\frac{\sin(\pi xy/a_+a_-)}{2s_-(x)s_-(y)}.
\end{gather}
More generally, the product can be written as
\begin{gather}\label{RNfor}
R_N(x,y)=(-i)^{N+1}(K_N(x,y)-K_N(x,-y))/P_N(x)P_N(y),
\end{gather}
where we have set
\begin{gather}
  K_N(x,y)  \equiv  \exp(ixy/a_+a_-)\prod_{l=1}^N 2\sin(\pi la_+/a_-)
\nonumber \\
\phantom{K_N(x,y)  \equiv}{}   \times \sum_{j=0}^N\frac{e_-((2j-N)x)\prod\limits_{k=j+1}^Ns_-(y-ika_+)\prod\limits_{k=N-j+1}^Ns_-(y+ika_+)}{\prod\limits_{k\ne j}\sin(\pi(j-k)a_+/a_-)}.\label{KN}
\end{gather}
Introducing the phase factor
\begin{gather}
q\equiv \exp(i\pi a_+/a_-),
\end{gather}
it is not hard to see from~\eqref{KN} that $K_N$ is of the form{\samepage
\begin{gather}\label{KNform}
K_N(x,y)  = \exp(ixy/a_+a_-)e_-(Nx+Ny)S_N(q;e_-(-2x),e_-(-2y)),
\\
\label{SN}
S_N(q;r,t)=\sum_{k,l=0}^Nc_{kl}^{(N)}(q)r^kt^l,
\end{gather}
with $c_{kl}^{(N)}(q)$ a rational function of $q$.}

Thus far, our conclusions about $R_N$ and $K_N$ were only based on an explicit evaluation of the representation~\eqref{Riii}. However, a lot more information follows upon using the features of $\cR_r(b;x,y)$. In particular, the function $K_N(x,y)/P_N(x)P_N(y)$ is a joint eigenfunction of the four A$\De$Os~\eqref{4A} (where $b=(N+1)a_+$) with eigenvalues~\eqref{4ev}, since this holds true for $R_N(x,y)$. Also, the self-duality and evenness of $R_N$ imply
\begin{gather}
K_N(x,y)=K_N(y,x)=K_N(-x,-y),
\end{gather}
and this entails that the coef\/f\/icients in~\eqref{SN} have the symmetry properties
\begin{gather}\label{csym}
c_{kl}^{(N)}=c_{lk}^{(N)}=c_{N-k,N-l}^{(N)},\qquad k,l=0,\ldots,N.
\end{gather}
Of course, this can be easily checked for small $N$, but for arbitrary integers~\eqref{csym} is not at all obvious from~\eqref{KN}.

One more feature of the coef\/f\/icients is that they are Laurent polynomials in $q$ with integer coef\/f\/icients. Like the symmetries~\eqref{csym}, it is not a routine matter to show this from~\eqref{KN}. The crux is, however, that the above functions $K_N$ coincide with those of~\cite{GLFII}, by virtue of the uniqueness argument used in~Subsection~\ref{section2.3}. The coef\/f\/icients were studied in detail in Section~II of~\cite{GLFII}, and there the interested reader can f\/ind explicit formulas for the coef\/f\/icients as Laurent polynomials in $q$. See also the paper by van Diejen and Kirillov~\cite{diki}, where yet dif\/ferent representations of the functions $K_N(x,y)$ were obtained.

With a little more ef\/fort, the elementary character of $\cR_r$ for the more general $b$-values
\begin{gather}\label{bNM}
b_{+-}\equiv (N+1)a_{+}-Ma_-,\qquad N,M\in\N,
\end{gather}
can also be understood from~\eqref{Riii}. Indeed, from~\eqref{Gades} it follows by a straightforward calculation that we have an identity
\begin{gather}\label{GbNM}
\frac{G(v-ib_{+-}/2)}{G(v+ib_{+-}/2)}=\frac{\prod\limits_{k=1}^M2c_+\left(v+\frac{i}{2}\big( (N+1)a_+ +(M+1-2k)a_-\big)\right)}{\prod\limits_{j=0}^N2c_-\left(v+\frac{i}{2}\big( Ma_- +(N-2j)a_+\big)\right)}.
\end{gather}
Using this identity several times (together with~\eqref{Geval} for the factor $G(ib_{+-}-ia)$), we deduce from~\eqref{Riii} and~\eqref{cRr} the representation
\begin{gather}
   \cR_r(b_{+-};x,y)=\frac{1}{a_-}\frac{\prod\limits_{j=1}^N 2\sin(j\pi a_+/a_-)}{\prod\limits_{k=1}^M 2\sin(k\pi a_-/a_+)}
\frac{\prod\limits_{k=-M}^M2s_+(x-ika_-)}{\prod\limits_{j=-N}^N2s_-(x-ija_+)}\label{cRbNM}
 \\
\phantom{\cR_r(b_{+-};x,y)=}{}   \times
\int_{\R}dze^{i\alpha zx}\frac{\prod\limits_{k=1}^M4c_+\left(z+\frac{1}{2}\big(\pm y+i(N+1)a_++i(M+1-2k)a_-\big)\right)}{\prod\limits_{j=0}^N4c_-\left(z+\frac{1}{2}\big(\pm y+iMa_-+i(N-2j)a_+\big)\right)}.\nonumber
\end{gather}
The denominator of the integrand has no zero for $z\in\R$ unless $M$ is odd and $N$ is even. The zeros of the corresponding factor $s_-(z\pm y/2)$ are then matched by the zeros of the factor $s_+(z\pm y/2)$ of the numerator. After a suitable contour shift, we can expand the numerator product into exponentials, yielding a sum of convergent integrals (recall we require $b_{+-}\in(0,2a)$). When $M$ is even or $N$ is odd, we can do the same without a contour shift. Each of the integrals is then basically of the same form as previously evaluated for the case $M=0$. (More in detail, the same $2N+2$ poles arise in the period strip  for the $c_-$-product.)

From these observations the general structure anticipated in~\eqref{cRM} and~\eqref{Mform} readily follows, provided $n>0$ and $m\le 0$, or $m>0$ and $n\le 0$. In Section~III of~\cite{GLFII} we studied the functions~\eqref{cRbNM} in considerable detail, but it is beyond our scope to derive the explicit form used there directly from their representation~\eqref{cRbNM}. We do add that it seems plausible that the factorization exhibited in equations~(3.3)--(3.4) of~\cite{GLFII} can be understood by a more ref\/ined analysis of the above sum of contour integrals. In any case, we repeat that equality of the pertinent functions follows from the uniqueness argument explained in~Subsection~\ref{section2.3}.

\subsection{The nonrelativistic limit}\label{section4.2}

We begin this subsection with a remark addressed to physicist readers, who may care about dimension issues. In our paper~\cite{GLFII}, which we had occasion to cite several times in the previous subsection, the variable~$y$ of the present paper was denoted by~$p$. This is a widely used notation for the spectral variable in nonrelativistic quantum mechanics, where~$p$  is viewed as a momentum. In our relativistic setting, however, it is far more natural to view the scale parameters~$a_{+}$ and~$a_-$ as having dimension [position], and then the `geometric' and `spectral' variables~$x$ and~$y$ both have dimension [position] as well. (To be more specif\/ic, one of the scale parameters can be viewed as an interaction range, and the other one as the Compton wave length $\hbar/mc$ of the relativistic particles under consideration.) This goes to explain our change from~$p$ to~$y$.

Of course, from a mathematical viewpoint such notational issues and the notion of dimension may be ignored. When taking the nonrelativistic limit, however, these physical considerations point the way. We need to let the speed of light~$c$ go to~$\infty$, so one of the scale parameters should go to~0. In particular, we cannot retain modular invariance. Accordingly, we f\/irst set
\begin{gather}\label{repar}
a_+= 2\pi/\mu,\qquad  a_-= \hbar \beta,\qquad \mu,\hbar,\beta >0.
\end{gather}
Here, we view~$\beta$ as $1/mc$, with~$m$ the particle mass, and we trade~$a_+$ for a parameter~$\mu$ with dimension [position]$^{-1}$ to avoid a great many factors~$\pi$. Next, the spectral variable~$y$ (dual position) is replaced by the momentum variable
\begin{gather}\label{pdef}
p= \mu y/\beta.
\end{gather}
Finally, the coupling parameter $b$ (with dimension [position]) is replaced by
\begin{gather}\label{bg}
b=g\beta,\qquad g>0,
\end{gather}
so that $g$ has dimension [action]=[position]$\times$[momentum] as well as Planck's constant~$\hbar$.

With these changes in ef\/fect, it is easy to verify the expansion
\begin{gather}
A_+(b;x)= 2+ \beta^2 A +O\big(\beta^4\big),\qquad \beta \to 0,
\end{gather}
where
\begin{gather}
A:= -\hbar^2 \partial_x^2-g\hbar \mu \coth (\mu x/2)\partial_x-g^2\mu^2/4,
\end{gather}
and the limits
\begin{gather}
\lim_{\beta \to 0}A_-(b;x)=\exp(-i\pi g/\hbar)T_{2i\pi/\mu}^x +(i\to -i)=:M,\qquad x>0,
\\
\lim_{\beta \to 0}A_+(b;y)=\frac{p-ig\mu}{p}T_{i\hbar \mu}^p+(i\to -i)=: \hat{A},
\end{gather}
while $A_-(b;y)$ has no sensible limit. The eigenvalue of $A_+(b;x)$ on $\cR(b;x,y)$ satisf\/ies
\begin{gather}\label{cpexp}
2c_+(y)=
2+\beta^2 p^2/4 +O\big(\beta^4\big),
\end{gather}
so that the eigenvalue of $A$ becomes $p^2/4$, while the eigenvalues of the monodromy operator $M$ and dual A$\De$O $\hat{A}$ remain $2\cosh(\pi p/\hbar\mu)$ and $2\cosh(\mu x/2)$, resp.

Likewise, the similarity-transformed A$\De$Os~$\cA_{\pm}$~\eqref{cAdef} and Hamiltonians~$H_{\pm}$~\eqref{H} yield
\begin{gather}
\cA_+(b;x),H_+(b;x)= 2+ \beta^2 H +O\big(\beta^4\big),\qquad \beta \to 0,
\\
H:=  -\hbar^2 \partial_x^2+\frac{g(g-\hbar)\mu^2}{4\sinh^2(\mu x/2)},
\\
\label{cM}
\lim_{\beta \to 0}\cA_-(b;x)=\lim_{\beta \to 0}H_-(b;x)= T_{2i\pi/\mu}^x+T_{-2i\pi/\mu}^x=:\cM,\qquad x>0,
\\
\lim_{\beta \to 0}\cA_+(b;y)=T_{i\hbar \mu}^p+\frac{p+ig\mu}{p}T_{-i\hbar \mu}^p\frac{p-ig\mu}{p}=:\hat{\cA},
\\
\lim_{\beta \to 0}H_+(b;y)=
\sum_{\tau=+,-}\left( \frac{p-\tau ig\mu}{p}\right)^{1/2}T^p_{\tau i\hbar\mu}\left(\frac{p+\tau ig\mu}{p}\right)^{1/2}=:\hat{H}.
\end{gather}
(We take this opportunity to point out that the `$BC_1$ version' of the monodromy operator~\eqref{cM} specif\/ied in~equation~(8.20) of~\cite{Rsurv} contains an error: the phase factor on the r.h.s.\ should be replaced by~1.)

The operators $A$, $M$ and $H$, $\cM$ are related by a similarity transformation with the limit function
\begin{gather}\label{cxlim}
\lim_{\beta \to 0}\phi(b)/c(b;x) =w_{\rm nr}(g/\hbar;\mu x/2)^{1/2},\qquad x>0,
\end{gather}
where
\begin{gather}\label{wnr}
w_{\rm nr}(\lambda;r)\equiv (2\sinh r)^{2\lambda},\qquad \lambda,r>0,
\end{gather}
the operators $\hat{A}$ and $\hat{\cA}$ by similarity with
\begin{gather}\label{cylim}
\lim_{\beta \to 0}G(ia-2ib)G(ib-ia)/c(b;y)=2/\hat{c}_{\rm nr}(g/\hbar;p/\hbar\mu),
\end{gather}
where the Harish-Chandra $c$-function is given by
\begin{gather}\label{cnr}
 \hat{c}_{\rm nr}(\lambda;k)\equiv \frac{2\Gamma(2\lambda)}{\Gamma(\lambda)}\frac{\Gamma(ik)}{\Gamma(\lambda+ik)},
 \end{gather}
and the operators $\hat{A}$ and $\hat{H}$ by similarity with
$\hat{w}_{\rm nr}(g/\hbar;p/\hbar\mu)^{1/2}$, where
\begin{gather}
\hat{w}_{\rm nr}(\lambda;k)\equiv 1/ \hat{c}_{\rm nr}(\lambda;\pm k).
\end{gather}
The limits \eqref{cxlim} and \eqref{cylim} can be readily verif\/ied via the def\/inition~\eqref{cb} of the relativistic $c$-function and the above reparametrizations~\eqref{repar}--\eqref{bg} by using the $G$-limits~\eqref{GhGr}--\eqref{Gnr}. The functions~\eqref{wnr} and \eqref{cnr} are normalized so that
\begin{gather}\label{cwnorm}
\lim_{\lambda\to 0}w_{\rm nr}(\lambda;r)=1,\qquad
\lim_{\lambda\to 0}\hat{c}_{\rm nr}(\lambda;k)=1.
\end{gather}

Next, we study the nonrelativistic limit of the $\cR$-function. To this end we f\/irst use the def\/inition~\eqref{Rdef}, the scaling property~\eqref{scale} and the limit
\begin{gather}\label{limRF}
\lim_{t \to 0}R(\pi,t,t {\bf c};v,t
u)={}_2F_1\big(\hat{c}_0+iu,\hat{c}_0-iu,c_0+c_2+1/2;-\sh^2v\big),
\\
\hat{c}_0\equiv (c_0+c_1+c_2+c_3)/2,
\end{gather}
established and discussed in~I. We write the limit
\begin{gather}\label{limcRF}
\lim_{\beta \to 0}\cR(2\pi/\mu,\hbar\beta,g\beta;x,\beta p/\mu)=:\psi_{\rm nr}(g/\hbar;\mu x/2,p/\hbar\mu),
\end{gather}
in terms of the dimensionless quantities
\begin{gather}\label{dless}
\lambda \equiv g/\hbar,\qquad  r\equiv \mu x/2,\qquad k\equiv p/\hbar\mu,
\end{gather}
already used in~\eqref{wnr} and~\eqref{cnr}. The result reads{\samepage
\begin{gather}\label{psi1}
\psi_{\rm nr}(\lambda;r,k)={}_2F_1\big((\lambda+ik)/2,(\lambda-ik)/2,\lambda +1/2;-\sinh^2 r\big).
\end{gather}
(See~\cite{koo} for a limit that is related to~\eqref{limcRF}, cf.~\eqref{qGeg}.)}

Likewise, the alternative representation~\eqref{cR2} entails
\begin{gather}\label{psi2}
\psi_{\rm nr}(\lambda;r,k)={}_2F_1\big(\lambda+ik,\lambda-ik,\lambda +1/2;-\sinh^2 (r/2)\big),
\end{gather}
whereas~\eqref{cR3}--\eqref{cR4} again give rise to~\eqref{psi1}. The equality of~\eqref{psi1} and~\eqref{psi2} can be rewritten as
\begin{gather}\label{Fq}
{}_2F_1(a,b,a+b+1/2;4w(1-w))={}_2F_1(2a,2b,a+b+1/2;w),
\end{gather}
which is a well-known quadratic transformation, cf.~e.g.~\cite[p.~125]{aar}.

Using other familiar features of the hypergeometric function, it is not dif\/f\/icult to verify that the operators $A$, $M$ and $\hat{A}$ do have the expected eigenvalues $p^2/4$, $2\cosh(p/\hbar\mu)$ and $2\cosh(\mu x/2)$ on the limit function~$\psi_{\rm nr}(g/\hbar;\mu x/2,p/\hbar\mu)$. More specif\/ically, for $A$ this amounts to the ODE satisf\/ied by ${}_2F_1$ and for $\hat{A}$ this involves the contiguous relations. The $M$-eigenfunction property follows by using the known analytic continuation of ${}_2F_1(a,b,c;w)$ across the logarithmic branch cut $w\in[1,\infty)$.

The above specialization of the hypergeometric function  basically yields the so-called conical (or Mehler) function. To be specif\/ic, the latter can be def\/ined by
\begin{gather}\label{conf}
P_{ik-1/2}^{1/2-\lambda}(\cosh r)\equiv \frac{(\sinh r)^{\lambda-1/2}}{2^{\lambda -1/2}\Gamma(\lambda +1/2)}
\ {}_2F_1(\lambda+ik,\lambda-ik,\lambda +1/2;(1-\cosh r)/2),
\end{gather}
cf.~\cite[equation~(14.3.15)]{dlmf} and the hypergeometric function occurring here equals the one in~\eqref{psi2}.

We now consider the nonrelativistic limit of the minimal representations of $\cR(b;x,y)$ derived in Section~\ref{section3}. We were not  able to obtain a sensible limit for the second one, given by~\eqref{Rii}. For the remaining four, however, the limit can be handled in a sense to be explained shortly. For expository reasons we f\/irst list the resulting representations for~$\psi_{\rm nr}(\lambda;r,k)$:
\begin{gather}\label{Rlimi}
\frac{\Gamma(2\lambda)}{2^{\lambda}\Gamma(\lambda)^2}\int_{\R}dt\frac{1}{(\cosh r+\cosh t)^{\lambda}}\cos\left(k\ln \left(\frac{\cosh((t+r)/2)}{\cosh((t-r)/2)}\right)\right),
\\
\label{Rlimiii}
\frac{\Gamma(2\lambda)(\sinh r)^{1-2\lambda}}{2^{2\lambda +1}\pi}
\int_{\R}dt\frac{\Gamma((it-\lambda \pm ik +c+1)/2)}{\Gamma((it+\lambda \pm ik +c+1)/2)}\exp((i(t-ic)r),\qquad c>\lambda -1,
\\
\label{Rlimiv}
\frac{\Gamma(2\lambda)}{4\pi\Gamma(\lambda)^2\Gamma(\lambda \pm ik)}\int_{\R}dt\Gamma((it+\lambda \pm ik)/2)\Gamma((-it+\lambda \pm ik)/2)\exp(itr),
\\
\label{Rlimv}
\frac{\Gamma(2\lambda)}{2^{\lambda}\Gamma(\lambda \pm ik)}\int_{\R}dt\frac{\exp(itk)}{(\cosh r +\cosh t)^{\lambda}}.
\end{gather}

We proceed to discuss these formulas. First, we note that they are derived under the assumption
\begin{gather}
\lambda,r,k\in(0,\infty),
 \end{gather}
and that this implies that the integrals in~\eqref{Rlimi}, \eqref{Rlimiv} and~\eqref{Rlimv} are absolutely convergent.
The integral in~\eqref{Rlimiii}, however, is only absolutely convergent for $\lambda>1/2$; For $\lambda\in(0,1/2]$ it should be viewed as a Fourier transform in the sense of tempered distributions.

Second, we compare these formulas to results in~\cite{dlmf}, where a host of representations for $_2F_1$ and its conical function specialization are listed. Formula~\eqref{Rlimv} can be readily found there: It can be obtained from equation~(14.12.4), which can be written
\begin{gather}\label{Prep}
P_{ik-1/2}^{1/2-\lambda}(\cosh r)=\sqrt{\frac{2}{\pi}} \frac{\Gamma(\lambda)(\sinh r)^{\lambda-1/2}}{\Gamma(\lambda\pm ik)}
\int_{0}^{\infty}dt\frac{\cos kt}{(\cosh r +\cosh t)^{\lambda}}.
\end{gather}
(This involves the duplication formula of the gamma function, equation~(5.5) in~\cite{dlmf}.) As they stand, the three remaining representations do not occur in~\cite{dlmf}. However, as was pointed out by a referee, they can also be tied in with results in the vast literature connected to the hypergeometric function. To begin with, formula~\eqref{Rlimiv} can be derived (with some ef\/fort) by combining equations~(15.8.14) and~(15.6.7) in~\cite{dlmf}. It seems that the formulas~\eqref{Rlimi} and~\eqref{Rlimiii} cannot be obtained by using only~\cite{dlmf} or some other standard reference book. Even so, they agree with known results. Indeed, \eqref{Rlimi} amounts to equation~(2.3) in the paper~\cite{bako}, whereas~\eqref{Rlimiii} can be derived by combining several sources. Specif\/ically, the integral in~\eqref{Rlimiii} can be viewed as a special case of the Meyer $G$-function, cf.~Section~16.17 in~\cite{dlmf} and p.~144, (2) of~\cite{luke}. After contour deformation, a residue calculation leads to a formula involving a linear combination of two~$_2F_1$'s with gamma function coef\/f\/icients, cf.~equations~(16.17.2) and (16.17.3) in~\cite{dlmf} or~(7) in~\cite{luke}. Finally, it follows from equation~(3.2(27)) in~\cite{erd} that the latter formula yields the conical function as represented by~\eqref{Prep}.

Third, none of the four representations~\eqref{Rlimi}--\eqref{Rlimv} has been obtained with complete rigor. The diff\/iculty is to obtain uniform tail bounds on the pertinent integrands that allow an application of the dominated convergence theorem. (In fact, to date a similar snag has not yet been obviated for the limit transition~\eqref{limRF} either.)
In this connection we should add that we were unable to verify directly that each of the four formulas gives rise to a joint eigenfunction of the operators $A$, $M$ and $\hat{A}$ with eigenvalues $p^2/4$, $2\cosh(p/\hbar\mu)$ and $2\cosh(\mu x/2)$. On the other hand, a direct check of the joint eigenfunction properties of the f\/ive relativistic representations~\eqref{Ri}--\eqref{Rv} seems not feasible either.

We continue by sketching the derivation of the four formulas~\eqref{Rlimi}--\eqref{Rlimv}. First, we observe that any factor of the form
\begin{gather}
G(a_{+},a_-;ia-ita_-),
\end{gather}
with $t$ not depending on $a_-$, can be treated via~\eqref{GhGr}. Indeed, a scaling by $a_+$ yields
\begin{gather}
G(1,\kappa;i/2+\kappa(i/2-it))\sim \frac{2\pi \sqrt{\kappa}}{\Gamma(t)\exp(t\ln(2\pi \kappa))},\qquad \kappa \to 0.
\end{gather}
In particular, this yields not only the asymptotics of the numerical prefactors, viz.,
\begin{gather}
\frac{1}{\sqrt{a_+a_-}G(ia-2ib)}\sim \frac{\Gamma(2\lambda)\exp(2\lambda \ln(\beta\hbar\mu))}{2\pi\beta\hbar},\qquad \beta\to 0,
\\
\frac{G(ia-ib)^2}{\sqrt{a_+a_-}G(ia-2ib)}\sim \frac{\mu\Gamma(2\lambda)}{\Gamma(\lambda)^2},\qquad \beta\to 0,
\end{gather}
but also that of the $y$-dependent prefactor in~\eqref{Riv}--\eqref{Rv}:
\begin{gather}
G(ia-ib\pm y)\sim \frac{2\pi\beta\hbar\mu\exp(-2\lambda\ln(\beta\hbar\mu))}{\Gamma(\lambda\pm ik)},\qquad \beta\to 0.
\end{gather}

To handle the $x$-dependent prefactor in~\eqref{Rii}--\eqref{Riii}, however, \eqref{GhGr} is of no help and~\eqref{Gnr} seems not to apply either. But in fact we can use the $G$-A$\De$Es~\eqref{Gades} to f\/irst write
\begin{gather}
\frac{G(x+ia-ib)}{G(x-ia+ib)}=2is_-(x-ib)\frac{G(x-ia_+/2+ia_-/2-ib)}{G(x-ia_+/2-ia_-/2+ib)},
\end{gather}
and then~\eqref{Gnr} can be invoked to deduce the dominant asymptotics
\begin{gather}\label{Gx}
G(\pm x+ia-ib)\sim \exp\left( \frac{\pi x}{\beta\hbar}\right)(2\sinh r)^{1-2\lambda},\qquad x>0,\qquad \beta\to 0.
\end{gather}

The plane waves in~\eqref{Rii}--\eqref{Rv} become
\begin{gather}
\exp(i\alpha zx)=\exp(i\mu zx/\beta\hbar),\qquad
\exp(i\alpha zy)=\exp(i zp/\hbar).
\end{gather}
For the f\/irst and second plane wave we now switch to a new variable $t$ given by
\begin{gather}
z\to \beta\hbar t/2,\qquad z\to t/\mu,
\end{gather}
so that they turn into
$\exp(itr)$ and $\exp(itk)$,
resp. Likewise, in~\eqref{Ri} we change $z$ to $t/\mu$ to get dimensionless $G$-function arguments.

With these variable changes in place, we proceed to look at the asymptotic behavior of the $G$-ratios in~\eqref{Ri}--\eqref{Rv}. For the f\/irst case this is immediate from~\eqref{Gnr}, and this easily yields~\eqref{Rlimi} when the pointwise limit is interchanged with the integration. (As alluded to above, an $L^1$-bound uniform for $\beta$ near~0 is needed to make the interchange rigorous. As well as in the next cases, no such bound is available for now.) For the second case~\eqref{Rii}, it seems not possible to get from the pointwise behavior of the integrand (with $x>0$) as $\beta$ goes to~0 a factor $\exp(-\pi x /\beta\hbar)$ that takes care of the diverging factor $\exp(\pi x/\beta\hbar)$ coming from the prefactor~\eqref{Gx}. By contrast, for the third case we may and will make a shift
\begin{gather}
t\to t +i/(2\mu \beta\hbar)-ic,
\end{gather}
with $c\in\R$ chosen so as to stay away from poles while shifting and taking $\beta $ to 0. This cancels the diverging factor and results in~\eqref{Rlimiii} via~\eqref{GhGr}. Finally, an application of~\eqref{GhGr} and~\eqref{Gnr} leads to the limits~\eqref{Rlimiv} and~\eqref{Rlimv}, resp.

To conclude this subsection, we point out that in view of~\eqref{ER} and \eqref{cRdef} the nonrelativistic limit of $\rE(b;x,y)$ is given by
\begin{gather}
 \frac{2w_{\rm nr}(\lambda;r)^{1/2}}{\hat{c}_{\rm nr}(\lambda;k)}\psi_{\rm nr}(\lambda;r,k)=:\rE_{\rm nr}(\lambda;r,k),
\end{gather}
cf.~\eqref{cxlim}--\eqref{cnr}. It has the unitary asymptotics
\begin{gather}
\rE_{\rm nr}(\lambda;r,k)\sim \exp(irk)-\hat{u}_{\rm nr}(\lambda;-k)\exp(-irk),\qquad r\to\infty,
\end{gather}
with the scattering function
\begin{gather}
\hat{u}_{\rm nr}(\lambda;k)\equiv -\frac{\hat{c}_{\rm nr}(\lambda;k)}{\hat{c}_{\rm nr}(\lambda;-k)}=-\frac{\Gamma(ik)\Gamma(\lambda-ik)}{\Gamma(-ik)\Gamma(\lambda+ik)}.
\end{gather}
The latter is normalized so that it equals 1 for $\lambda =1$, just as $u(b;z)$~\eqref{u} is normalized to equal 1 for $b=a_{\pm}$. In this connection we would like to add that from~\eqref{Rzero} and~\eqref{ER}--\eqref{cb} one readily deduces
\begin{gather}
\rE(a_{\pm};x,y)=2i\sin(\pi xy/a_+a_-).
\end{gather}
Hence the reparametrizations~\eqref{repar}--\eqref{bg} and~\eqref{dless} yield
\begin{gather}
E_{\rm nr}(1;r,k)=2i\sin kr.
\end{gather}
Accordingly, the `free' theory with which the scattering is compared is given by the sine transform (and not by the cosine transform, which arises for $b=\lambda=0$, cf.~\eqref{ER0}).

\section{The relativistic Toda case}\label{section5}

\subsection{Taking the relativistic Toda limit}\label{section5.1}

Throughout this section, we require that the parameters $a_+$ and $a_-$ be positive. It is also convenient to require
\begin{gather}\label{xyp}
(x,y)\in(0,\infty)^2,
\end{gather}
until further notice.
In keeping with our outline in the Introduction, we begin by considering the $b$-values~\eqref{gam}. In this case the $w$-function~\eqref{w} reads
\begin{gather}\label{wmod}
w(a-i\gamma;z)=G(\pm z +ia)/G(\pm z -\gamma),\qquad a=(a_{+}+a_{-})/2,\qquad \gamma\in\R,
\end{gather}
and hence is no longer real-valued for real $z$ and $\gamma\ne 0$. By contrast, the $u$-function~\eqref{u} is given by
\begin{gather}\label{ug}
u(a-i\gamma;z)=-G(z\pm \gamma)/G(z\pm ia),
\end{gather}
so it is still unitary for real $z$; moreover, it is even in $\gamma$. The Hamiltonians~\eqref{H} can be written
\begin{gather}
H_{\de}(a-i\gamma;z)=\left(\frac{c_{\de}(z+ia_{-\de}/2\pm\gamma)}{s_{\de}(z)s_{\de}(z+ia_{-\de})}\right)^{1/2}\exp(ia_{-\de}\partial_z)+(i\to -i),
\end{gather} so they remain formally self-adjoint for real $z$; they are also  even in $\gamma$.

Next, we consider the f\/ive representations of the joint eigenfunction ${\rm F}(a-i\gamma;x,y)$ of the four Hamiltonians $H_{\pm}(a-i\gamma;x)$ and $H_{\pm}(a-i\gamma;y)$. Combining~\eqref{FR} and \eqref{cRr} with~\eqref{Ri}--\eqref{Rv}, these are given by
\begin{gather}\label{Fi}
\frac{G(-\gamma)}{\sqrt{a_{+}a_-}}\left(\frac{G(\pm x +ia)}{G(\pm x -\gamma)}\frac{G(\pm y +ia)}{G(\pm y -\gamma)}\right)^{1/2}\iR dz \frac{G(z\pm (x-y)/2-ia/2-\gamma/2)}{G(z\pm (x+y)/2+ia/2+\gamma/2)},
\\
\label{Fii}
\frac{G(-\gamma)}{\sqrt{a_{+}a_-}}\left(\frac{G(\pm x +ia)}{G(\pm x +\gamma)}\frac{G(\pm y +ia)}{G(\pm y -\gamma)}\right)^{1/2}\iR dz G(\pm z\pm x/2-ia/2+\gamma/2)\exp(i\alpha zy),
\\
\label{Fiii}
\frac{G(\gamma)}{\sqrt{a_{+}a_-}}\left(\frac{G(\pm x +ia)}{G(\pm x +\gamma)}\frac{G(\pm y +ia)}{G(\pm y -\gamma)}\right)^{1/2}\iR dz G(\pm z\pm y/2-ia/2-\gamma/2)\exp(i\alpha zx),
\\
\label{Fiv}
\frac{G(-\gamma)}{\sqrt{a_{+}a_-}}\left(\frac{G(\pm x +ia)}{G(\pm x -\gamma)}\frac{G(\pm y +ia)}{G(\pm y +\gamma)}\right)^{1/2}\iR dz G(\pm z\pm y/2-ia/2+\gamma/2)\exp(i\alpha zx),
\\
\label{Fv}
\frac{G(\gamma)}{\sqrt{a_{+}a_-}}\left(\frac{G(\pm x +ia)}{G(\pm x -\gamma)}\frac{G(\pm y +ia)}{G(\pm y +\gamma)}\right)^{1/2}\iR dz G(\pm z\pm x/2-ia/2-\gamma/2)\exp(i\alpha zy).
\end{gather}
(The square roots are positive for $\gamma=0$.) As they stand, none of these representations yields a manifestly real-valued function for $\gamma\ne 0$. However, comparing~\eqref{Fiii} and~\eqref{Fiv}, we see that these formulae are related by a complex conjugation (take $z\to -z$ in one of them to check this).
Likewise, \eqref{Fii} and~\eqref{Fv} are related by a complex conjugation. Since the f\/ive formulae yield the same function
${\rm F}(a-i\gamma;x,y)$, this function is in fact real-valued.

This reality feature can be tied to the $b\to 2a-b$ symmetry of the E-function, cf.~\eqref{Esym}. Indeed, the latter invariance implies that E$(a-i\gamma;x,y)$ is even in~$\gamma$. Now since the $u$-function~\eqref{ug} and the phase $\phi(a-i\gamma)$ (given by~\eqref{phib}) are manifestly even in $\gamma$, it follows that F$(a-i\gamma;x,y)$ is even in $\gamma$, cf.~\eqref{Fdef}. Comparing once again~\eqref{Fiii} with~\eqref{Fiv}, and~\eqref{Fii} with~\eqref{Fv}, we see that these formulae are also related by f\/lipping the sign of $\gamma$, in accord with evenness.

Substituting
\begin{gather}\label{Tsub}
\gamma\to \eta+\Lambda,\qquad x\to x+\Lambda,
\end{gather}
we are now prepared to study the Toda limit $\Lambda\to\infty$. First, for $r\in\R$ we have
\begin{gather}\label{Vlim}
\lim_{\Lambda\to\infty}\frac{c_{\de}(x+ ira_{-\de}/2+\eta+2\Lambda)c_{\de}(x+ ira_{-\de}/2-\eta)}{s_{\de}(x+\Lambda)s_{\de}(x+ ira_{-\de}+\Lambda)}=1+e_{\de}(-2x- ira_{-\de}+2\eta).
\end{gather}
Thus we obtain relativistic nonperiodic Toda Hamiltonians given by
\begin{gather}
H_{\de}^T(\eta;x)     \equiv     \lim_{\Lambda\to\infty}H_{\de}(a-i\eta-i\Lambda;x+\Lambda)
\nonumber \\
\phantom{H_{\de}^T(\eta;x)}{}   =    [1+e_{\de}(-2x-ia_{-\de}+2\eta)]^{1/2}\exp(ia_{-\de}\partial_x)+(i\to -i).\label{HT}
\end{gather}
(The square roots are positive for $x\to\infty$.) Clearly, these are formally self-adjoint on $L^2(\R,dx)$. In this connection we point out that in view of the diverging $x$-shift, we may and will from now on allow $x$ to vary over $\R$ in the Toda quantities, whereas we continue to require that $y$ be positive.

Next, we note
\begin{gather}
\lim_{\Lambda\to\infty}e_{\de}(-2\Lambda)c_{\de}(y+ ira_{-\de}/2+\eta+\Lambda)c_{\de}(y+ ira_{-\de}/2-\eta-\Lambda)
=e_{\de}(2\eta)/4.
\end{gather}
Hence we get dual Hamiltonians
\begin{gather}
\hat{H}_{\de}^T(\eta;y)     \equiv     \lim_{\Lambda\to\infty}e_{\de}(-\Lambda)H_{\de}(a-i\eta-i\Lambda;y)
\nonumber \\
\phantom{\hat{H}_{\de}^T(\eta;y)}{}   =    \frac{e_{\de}(\eta)}{2}s_{\de}(y)^{-1/2}\big(\exp(ia_{-\de}\partial_y)+\exp(-ia_{-\de}\partial_y)\big)s_{\de}(y)^{-1/2},\label{HTd}
\end{gather}
which are formally positive on $L^2((0,\infty),dy)$.

To obtain the Toda limit of the joint eigenfunction F$(a-i\gamma;x,y)$ involves a greater ef\/fort.
The key tool is the asymptotics~\eqref{Gas} of the hyperbolic gamma function. This enables us to show that the limit
\begin{gather}\label{Flim}
{\rm F}^T(\eta;x,y)\equiv \lim_{\Lambda\to\infty}{\rm F}(a-i\eta-i\Lambda;x+\Lambda,y)
\end{gather}
exists for each of the f\/ive representations~\eqref{Fi}--\eqref{Fv}. The details now follow.

To start with, the asymptotic behavior for $\Lambda\to\infty$ of the f\/ive prefactors can be assembled from the three formulae
\begin{gather}\label{as1}
G(x+\eta+\Lambda)\sim e^{-i\chi}\exp\left(\frac{-i\pi\alpha}{4}\big((\eta+\Lambda)^2+x^2+2x(\eta+\Lambda)\big)\right),
\\
\label{as2}
G(\pm y-\eta-\Lambda)\sim e^{2i\chi}\exp\left(\frac{i\pi\alpha}{2}\big((\eta+\Lambda)^2+y^2\big)\right),
\\
\label{as3}
G(\pm (x+\Lambda)+ia)\sim \exp(\pi\alpha a(x+\Lambda)).
\end{gather}
Next, consider the integrand of~\eqref{Fi} with the substitutions~\eqref{Tsub}. Two of the four $G$-functions are invariant, and the remaining two yield
\begin{gather}
\frac{G(z-(x-y)/2-ia/2-\eta/2-\Lambda)}{G(z+(x+y)/2+ia/2+\eta/2+\Lambda)}\nonumber\\
\qquad{} \sim e^{2i\chi}\exp\left(\frac{i\pi\alpha}{2}\big([z+y/2]^2+[(x+ia+\eta)/2+\Lambda]^2\big)\right).\label{Intisub}
\end{gather}
If we now combine the $\Lambda$-dependent terms coming from the prefactor in~\eqref{Fi}, then we see that they cancel the $\Lambda$-dependent term in~\eqref{Intisub}. The product of the remaining terms is readily verif\/ied to be given by
\begin{gather}\label{aux}
\left(\frac{G(\pm y+ia)}{a_{+}a_-G(x-\eta)}\right)^{1/2} \exp(3i\chi/2)\frac{G(z+(x-y)/2-ia/2-\eta/2)}{G(z-(x+y)/2+ia/2+\eta/2)}\exp\left( \frac{i\alpha}{4}M\right),
\end{gather}
with
\begin{gather}
M\equiv 2z^2+2zy-y^2/2-iax-a^2/2+ia\eta.
\end{gather}
Finally, since we integrate $z$ in~\eqref{aux} over~$\R$, we may shift $z$ by $y/2$, yielding the limit function
\begin{gather}
{\rm F}^T(\eta;x,y)  =   \left(\frac{G(\pm y+ia)}{a_{+}a_-G(x-\eta)}\right)^{1/2} \exp(3i\chi/2)\exp\left( \frac{i\alpha}{4}\big(y^2-ia(x-\eta)-a^2/2\big)\right)
\nonumber \\
\phantom{{\rm F}^T(\eta;x,y)  =}{}  \times \iR dz G(\pm z +(x-\eta-ia)/2)
 \exp\left( \frac{i\alpha}{2}\big(z^2+2zy\big)\right).\label{T1}
 \end{gather}

Turning to the integrand of~\eqref{Fii}, the substitution~\eqref{Tsub} again leaves two of the four $G$-functions unchanged, as well as the plane wave factor. The remaining $G$-product has asymptotics
\begin{gather}
G(\pm z +(x-ia+\eta)/2+\Lambda)\sim e^{-2i\chi}\exp\left(\frac{-i\pi\alpha}{2}\big(z^2+[(x-ia+\eta)/2+\Lambda]^2\big)\right).
\end{gather}
Combining this with the asymptotics of the prefactor following from~\eqref{as1}--\eqref{as3}, the $\Lambda$-dependent terms drop out. Taking $z\to -z$ in the resulting limit function yields
\begin{gather}
{\rm F}^T(\eta;x,y)  =   \left(\frac{G(\pm y+ia)}{a_{+}a_-G(-x+\eta)}\right)^{1/2} \exp(-3i\chi/2)\exp\left( \frac{-i\alpha}{4}\big(y^2+ia(x-\eta)-a^2/2\big)\right)
\nonumber \\
\phantom{{\rm F}^T(\eta;x,y)  =}{}  \times \iR dz G(\pm z -(x-\eta+ia)/2)
 \exp\left( \frac{-i\alpha}{2}\big(z^2+2zy\big)\right).\label{T2}
 \end{gather}
Comparing this representation to~\eqref{T1}, we see that each of the factors on the right-hand side is matched by its complex-conjugate. Thus, the equality of~\eqref{T1} and~\eqref{T2} is in keeping with the real-valuedness of ${\rm F}^T(\eta;x,y)$ following from its being the limit of a real-valued function.

Proceeding in the same way for~\eqref{Fv}, we obtain as its limit again~\eqref{T1}. Of course this should be expected, since the factors on the right-hand side of~\eqref{Fv} and~\eqref{Fii} are related by complex conjugation. On the other hand, the equality of the limits of~\eqref{Fi} and~\eqref{Fv} yields a~nontrivial check of the substantial limit calculations.

At face value, the representations~\eqref{Fiii} and~\eqref{Fiv} seem not to give rise to a sensible Toda limit. In fact, however, they do, but it is expedient to postpone the details. First, we rewrite the representations~\eqref{T1} and \eqref{T2} in a more telling form, by bringing in the $G$-cousins~$G_L$ and~$G_R$, cf.~\eqref{GRDef}--\eqref{GLDef}. Indeed, a straightforward calculation yields the equivalent representations
\begin{gather}\label{T1eq}
{\rm F}^T(\eta;x,y) = \left(\frac{G(\pm y+ia)}{a_{+}a_-G_R(x-\eta)}\right)^{1/2}e^{i\alpha y^2/4}\iR dz G_R(\pm z +(x-\eta)/2-ia/2)
 e^{i\alpha zy},
 \end{gather}
and
\begin{gather}\label{T2eq}
{\rm F}^T(\eta;x,y) =  \left(\frac{G(\pm y+ia)}{a_+a_-G_L(\eta-x)}\right)^{1/2}e^{-i\alpha y^2/4}
\iR dz G_L(\pm z -(x-\eta)/2-ia/2)
 e^{i\alpha zy}.
 \end{gather}
 (Here, the square roots are positive for $x\to\infty$, cf.~\eqref{GRLas}.)

Now we use once again the Plancherel relation, as encoded in~\eqref{Plan}--\eqref{Four}. Taking f\/irst
\begin{gather}
f(p)=G_R(p+(x-\eta)/2-ia/2)\exp(i\alpha [p+(x-\eta+ia)/2][y/2+ia/4]),
\\
g(p)=G_R(-p+(x-\eta)/2-ia/2)\exp(i\alpha [-p+(x-\eta+ia)/2][- y/2+ia/4]),
\end{gather}
we can use~\eqref{cor} with $s=a/4$ to compute the Fourier transforms. This yields
\begin{gather}
\hat{f}(q)=\left(\frac{2\pi}{\alpha}\right)^{1/2}e^{-i\pi/4-2i\chi}
\exp\left( \frac{-i\alpha}{2}q(x-\eta+ia)\right)
G_L(q+y/2-3ia/4),
\\
\hat{g}(q)=\left(\frac{2\pi}{\alpha}\right)^{1/2}e^{-i\pi/4-2i\chi}
\exp\left( \frac{i\alpha}{2}q(x-\eta+ia)\right)
G_L(-q-y/2-3ia/4).
\end{gather}
The integral in~\eqref{T1eq} is therefore equal to
\begin{gather}
  \exp\left( \frac{-i\alpha}{4}(ia)(x-\eta+ia)\right)\frac{\alpha}{2\pi}\iR \hat{f}(q)\hat{g}(-q)dq
  \nonumber \\
\qquad {} =e^{-i\pi/2-4i\chi}\iR dq \exp\big({-}i\alpha(q+ia/4)(x-\eta+ia)\big)G_L(q\pm y/2-3ia/4).
\end{gather}

The new representation thus obtained can be somewhat simplif\/ied by reverting to the $G$-function, and by shifting the contour down by $a/4$ (recall $y>0$). In this way we obtain from~\eqref{T1eq} an alternative representation
\begin{gather}
{\rm F}^T(\eta;x,y)    =    \left(\frac{G(\pm y+ia)}{a_{+}a_-G_R(x-\eta)}\right)^{1/2}\exp\big(i\alpha y^2/8\big)
\nonumber \\
\phantom{{\rm F}^T(\eta;x,y)    =}{}  \times  \int_{\R +i0}dzG(z\pm y/2-ia)\exp\big({-}i\alpha[z(x-\eta)+z^2/2]\big).\label{T3}
\end{gather}

Proceeding in the same way for~\eqref{T2eq}, we arrive at a fourth representation, namely
\begin{gather}
{\rm F}^T(\eta;x,y)    =      \left(\frac{G(\pm y+ia)}{a_+a_-G_L(\eta-x)}\right)^{1/2}\exp\big({-}i\alpha y^2/8\big)
\nonumber \\
\phantom{{\rm F}^T(\eta;x,y)    =}{}   \times \int_{\R +i0} dzG(z\pm y/2-ia)\exp\big({-}i\alpha[z(x-\eta)-z^2/2]\big).\label{T4}
\end{gather}
Comparing it to~\eqref{T3}, we deduce once again real-valuedness of the function ${\rm F}^T(\eta;x,y)$ on $\R^2\times(0,\infty)$.

Having these two new representations at hand, we can see with hindsight that they can also be obtained from~\eqref{Fiv} and~\eqref{Fiii}, respectively. Indeed, when we let
\begin{gather}
z\to z+ia/2-\eta/2-\Lambda/2
\end{gather}
in the integrand of~\eqref{Fiv}, then we obtain the two $G$-functions featuring in~\eqref{T3}, times two $\Lambda$-dependent ones. If we now proceed in the same way as before, using~\eqref{as1}--\eqref{as3} to handle the asymptotics of the prefactor, then we arrive once more at~\eqref{T3}, yielding a check on the rather extensive calculations. Likewise, \eqref{Fiii} gives rise to~\eqref{T4}.

\subsection[Asymptotic and analytic properties of F$^T(\eta;x,y)$]{Asymptotic and analytic properties of $\boldsymbol{{\rm F}^T(\eta;x,y)}$}\label{section5.2}

With the various representations of the function F$^T(\eta;x,y)$ at our disposal, several salient features can be readily derived. First,  it is remarkably easy to show from~\eqref{T1} that it has exponential decay for $x\to -\infty$ (as might be expected from the exponential divergence of the `potential' factors in the Hamiltonians~$H^T_{\pm}(\eta;x)$~\eqref{HT}). To be specif\/ic, we have a bound
\begin{gather}\label{expdec}
{\rm F}^T(\eta;x,y)=O(\exp(\alpha ax/4)),\qquad x\to -\infty.
\end{gather}

Inspecting~\eqref{T1}, it is clear that we need only show that the integral  yields a function that is $O(1)$ for $x\to -\infty$. To this end we point out that from~\eqref{Gas} we have estimates
\begin{gather}\label{Gb}
G(v-ia/2)=O(\exp(\mp \alpha av/4),\qquad v\to \pm\infty,
\end{gather}
and that no poles arise for real $v$. Hence the function $v\mapsto G(v-ia/2)$ is bounded on $\R$. If we now take $z\to z+(x-\eta)/2$ in~\eqref{T1}, then we obtain a factor $G(z-ia/2)$ times a factor that is bounded for $x,y,z,\eta\in\R$. Thus we can invoke the bound~\eqref{Gb} on the f\/irst factor to deduce that the integral is in fact bounded for $x$, $y$, $\eta$ varying over $\R$, completing the proof of~\eqref{expdec}.

It is also not hard to obtain the $x\to\infty$ asymptotics. To this end we start from the representation~\eqref{T3} and follow the reasoning below~\eqref{Enew}. Thus we shift the contour down by $\epsilon$, where $\epsilon>0$ is small enough so that only the simple poles at $z=\pm y/2$ are passed (recall our standing assumption $y>0$). The residues of the integral then follow from~\eqref{Gres}, yielding a residue sum
\begin{gather}
\left(\frac{G(\pm y+ia)}{G_R(x-\eta)}\right)^{1/2}\!\!\big(G(-y-ia)\exp(i\alpha y(x-\eta)/2) +
G(y-ia)\exp(-i\alpha y(x-\eta)/2)\big).\!\!\!
\end{gather}
Using the $G$-asymptotics~\eqref{Gas}, it is easily verif\/ied that the remainder integral vanishes for $x\to\infty$, so that we deduce
\begin{gather}
{\rm F}^T(\eta;x,y)\sim u^T(\eta;y)^{1/2}\exp(i\alpha xy/2)+
u^T(\eta;-y)^{1/2}\exp(-i\alpha xy/2),\qquad x\to\infty.
\end{gather}
Here we have introduced the Toda $u$-function
\begin{gather}
u^T(\eta;y)\equiv \exp(-i\alpha \eta y)G(-y+ia)/G(y+ia),
\end{gather}
which can also be written
\begin{gather}
u^T(\eta;y)=c^T(\eta;y)/c^T(\eta;-y),
\end{gather}
with the Toda $c$-function def\/ined by
\begin{gather}
c^T(\eta;y)\equiv \exp(-i\alpha \eta y/2)/G(y+ia).
\end{gather}
The corresponding weight function is given by
\begin{gather}\label{wTy}
w^T(y)\equiv 1/c^T(\eta;\pm y)=G(\pm y+ia)=4s_+(y)s_-(y),
\end{gather}
where we used the $G$-A$\De$Es~\eqref{Gades} in the last step.

The dual counterparts of these formulae are not obvious. To begin with, we have been unable to establish the large-$y$ asymptotics of F$^T(\eta;x,y)$. We conjecture, however, that this is given by
\begin{gather}
  {\rm F}^T(\eta;x,y)  \sim  e^{-i\chi/2-i\pi/8}G_R(x-\eta)^{1/2}\exp(i\alpha[y^2/4+ (x-\eta)y/2])
 \nonumber\\
\qquad{}   +
e^{i\chi/2+i\pi/8}G_R(x-\eta)^{-1/2}\exp(-i\alpha[y^2/4+ (x-\eta)y/2]),\qquad y\to\infty.\ \ (?)\label{yas}
\end{gather}
Even when this can be shown, it is not clear whether the function $G_R(x-\eta)$ can be viewed as an $S$-matrix for the dual dynamics. Indeed, the dual scattering theory seems quite unusual, just as at the classical level~\cite{toda}. Moreover, like the $w$-function  $w(a-i\eta-i\Lambda;x+\Lambda)$, the $u$-function $u(a-i\eta-i\Lambda;x+\Lambda)$ has no limit for $\Lambda\to \infty$.

On the other hand, the similarity transforms
\begin{gather}
\cA^{T}_{\de}(\eta;x)     \equiv       G_R(x-\eta)^{- 1/2} H^T_{\de}(\eta;x)  G_R(x-\eta)^{ 1/2}
\nonumber \\
\phantom{\cA^{T}_{\de}(\eta;x)}{}    =    \exp(-ia_{-\de}\partial_x)+[1+e_{\de}(-2x-ia_{-\de}+2\eta)]\exp(ia_{-\de}\partial_x),\label{cAT}
\end{gather}
can also be obtained as the limits
\begin{gather}
\cA^{T}_{\de}(\eta;x)= \lim_{\Lambda\to\infty} \cA_{\de}(a-i\eta-i\Lambda;x+\Lambda),
\end{gather}
cf.~\eqref{cAdef} and~\eqref{Vlim}. Note that they have holomorphic coef\/f\/icients, whereas the dual A$\De$Os
\begin{gather}
\hat{A}^T_{\de}(\eta;y)     \equiv    w^T(y)^{-1/2}\hat{H}^T_{\de}(\eta;y)w^T(y)^{1/2}
=    \frac{ie_{\de}(\eta)}{2s_{\de}(y)}\big(\exp(ia_{-\de}\partial_y)-\exp(-ia_{-\de}\partial_y)\big),\label{ATd}
\end{gather}
(which are the counterparts of $A_{\de}(b;y)$), have meromorphic coef\/f\/icients. Surprisingly, the A$\De$Os
$\cA_{\pm}(a-i\eta-i\Lambda;y)$ have no limit, whereas they do have obvious Toda counterparts, namely
\begin{gather}\label{cATd}
\hat{\cA}^T_{\de}(\eta;y)\equiv c^T(\eta;y)^{-1}\hat{A}^T_{\de}(\eta;y)c^T(\eta;y)=\exp(-ia_{-\de}\partial_y)
+\frac{e_{\de}(\eta)}{2s_{\de}(y)}\exp(ia_{-\de}\partial_y)\frac{e_{\de}(\eta)}{2s_{\de}(y)}.\!\!\!
\end{gather}

Even though $w(a-i\eta-i\Lambda;x+\Lambda)$ has no  limit either, there exists a function
\begin{gather}\label{wT}
w_T(\eta;x)\equiv 1/E(\pm (x-\eta)),
\end{gather}
 that may be viewed as a weight function. Here, $E(x)$ is the function  featuring in~\eqref{GE}, which we already had occasion to use in~Section~\ref{section3}, cf.~the paragraph containing~\eqref{P}.

To explain this interpretation, we f\/irst invoke~\eqref{EwT}: This representation makes clear that $w_T(\eta;x)$ is a real-analytic positive function on $\R$. Secondly, we note that when we set
 \begin{gather}
 c_T(\eta;x)\equiv E(x-\eta),
 \end{gather}
 then we get
 \begin{gather}
 u_T(\eta;x)\equiv c_T(\eta;x)/c_T(\eta;-x)=G(x-\eta).
 \end{gather}
 Hence, ignoring phases and quadratic exponentials, this $u$-function encodes the conjectured large-$y$ asymptotics~\eqref{yas}. Finally, the similarity-transformed A$\De$Os
 \begin{gather}
A_{\de}^T(x-\eta)\equiv w_T(\eta;x)^{-1/2}H_{\de}^T(\eta;x)w_T(\eta;x)^{1/2},
\end{gather}
have holomorphic coef\/f\/icients. Specif\/ically, from~\eqref{Eades} we compute
\begin{gather}\label{ATnew}
\frac{A_{\de}^T(x)}{\sqrt{2\pi}}
=    \frac{e_{\de}(-x-ia_{-\de}/2)\exp(-i(x+ia_{-\de}/2)K_{-\de})}{\Gamma(-i(x+ia_{-\de}/2)/a_{\de}+1/2)}\exp(ia_{-\de}\partial_x)+(i\to -i).
\end{gather}

Thus far, we have kept $x$ real and $y$ positive in the function F$^T(\eta;x,y)$. We proceed to study its analyticity features. To this end, consider the function
\begin{gather}\label{cH}
\cH(x-\eta,y)\equiv \sqrt{a_+a_-}w_T(\eta;x)^{-1/2}w^T(y)^{-1/2}{\rm F
}^T(\eta;x,y).
\end{gather}
The weight functions $w^T$ and $w_T$ (given by~\eqref{wTy} and~\eqref{wT}) are well understood from an analytic viewpoint, so we need only clarify the character of $\cH(x,y)$. The f\/irst point to note is that for each of the four integral representations~\eqref{T1}, \eqref{T2}, \eqref{T3} and~\eqref{T4} of~F$^T(\eta;x,y)$ the weight function factors on the right-hand side of~\eqref{cH} ensure that the prefactors of the $z$-integrals become entire functions of $x$ and $y$. Indeed, this is clear from the corresponding representations
\begin{gather}
\cH(x,y)     =    E(-x)\exp(3i\chi/2)\exp\left( \frac{i\alpha}{4}\big(y^2-iax-a^2/2\big)\right)
\nonumber \\
\phantom{\cH(x,y)     =}{}   \times \iR dz G(\pm z +x/2-ia/2)
 \exp (i\alpha[zy+z^2/2]),\label{cH1}
\\
\cH(x,y)    =   E(x) \exp(-3i\chi/2)\exp\left( \frac{-i\alpha}{4}\big(y^2+iax-a^2/2\big)\right)
\nonumber \\
\phantom{\cH(x,y)    =}{}  \times \iR dz G(\pm z -x/2-ia/2) \exp \big({-}i\alpha[zy+z^2/2]\big),\label{cH2}
\\
\cH(x,y)     =    E(-x)\exp(-i\chi/2)\exp(i\alpha (y^2-x^2)/8)
\nonumber \\
\phantom{\cH(x,y)     =}{}  \times \int_{\R +i0} dzG(z\pm y/2-ia)\exp\big({-}i\alpha[zx+z^2/2]\big),\label{cH3}
\\
\cH(x,y)    =   E(x)\exp(i\chi/2)\exp\big(i\alpha (x^2- y^2)/8\big)
\nonumber \\
\phantom{\cH(x,y)    =}{}  \times \int_{\R +i0} dzG(z\pm y/2-ia)\exp\big({-}i\alpha[zx-z^2/2]\big).\label{cH4}
\end{gather}

We have singled out the function~$\cH(x,y)$, because it extends from the real $x$-axis and the positive $y$-axis (where it takes real values) to a holomorphic (i.e., entire) function in $x$ and $y$. Taking this assertion for granted, the analytic character of F$^T(\eta;x,y)$ can be read of\/f from~\eqref{cH}. We continue by proving the holomorphy claim.

Consider f\/irst the function def\/ined by the integral in~\eqref{cH1}. The integrand is a meromorphic function~$I(z)$, whose asymptotics for $\re z\to\pm\infty$ readily follows from the $G$-asymptotics~\eqref{Gas}. Specif\/ically, we get
\begin{gather}
I(z)=O(\exp(-\alpha \re z[-\im x/2+a/2+\im z+\im y]),\qquad \re z\to\infty,
\\
I(z)=O(\exp(\alpha \re z[-\im x/2+a/2-\im z-\im y]),\qquad \re z\to -\infty.
\end{gather}
Therefore, exponential decay for $\re z\to\infty$ can be achieved by taking
\begin{gather}\label{tailr1}
\im z>\im x/2 -a/2-\im y,
\end{gather}
and for $\re z\to -\infty$ by taking
\begin{gather}
\im z<-\im x/2 +a/2-\im y.
\end{gather}
Since the two $G$-functions do not depend on $y$, this already implies that $\cH(x,y)$ extends to a holomorphic function of $y$. Indeed, when we continue $y$ of\/f the positive axis, we need only move the contour~$\R$ up on the right and down on the left (whenever need be) so as to retain exponential decay.

For the $x$-continuation there is also no problem coming from the tail ends of the contour, but we need to avoid that the contour gets pinched between the upward and downward pole sequences (cf.~\eqref{zkl}--\eqref{Gpo}),
\begin{gather}
z=-x/2-ia/2-z_{kl},\qquad z=x/2+ia/2+z_{kl},\qquad k,l\in\N,
\end{gather}
as $x$ is continued of\/f the real axis. This can be achieved by requiring
\begin{gather}\label{xres1}
x\notin -i[a,\infty).
\end{gather}
As a result, $\cH(x,y)$ extends to a holomorphic function of $x$ and $y$ outside the half line~\eqref{xres1}.

Turning to the representation~\eqref{cH2}, we can argue in the same way to conclude that
exponential decay for $\re z\to\infty$ can be achieved by taking
\begin{gather}\label{br}
\im z<\im x/2 +a/2-\im y,
\end{gather}
and for $\re z\to -\infty$ by taking
\begin{gather}\label{bl}
\im z>-\im x/2 -a/2-\im y.
\end{gather}
Here we get poles for
\begin{gather}\label{zxpo}
z=x/2-ia/2-z_{kl},\qquad z=-x/2+ia/2+z_{kl},\qquad k,l\in\N,
\end{gather}
as $x$ is continued of\/f the real axis. Thus we should require
\begin{gather}\label{xres2}
x\notin i[a,\infty),
\end{gather}
so as to avoid contour pinching. It therefore follows that~$\cH(x,y)$ extends to a holomorphic function outside the half line~\eqref{xres2}.

Combining these two conclusions, we deduce that $\cH(x,y)$ extends to a holomorphic function on~$\C^2$, as asserted. It also follows that the contour integral in~\eqref{cH1} extends to a meromorphic function of $x$ and $y$, with poles only at $x=-ia-z_{kl}$. Likewise, the contour integral in~\eqref{cH2} yields a meromorphic function with poles only at $x=ia+z_{kl}$.

To conclude this account of analyticity features, we point out that the holomorphy of $\cH(x,y)$ can also be derived in a somewhat dif\/ferent way from the two representations~\eqref{cH3}--\eqref{cH4}. For these cases the two downward pole sequences at $z=\mp y/2-z_{kl}$ can always be avoided by moving the contour up. Now, however, another type of restriction arises from the requirement of exponential decay on the contour tails. For the integrand in~\eqref{cH3} we need $\im z<a/2-\im x/2$ on the right tail, but to obtain exponential decay on the left tail we must require $\im x>-a$. Thus we can only deduce holomorphy for $\im x>-a$. Likewise, for~\eqref{cH4} we need~$\im x<a$ and~$\im z<a/2+\im x/2$ on the left tail, so we can only infer holomorphy for~$\im x<a$.

Even so, from these two f\/indings we can again conclude holomorphy on~$\C^2$.  Moreover, it follows that the contour integrals in~\eqref{cH3} and \eqref{cH4} give rise  to meromorphic functions of $x$ and $y$ with poles only at $x=-ia-z_{kl}$ and  $x=ia+z_{kl}$, respectively.

\subsection{Joint eigenfunction properties}\label{section5.3}

To complete this section, we verify that the joint eigenfunction properties have survived the \mbox{$\Lambda\to\infty$} limit. Due to the simpler analyticity properties of the pertinent contour integrals (compared to the hyperbolic case), this is rather straightforward. First, to show that F$^T(\eta;x,y)$ is an eigenfunction of the Hamiltonians~$H^T_{\pm}(\eta;x)$~\eqref{HT} with eigenvalues~$2c_{\pm}(y)$, we need only show that the A$\De$Os $\cA^T_{\pm}(0;x)$~\eqref{cAT} have the latter eigenvalues on the function~$G_R(x)^{-1/2}\rF^T (0;x,y)$.

{\sloppy To this end we invoke the representation~\eqref{T2eq}.  It follows from our analysis of the contour integral in~\eqref{cH2} that the contour integral in~\eqref{T2eq} with $\eta=0$ def\/ines a function $\cM(x,y)$ that extends to a meromorphic function of $x$ and $y$ with poles occurring solely for $x=ia+z_{kl}$. We may write this function for $x\notin i[a,\infty)$ as
\begin{gather}
\cM(x,y)=\int_{\cC}dzK_T(x,z)\exp(i\alpha zy),
\end{gather}
where we have introduced
\begin{gather}\label{KT}
K_T(x,z)\equiv G_L(\pm z-x/2-ia/2).
\end{gather}
Also, the choice of $\im z$ on the horizontal tails of the contour $\cC$ depends on $x$ and $y$ via \eqref{br}/\eqref{bl} on the right/left tail, while the tails are connected by a curve separating the upward and downward pole sequences~\eqref{zxpo}.

}

Now our task is to prove
\begin{gather}\label{xade}
\cA^T_{\de}(0;x)\cM(x,y)=2c_{\de}(y)\cM(x,y).
\end{gather}
By virtue of the already known meromorphy properties, it suf\/f\/ices to show~\eqref{xade} for $x$ varying over a rectangle $\re x \in (-a,a)$, $\im x \in (-4a,-3a)$ (say), while keeping $y$ positive. The restriction on $\im x$ ensures that the upward and downward $z$-pole sequences~\eqref{zxpo} are at a distance at least $2a$ from the real axis. Accordingly, we choose the contour $\cC$ to coincide with the real axis for $\re z\in [-2a,2a]$ (say). Furthermore, we need to  push down the right tail end and push up the left tail end suf\/f\/iciently far to retain exponential decay when we act with the $x$-shifts under the integral sign.

The crux is now that we have a kernel identity
\begin{gather}
K_T(x-ia_{-\de},z)+[1+e_{\de}(-2x-ia_{-\de})]K_T(x+ia_{-\de},z)\nonumber\\
\qquad{} =
K_T(x,z-ia_{-\de}/2)+K_T(x,z+ia_{-\de}/2).\label{Kid1}
\end{gather}
(This identity can be readily checked by dividing f\/irst by $K_T(x,z-ia_{-\de}/2)$ and then using the $G_L$-A$\De$Es~\eqref{GLDE}.) Thus we obtain
\begin{gather}
\cA^T_{\de}(0;x)\int_{\cC} dzK_T(x,z)e^{i\alpha zy}
=\int_{\cC} dz\big(K_T(x,z-ia_{-\de}/2)+K_T(x,z+ia_{-\de}/2)\big)e^{i\alpha zy}.
\end{gather}
Shifting $\cC$ up and down by $a_{-\de}/2$, no poles are met, and on the resulting contours $\cC_{+}$ and $\cC_{-}$ there is still exponential decay at the left and right. Hence we get
\begin{gather}
e_{\de}(-y)\int_{\cC_+}dzK_T(x,z)e^{i\alpha zy}+
e_{\de}(y)\int_{\cC_-}dzK_T(x,z)e^{i\alpha zy}.
\end{gather}
Since the two integrands are now equal and the integrals yield the same value $\cM(x,y)$, the eigenvalue property~\eqref{xade} follows.

We point out that the kernel identity~\eqref{Kid1} plays a role similar to the kernel identity~\eqref{dv} of the hyperbolic case. In the Toda case, however, we can ensure that the $z$-poles stay at an arbitrary distance from the real axis by choosing $\im x$ appropriately, by contrast to the $v$-poles in~\eqref{defB}, cf.~\eqref{Kpo}. Moreover, in the Toda case the decay properties on the horizontal tails of the contour depend on $\im z$, whereas the choice of $\im z$ is irrelevant in the hyperbolic case (since the `$z^2$-terms' drop out in the $|\re z|\to\infty$ asymptotics).

Next, we note that to prove that F$^T(\eta;x,y)$ is a joint eigenfunction of the dual Hamilto\-nians~$\hat{H}^T_{\pm}(\eta;y)$~\eqref{HTd} with eigenvalues~$e_{\pm}(x)$, we need only show that the A$\De$Os $\hat{A}^T_{\pm}(0;y)$~\eqref{ATd} have these eigenvalues on the function~$w^T(y)^{-1/2}\rF^T (0;x,y)$.

To this end we use the representation~\eqref{T3}. Accordingly we introduce
\begin{gather}\label{cMd}
\hat{\cM}(x,y)=\int_{\hat{\cC}}dz\hat{K}_T(x,z)\exp(-i\alpha zx),
\end{gather}
where
\begin{gather}\label{KTd}
\hat{K}_T(y,z)\equiv G( z\pm y/2-ia)\exp(i\alpha(y^2/8-z^2/2)).
\end{gather}
Recalling our analysis of the contour integral in~\eqref{cH3}, we see that we should require f\/irst of all $\im x>-a$ in~\eqref{cMd}. Then the integral is well def\/ined when we choose the horizontal tails of~$\hat{\cC}$ equal to $\R$ (say) on the left and having $\im z<a/2-\im x/2$ on the right, while the middle part is above the pole sequences at $z=\pm y/2-z_{kl}$. Furthermore, the function $\hat{\cM}(x,y)$ is holomorphic in $x$ and $y$ for $\im x>-a$ and extends to a meromorphic function with poles at $x=-ia-z_{kl}$.

In view of these features, we need only prove
\begin{gather}\label{ATM}
\hat{A}^T_{\de}(0;y)\hat{\cM}(x,y)=e_{\de}(x)\hat{\cM}(x,y),
\end{gather}
for $y$ varying over a square $\re y \in (-a,a)$, $\im y \in (-a,a)$ (say), while keeping $x$ real. To do so, we choose the middle part of the contour $\hat{\cC}$ equal to $2ia+(-2a,2a)$,  and connect this part to $(-\infty,-3a)$ and $(3a,\infty)$ in the obvious way. Then the $y$-shifts can be taken under the integral sign without $z$-poles hitting the contour.

With these analytic preliminaries in place, the key algebraic point is once more a kernel identity, namely,
\begin{gather}\label{Kdid}
\frac{i}{2s_{\de}(y)}\big(\hat{K}_T(y+ia_{-\de},z)-\hat{K}_T(y-ia_{-\de},z)\big)=\hat{K}_T(y,z-ia_{-\de}/2).
\end{gather}
(To check it, one need only divide by the r.h.s.\ and use the $G$-A$\De$Es~\eqref{Gades}.) Thus the l.h.s.\ of~\eqref{ATM}
equals
\begin{gather}
\int_{\hat{\cC}}dz\hat{K}_T(x,z-ia_{-\de}/2)\exp(-2i\pi zx/a_+a_-).
\end{gather}
Shifting the contour
up by $a_{-\de}/2$, no poles are met, and so the joint eigenvalue equations~\eqref{ATM} follow.

{\sloppy We conclude this subsection with some remarks.  It follows from the eigenvalue features just proved that the holomorphic function $\cH(x,y)$ is a joint eigenfunction of the four A$\De$Os $A^T_{\pm}(x)$~\eqref{ATnew} and $\hat{A}^T_{\pm}(0;y)$~\eqref{ATd} with eigenvalues $2c_{\pm}(y)$ and $e_{\pm}(x)$, respectively. The coef\/f\/icients of the former A$\De$Os are entire, whereas the coef\/f\/icients of the latter are meromorphic. The coef\/f\/icients of the A$\De$Os $\cA^T_{\pm}(0;x)$ \eqref{cAT} are entire as well, but their joint eigenfunction $G_R(x)^{-1/2}\rF^T(0;x,y)$ is meromorphic in~$x$, with poles for $x=ia+z_{kl}$.

}

This shows by example that the meromorphic vs.~entire character of the coef\/f\/icients of the type of commuting A$\De$O pairs at issue in this paper is compatible both with entire and with meromorphic joint eigenfunctions. In this connection, it should be noted that when the ratio $a_+/a_-$ is irrational, then the only multipliers that do no destroy the joint eigenfunction property are the constants.

Another consequence worth pointing out consists of the relations
\begin{gather}
\cH(x,ka_{\de}+ia_{-\de})=\cH(x,ka_{\de}-ia_{-\de}),\qquad \forall \, (x,k,\de)\in\C\times\Z\times\{ +,-\}.
\end{gather}
Indeed, these follow from the eigenvalue equations
\begin{gather}
\cH(x,y-ia_{-\de})-\cH(x,y+ia_{-\de})=2is_{\de}(y)e_{\de}(x)\cH(x,y).
\end{gather}

\section{The nonrelativistic Toda case}\label{section6}

In this section we obtain the nonrelativistic counterparts of the quantities in Section~\ref{section5}, along the lines laid out in Subsection~\ref{section4.2} for the hyperbolic case. Thus we switch to the parameters~\eqref{repar} and momentum variable~\eqref{pdef}, whereas the Toda analog of~\eqref{bg} is the substitution
\begin{gather}
\eta=\frac{2}{\mu}\ln (\beta\mu g),\qquad g>0.
\end{gather}

For the operators~$H^T_+(\eta;x)$~\eqref{HT} and $\cA^T_+(\eta;x)$~\eqref{cAT}, these substitutions entail the expansion
\begin{gather}\label{HcAp}
H^T_+,\cA^T_+= 2+ \beta^2 H +O\big(\beta^4\big),\qquad \beta \to 0,
\end{gather}
where
\begin{gather}
H\equiv -\hbar^2 \partial_x^2+ \mu^2g^2 \exp(-\mu x),
\end{gather}
is the nonrelativistic Toda Hamiltonian. Moreover, we clearly get
\begin{gather}\label{HcAm}
\lim_{\beta \to 0}H^T_-=\lim_{\beta \to 0}\cA^T_-=
\exp\big(2i\pi\mu^{-1}\partial_x\big) +(i\to -i)=:\cM.
\end{gather}
In view of \eqref{cpexp} the eigenvalue of the Hamiltonian $H$ becomes $p^2/4$, while the eigenvalue of the monodromy operator $\cM$ remains $2\cosh(\pi p/\hbar\mu)$.

Turning to the dual operators, for $\hat{H}^T_+(\eta;y)$~\eqref{HTd}, $\hat{A}^T_+(\eta;y)$~\eqref{ATd} and~$\hat{\cA}^T_+(\eta;y)$~\eqref{cATd} we obtain (with $p>0$)
\begin{gather}
\lim_{\beta\to 0}\hat{H}^T_+=\mu gp^{-1/2}\big(\exp(i\hbar \mu\partial_p)+(i\to -i)\big)p^{-1/2}=: \hat{H},
\\
\label{hatA}
\lim_{\beta\to 0}\hat{A}^T_+=i\mu gp^{-1}\big(\exp(i\hbar \mu\partial_p)-(i\to -i)\big)=:\hat{A},
\\
\lim_{\beta\to 0}\hat{\cA}^T_+=\exp(-i\hbar \mu\partial_p)+
\frac{\mu g}{p}
\exp(i\hbar \mu\partial_p)\frac{\mu g}{p}=:\hat{\cA},
\end{gather}
and their eigenvalue becomes $\exp(\mu x/2)$. The dual operators $\hat{H}^T_-$, $\hat{A}^T_-$ and $\hat{\cA}^T_-$ have no limits, however. The operators $\hat{A}$ and $\hat{H}$ are related by a similarity transformation with the square root of the limit function
\begin{gather}\label{limwT}
 \lim_{\beta \to 0}(\hbar\beta\mu)^{-1}w^T(2\pi/\mu,\hbar\beta;\beta p/\mu)= \hat{w}_{\rm nr}(p/\hbar\mu),
 \end{gather}
 where
  \begin{gather}
 \hat{w}_{\rm nr}(k)\equiv 2k\sinh (\pi k)=2\pi/\Gamma(\pm ik),
  \end{gather}
and the operators $\hat{A}$ and $\hat{\cA}$ by similarity with
\begin{gather}\label{limcT}
\lim_{\beta \to 0}(\hbar\beta\mu)^{1/2}c^T(2\pi/\mu,\hbar\beta,2\mu^{-1}\ln(\beta\mu g);\beta p/\mu)= \hat{c}_{\rm nr}(g/\hbar;p/\hbar\mu),
 \end{gather}
where
\begin{gather}
 \hat{c}_{\rm nr}(\lambda;k)\equiv (2\pi)^{-1/2}\exp(-ik\ln \lambda)\Gamma(ik).
   \end{gather}
(The limit~\eqref{limwT} is clear from~\eqref{wTy}, whereas~\eqref{limcT} follows by using~\eqref{GhGr}.)

We proceed to obtain the nonrelativistic limit of the joint eigenfunction $\rF^T(\eta;x,y)$ for each of the representations~\eqref{T1eq}, \eqref{T2eq}, \eqref{T3} and~
\eqref{T4}, taking $x$ real and $y$ positive.

 First, we have
\begin{gather}\label{GRlim}
G_R(x-\eta)\to G_R\left(\frac{2\pi}{\mu},\hbar\beta;x+\frac{2}{\mu}\ln\left(\frac{\hbar}{\mu g}\right)
+\frac{2}{\mu}\ln\left(\frac{1}{\hbar\beta}\right)\right).
\end{gather}
Comparing this to the limit~\eqref{gLRLimit}, we see that it applies to the $\beta\to 0$ limit of~\eqref{GRlim}, with the parameter $\lambda$ equal to 2. Therefore, the limit of $G_R(x-\eta)$ equals 1, a circumstance that also explains why we get coinciding limits for $H_{\de}^T$ and $\cA_{\de}^T$ in~\eqref{HcAp}--\eqref{HcAm}, cf.~\eqref{cAT}.

Next, from~\eqref{limwT} we see that the above substitutions imply
\begin{gather}\label{prelim}
\lim_{\beta\to 0}\frac{G(\pm y+ia)}{a_+a_-}=\frac{\mu p}{\pi\hbar}\sinh(\pi p/\hbar\mu).
\end{gather}
Moreover, $\alpha y^2$ clearly vanishes for $\beta\to 0$, so it now follows that the limits of the prefactors of the four representations are all equal to the square root of the right-hand side of~\eqref{prelim}.

Turning to the integrand in~\eqref{T1eq}, the plane wave $\exp(i\alpha zy)$ becomes $\exp(izp/\hbar)$. Furthermore, the substitutions on the two $G_R$-functions yield
\begin{gather}\label{GR2lim}
G_R\left(\frac{2\pi}{\mu},\hbar\beta;\pm z +\frac{x}{2}-\frac{i\pi}{2\mu}-\frac{i\hbar\beta}{4}+\frac{1}{\mu}\ln\left(\frac{\hbar}{\mu g}\right)
+\frac{1}{\mu}\ln\left(\frac{1}{\hbar\beta}\right)\right).
\end{gather}
Comparing once more to the limit~\eqref{gLRLimit}, we see that it now applies to the $\beta\to 0$ limit of~\eqref{GR2lim} with $\lambda$ equal to 1. Therefore, a short calculation gives the limit
\begin{gather}
\exp\big({-}2g\hbar^{-1}\exp(-\mu x/2)\cosh (\mu z)\big).
\end{gather}
Taking $z\to t/\mu$ in the integral, we wind up with the limit function
\begin{gather}\label{FnrT}
\rF_{\rm nr}^T(g/\hbar;\mu x/2,p/\hbar\mu),
\end{gather}
where
\begin{gather}\label{FnrT1}
\rF_{\rm nr}^T(\lambda;r,k)=2\big(\pi^{-1}k\sinh \pi k\big)^{1/2}
\int_0^{\infty} dt \exp(-2\lambda e^{-r}\cosh t)\cos(tk).
\end{gather}

It is readily verif\/ied that~\eqref{T2eq} also leads to the representation~\eqref{FnrT1} for the limit function. Proceeding with the above substitutions for~\eqref{T3}, we see that we can only get convergence of the exponentials for $\beta \to 0$ when we f\/irst replace the integration variable $z$ by $\hbar \beta w$ (say). Hence we get a factor $\hbar\beta$ up front, and the plane wave $\exp(-i\alpha zx)$ becomes $\exp(-i\mu wx)$; moreover, the quadratic exponential converges to 1 for $\beta\to 0$. It therefore remains to consider the ratio
\begin{gather}\label{intT3}
\hbar\beta\exp(2iw\ln (\beta\mu g))/G\left(\frac{2\pi}{\mu},\hbar\beta;-\hbar\beta w\pm \frac{\beta p}{2\mu}+\frac{i}{2}\left(\frac{2\pi}{\mu}+\hbar\beta\right)\right).
\end{gather}
Scaling the $G$-functions by $\mu/2\pi$, we can invoke~\eqref{GhGr} to deduce that~\eqref{intT3} has $\beta\to 0$ limit
\begin{gather}
(2\pi\mu)^{-1}\exp(2iw\ln (g/\hbar))\Gamma(-iw\pm ip/2\hbar\mu).
\end{gather}
Putting the pieces together, we now get the limit function~\eqref{FnrT} represented as
\begin{gather}\label{FnrT2}
\rF_{\rm nr}^T(\lambda;r,k)=\frac{1}{4\pi}\big(\pi^{-1}k\sinh \pi k\big)^{1/2}
\int_{\R+i0}dt\Gamma(i(-t\pm k)/2)\exp(it(\ln\lambda -r)).
\end{gather}

The representation~\eqref{T4} also leads to~\eqref{FnrT2}. Thus we obtain two dif\/ferent representations for the limit function~\eqref{FnrT}. The resulting identity
\begin{gather}
\int_0^{\infty} dt \exp\big({-}2 e^{-v}\cosh t\big)\cos(tk)=\frac{1}{8\pi}\int_{\R+i0}dt\Gamma(i(-t\pm k)/2)\exp(-itv),
\end{gather}
is known and not hard to verify. Indeed, the eigenfunction transform associated with~\eqref{FnrT} amounts to the Kontorovich--Lebedev transform, and the integrals yield distinct representations of the modif\/ied Bessel function $K_{ik}(2e^{-v})$, cf.\ e.g.~\cite[10.32.9, 10.32.13 and 10.43(v)]{dlmf}.

From~\eqref{FnrT1} it is obvious that $\hat{w}_{\rm nr}(k)^{-1/2}\rF_{\rm nr}^T(\lambda;r,k)$ extends from the positive $k$-axis to an entire function of~$k$. Neither from~\eqref{FnrT1} nor from~\eqref{FnrT2} entireness in~$r$ is manifest, but this is well known (it can be inferred, e.g., from ODE theory). As already mentioned in the Introduction, the eigenfunction property for the dual operators seems not to occur in the standard sources. It is most easily checked for~$\hat{A}$~\eqref{hatA} by proceeding as in the relativistic case. Here it follows from the kernel identity
\begin{gather}
ik^{-1}\big(\hat{K}_{\rm nr}(k+i,t)-\hat{K}_{\rm nr}(k-i,t)\big)
=\hat{K}_{\rm nr}(k,t-i),
\end{gather}
where
\begin{gather}
\hat{K}_{\rm nr}(k,t)\equiv \Gamma (-it/2\pm ik/2).
\end{gather}
(This identity is the counterpart of~\eqref{Kdid}.)

The large-$r$ asymptotics of $\rF_{\rm nr}^T(\lambda;r,k)$ is well known. It is given by
\begin{gather}
\rF_{\rm nr}^T(\lambda;r,k)\sim \hat{u}_{\rm nr}(\lambda;k)^{1/2}e^{irk}+\hat{u}_{\rm nr}(\lambda;-k)^{1/2}e^{-irk},\qquad r\to\infty,
\end{gather}
where
\begin{gather}
\hat{u}_{\rm nr}(\lambda;k)=\hat{c}_{\rm nr}(\lambda;k)/\hat{c}_{\rm nr}(\lambda;-k),
\end{gather}
and readily verif\/ied from~\eqref{FnrT2}.
It seems much harder to obtain the large-$k$ asymptotics from the above representations. Assuming plane-wave behavior, a consideration of the dual operators leads to the expectation
\begin{gather}
\rF_{\rm nr}^T(\lambda;r,k)\sim \exp( i\phi +ik\ln k -ik)e^{i(r-\ln \lambda)k} \nonumber\\
\phantom{\rF_{\rm nr}^T(\lambda;r,k)\sim}{} +\exp(- i\phi -ik\ln k +ik)e^{-i(r-\ln \lambda)k},\qquad  k\to\infty,\label{kas}
\end{gather}
where $\phi\in[-\pi,\pi)$. Indeed, just like for its relativistic counterpart~\eqref{yas}, this seems the simplest behavior that is consistent with the eigenvalues and dual eigenvalues. To be sure,  for neither case it is a priori clear that the asymptotic behavior must involve plane waves.

At any rate, a result pertinent to~\eqref{kas} can be found in the literature: An asymptotic expansion for $K_{ip}(x)$ with $p>x>0$ occurs in~\cite[Section~7.13.2, formula~(19)]{erde}. The dominant asymptotics does give rise to~\eqref{kas} with $\phi=-\pi/4$, but an unsettling $O(x^{-1})$ error term is present. (If dependence on~$x$ is included, then one would rather expect an~$O(x)$ error term. Indeed, $x\to 0$ corresponds to~$v\to \infty$, a limit for which the Toda potential is exponentially vanishing.)

\appendix


																																												 
\section{The hyperbolic gamma function}\label{appendixA}

The  hyperbolic gamma function was introduced and studied in~\cite{fo} as a so-called minimal solution of a special f\/irst order analytic dif\/ference equation. It is basically the same as Kurokawa's double sine~\cite{kuro}, Faddeev's quantum dilogarithm~\cite{fadd2}, and Woronowicz's quantum exponential function~\cite{woro}. (The precise connections between these functions are spelled out in Appendix~A of our paper~\cite{Hels}.) In this appendix we review features of the hyperbolic gamma function $G(\aaa;z)$ that are used in the present paper; if need be, see~\cite{fo} for proofs.

Unless specif\/ied otherwise, we choose
\begin{gather}
\aaa >0,
\end{gather}
and suppress the dependence of $G$ on $a_+$, $a_-$.

 To begin with, $G(z)$ can be def\/ined as the unique minimal solution of one of the two analytic dif\/ference equations
\begin{equation}\label{Gades}
	\frac{G(z+ia_\delta/2)}{G(z-ia_\delta/2)} = 2c_{-\delta}(z),\qquad  \delta=+,-,
\end{equation}
that has modulus 1 for real $z$ and satisf\/ies $G(0)=1$ (recall~\eqref{denot} for the notation used here); remarkably, this entails that the other one is then satisf\/ied as well.
It is meromorphic in~$z$,
and for $z$ in the strip
\begin{gather}\label{strip}
S\equiv \{z\in\C \mid |\im (z)|<a\},
\end{gather}
no poles and zeros occur. Hence we have
\begin{gather}\label{Gg}
G(z)=\exp(ig(z)),\qquad z\in S,
\end{gather}
with the function $g(z)$ being holomorphic in $S$.
 Explicitly, $g(z)$ has the integral representation
\begin{gather}\label{ghyp}
g(a_{+},a_{-};z) =\int_0^\infty\frac{dy}{y}\left(\frac{\sin 2yz}{2\sinh(a_{+}y)\sinh(a_{-}y)} - \frac{z}{a_{+}a_{-} y}\right),\qquad z\in S.
\end{gather}
From this, the following properties of the hyperbolic gamma function are immediate:
\begin{gather}\label{refl}	
G(-z) = 1/G(z),\qquad ({\rm ref\/lection\ equation}),
\\
\label{modinv}
G(a_-,a_+;z) = G(a_+,a_-;z),\qquad  ({\rm modular\ invariance}),
\\
\label{sc}
	G(\lambda a_+,\lambda a_-;\lambda z) = G(a_+,a_-;z),\qquad \lambda\in(0,\infty),\qquad ({\rm scale\ invariance}),
\\
\label{Gcon}
\overline{G(\aaa;z)}=G(\aaa;-\overline{z}).
\end{gather}

We have occasion to use a few more features that are less obvious, including the duplication formula
\begin{gather}\label{dupl}
G(\aaa;2z)=G(\aaa;z\pm ia_{+}/4 \pm ia_{-}/4),
\end{gather}
the closely related formula
\begin{gather}\label{G2}
G(a_{+},2a_{-};2z)=G(\aaa;z\pm ia_{+}/4),
\end{gather}
the explicit evaluation
\begin{gather}\label{Geval}
G(ia_{+}/2-ia_{-}/2)=(a_{+}/a_{-})^{1/2},
\end{gather}
and the asymptotic behavior of $G(z)$ for $\re (z)\to\pm \infty$. The latter is given by
\begin{gather}\label{Gas}
	G(\aaa;z) = \exp \big({\mp} i\left(\chi+\alpha z^2/4\right)\big)\big(1 + O(\exp(-r |\re(z)|))\big),\qquad \re(z)\to\pm\infty,
\end{gather}
where the decay rate $r$ can be any positive number satisfying
\begin{gather}
r <\alpha \min(a_+,a_-),
\end{gather}
and where
\begin{equation}\label{chi}
 \chi \equiv \frac{\pi}{24}\left(\frac{a_+}{a_-} + \frac{a_-}{a_+}\right).
\end{equation}

Def\/ining
\begin{gather}\label{zkl}
z_{kl}\equiv ika_{+}+ila_{-},\qquad k,l\in\N\equiv \{ 0,1,2,\ldots\},
\end{gather}
the hyperbolic gamma function has its poles at
\begin{gather}\label{Gpo}
z=z_{kl}^-,\qquad z_{kl}^-\equiv -ia -z_{kl},\qquad k,l\in\N,\qquad (G\mbox{-poles}),
\end{gather}
and its zeros at
\begin{gather}\label{Gze}
z=z_{kl}^+,\qquad z_{kl}^+\equiv ia +z_{kl},\qquad k,l\in\N,\qquad (G\mbox{-zeros}).
\end{gather}
The pole at $-ia$ is simple and has residue
\begin{gather}\label{Gres}
\lim_{z\to -ia}(z+ia)G(z)=\frac{i}{2\pi}(a_{+}a_{-})^{1/2}.
\end{gather}
In view of these features, $G(z)$ can be written as a ratio of entire functions,
\begin{gather}\label{GE}
G(\aaa;z)=E(\aaa;z)/E(\aaa;-z),
\end{gather}
where $E(\aaa;z)$ has its zeros at
\begin{gather}\label{Eze}
z=z^{+}_{kl},\qquad k,l\in\N,\qquad (E\mbox{-zeros}).
\end{gather}

The function $E(\aaa;z)$ we have occasion to employ is def\/ined in Appendix~A of~I; it is closely related to Barnes' double gamma function~\cite{barn}. We need two more of its properties. First, from equations~(A.41) and~(A.43) in~I we have
 \begin{gather}\label{EwT}
 E(z)E(-z)=\exp\left(\frac12\int_0^{\infty}\!\frac{dy}{y}\left(
 \frac{1-\cos(2yz)}{\sinh(a_{+}y)\sinh(a_{-}y)}-\frac{z^2}{a_+a_-}\big(e^{-2a_{+}y}+e^{-2a_{-}y}\big)\right)\right),\!\!\!
 \end{gather}
where $z$ belongs to the strip $S$~\eqref{strip}. Second, we need the A$\De$Es it satisf\/ies, namely,
\begin{gather}\label{Eades}
\frac{E(z+ia_{\de}/2)}{E(z-ia_{\de}/2)}=\sqrt{2\pi}\exp(izK_{\de})/\Gamma(iz/a_{-\de}+1/2),\qquad \de=+,-,
\\
K_{\de}\equiv \frac{1}{2a_{-\de}}\ln \left( \frac{a_{\de}}{a_{-\de}}\right),
\end{gather}
cf.~equations~(A.46)--(A.47) in~I.

We also state two zero step size limits of the hyperbolic gamma function, which we need for taking nonrelativistic limits. The f\/irst one yields the relation to the Euler gamma function:
\begin{gather}\label{GhGr}
\lim_{\kappa\downarrow 0}G(1,\kappa;
 \kappa z+i/2)\exp \big(iz\ln(2\pi\kappa)-\ln(2\pi)/2\big) = 1/\Gamma(iz+1/2).
\end{gather}
For the second one we need to require that $z$ stay away from  cuts given by $\pm i[a_+/2,\infty)$. Then we have
\begin{gather}\label{Gnr}
\lim_{a_{-}\downarrow 0}\frac{G(a_+, a_{-};z+iua_{-})}{G(a_+, a_{-};z+ida_{-})}
=\exp((u-d)\ln [2c_+(z)]),\qquad u,d\in\R,
\end{gather}
uniformly on compact subsets of the cut plane.

For the relativistic Toda setting it is expedient to employ two slightly dif\/ferent hyperbolic gamma functions def\/ined by
\begin{gather}\label{GRDef}
	G_R(z) \equiv \exp(ig_R(z)),\qquad g_R(z)\equiv g(z)+\chi + \alpha z^2/4,
\\
\label{GLDef}
	G_L(z) \equiv \exp(ig_L(z)),\qquad  g_L(z)\equiv g(z)-\chi - \alpha z^2/4.
\end{gather}
These functions are the unique minimal solutions of the analytic dif\/ference equations
\begin{gather}\label{GRDE}
	\frac{G_R(z+ia_{-\delta}/2)}{G_R(z-ia_{-\delta}/2)} = 1+e_\delta(-2z),\qquad \de=+,-,
\\
\label{GLDE}
	\frac{G_L(z+ia_{-\delta}/2)}{G_L(z-ia_{-\delta}/2)} = 1+e_\delta(2z),\qquad \de=+,-,
\end{gather}
with asymptotic behavior
\begin{gather}\label{GRLas}
	G_{\substack{R \\ L}}(z) = 1 + O\left(\exp(-r |\re(z)|)\right),\qquad \re(z)\to\pm\infty.
	\end{gather}
Furthermore, they are related by
\begin{gather}\label{GRGL}
G_R(z)G_L(-z)=1.
\end{gather}
 The properties of the functions $G_R$ and $G_L$ just stated are easy to infer from the corresponding properties of the hyperbolic gamma function. (In Appendix~A of~\cite{Hels} we already introduced functions $S_R$ and $S_L$ that dif\/fer from $G_R$ and $G_L$ by the shift $z\to z-ia$.)

 Finally, we have occasion to use  the limits
\begin{gather}\label{gLRLimit}
	\lim_{a_-\downarrow 0}g_{\substack{R\\L}}\left(a_{+},a_{-};z\pm \lambda s(\aaa)\right) =
	\left\lbrace \begin{array}{ll}	\pm\frac{a_+}{2\pi}e_+(\mp 2z), & \lambda=1,\\  0 , & \lambda>1,\end{array}\right.
	\end{gather}
where
\begin{gather}\label{defs}
s(\aaa)\equiv \frac{a_+}{2\pi}\ln\frac{1}{a_-},	
\end{gather}
which hold uniformly for $z$ varying over arbitrary compact subsets of $\C$. To our  knowledge, these limits have not been obtained before. We present their proof in the next appendix.



\section{Proof of~(\ref{gLRLimit})}\label{appendixB}

Since we have
\begin{gather}
g_L(z)=-g_R(-z),
\end{gather}
we need only show~\eqref{gLRLimit} for $g_R$. Our proof actually yields a stronger result, which may be useful in other contexts. To state this result, we f\/ix
\begin{gather}
\lambda \in \big(1/\sqrt{2},\infty\big),
\end{gather}
and choose $\de$ satisfying
\begin{gather}
\de\in\big(1,\sqrt{2}\big),\qquad \de\lambda>1.
\end{gather}
Then we shall show
\begin{gather}\label{gRb}
g_R(\aaa;z+\lambda s(\aaa))=\frac{a_-^{\lambda }}{2\sin(\pi a_-/a_+)}e_+(-2z)+O\big(a_-^{\de \lambda -1}\big),\qquad
 a_-\downarrow 0,
\end{gather}
where the implied constant can be chosen uniformly for $z$ varying over compact subsets of $\C$.
A key ingredient of our proof is the comparison function we used in~Subsection~III~A of~\cite{fo} to obtain the $G$-asymptotics~\eqref{Gas}. Specif\/ically, we f\/irst focus on the dif\/ference
\begin{gather}\label{defd}
d(\aaa;z)\equiv a_+a_-g_R(\aaa;z)-A^2g_R(A,A;z),\qquad z\in S,
\end{gather}
 for the special $A$-choice
\begin{gather}\label{Ac}
A\equiv \left(\frac{a_+^2+a_-^2}{2}\right)^{1/2}.
\end{gather}
Observing that $A\ge a$ (with equality for $a_+=a_-$), we deduce that $d(z)$ is well def\/ined and holomorphic for $z$ in the strip $S$~\eqref{strip}.
Recalling~\eqref{aconv}, \eqref{chi} and~\eqref{GRDef}, we see that we may rewrite $d(z)$ as
\begin{gather}
d(\aaa;z)=\frac{1}{2i}\int_{\R}dyI(\aaa;y)\exp(2iyz),
\end{gather}
where
\begin{gather}\label{defint}
I(\aaa;y)\equiv \frac{1}{2y}\left(\frac{a_+a_-}{\sinh(a_+ y)\sinh(a_- y)}-\frac{A^2}{\sinh^2(Ay)}\right).
\end{gather}
The $A$-choice \eqref{Ac} is the unique one guaranteeing that $I(y)$ has no pole at $y=0$. Specif\/ically, we easily calculate
\begin{gather}
I(y)=c(\aaa)y+O\big(y^3\big),\qquad y\to 0,
\end{gather}
where $c(\aaa)$ is a polynomial in $a_+$ and $a_-$ of degree~4.

Since we let $a_-$ go to 0, we may and will assume from now on
\begin{gather}\label{Aass}
 \frac{1}{A} > \de \frac{1}{a_+}.
\end{gather}
Next, we shift the $y$-contour up by
\begin{gather}\label{defr}
r\equiv \pi\de/a_+.
\end{gather}
On account of~\eqref{Aass}, this ensures that only the simple pole at $y=i\pi/a_+$ is passed. The residue at this pole is readily calculated, yielding the representation
\begin{gather}
d(z)=\frac{a_+a_-}{2\sin(\pi a_-/a_+)}e_+(-2z)+\rho(z), \qquad z\in S,
\end{gather}
where $\rho$ is the remainder integral
\begin{gather}
\rho(z)\equiv \frac{1}{2i}\exp(-2rz)\int_{\R}duI(u+ir)\exp(2iuz),\qquad r=\pi\de/a_+,\qquad z\in S.
\end{gather}

We are now prepared to replace $z$ by
\begin{gather}
z+\lambda s(\aaa),\qquad |\im z|\le a_+/2,
\end{gather}
so that we get, using~\eqref{defr},
\begin{gather}
d(z+\lambda s)=\frac{a_+a_-a_-^{\lambda}}{2\sin (\pi a_-/a_+)}e_+(-2z)\nonumber\\
\phantom{d(z+\lambda s)=}{}
+\frac{1}{2i}a_-^{\de\lambda}e_+(-2\de z)\int_{\R}duI(u+i\pi \de/a_+)\exp(2iu(z+\lambda s)).\label{dla}
\end{gather}
Next, we note
\begin{gather}
\lim_{a_-\downarrow 0}I(w)=\frac{1}{2w}\left( \frac{a_+}{w\sinh(a_+w)}-\frac{a^2_+/2}{\sinh^2(a_+w/\sqrt{2})}\right).
\end{gather}
Since $z$ is required to satisfy $|\im z|\le a_+/2$, it follows that the $u$-integrand in~\eqref{dla} remains bounded by a f\/ixed $L^1(\R)$-function as $a_-\downarrow 0$. Thus we readily deduce the bound{\samepage
\begin{gather}\label{dlim}
\frac{d(z+\lambda s)}{a_+a_-}=\frac{a_-^{\lambda }}{2\sin(\pi a_-/a_+)}
e_+(- 2z)+O\big(a_-^{\de\lambda -1}\big),\qquad |\im z|\le a_+/2,\qquad
 a_-\downarrow 0,
\end{gather}
with the implied constant uniform on compact subsets of the strip $|\im z|\le a_+/2$.}

Next, we claim that we have
\begin{gather}\label{gRsp}
\frac{A^2}{a_+a_-}g_R(A,A;z+\lambda s)=O(a_-^{\de\lambda -1}),\qquad |\im z|\le a_+/2,\qquad
 a_-\downarrow 0,
\end{gather}
 uniformly on compacts of $|\im z|\le a_+/2$. Taking this claim for granted, we see from~\eqref{defd} and \eqref{dlim} that \eqref{gRb} holds true, uniformly on compacts of the latter strip. Now the A$\De$E~\eqref{GRDE} with $\de=-$ implies
\begin{gather}
g_R(\aaa; z+ia_+/2 +\lambda s)-g_R(\aaa; z-ia_+/2 +\lambda s)\nonumber\\
\qquad{}=-i\ln \left(1+\exp\left(-\frac{2\pi}{a_-}(z+\lambda s)\right)\right).
\end{gather}
Hence we have (using the def\/inition~\eqref{defs} of $s$)
\begin{gather}
g_R(z+ia_+/2 +\lambda s)-g_R( z-ia_+/2 +\lambda s)=
O\left(\exp\left(-\frac{a_+}{a_-}\ln \frac{1}{a_-}\right)\right),\qquad a_-\downarrow 0,
\end{gather}
uniformly on compacts of $\C$. From this it is routine to infer that~\eqref{gRb} also holds for $z\in\C$, with the bound uniform for $z$ varying over arbitrary $\C$-compacts.

It remains to prove the claim. To this end we start from the identity
\begin{gather}
g_R(A,A;z)=\frac{1}{\pi}b_+(\pi z/A),\qquad b_+(w)\equiv\int_{w}^{\infty}dt \frac{te^{-t}}{\sinh(t)},
\end{gather}
which follows from equations~(3.41)--(3.46) in~\cite{fo}. Next, we note the bound
\begin{gather}
b_+(w+R)=R\exp(-2(w+R))(1+O(1/R)),\qquad R\to\infty,
\end{gather}
where the implied constant can be chosen uniform on $\C$-compacts. (One can use the elementary integral
\begin{gather}
\int_x^{\infty}te^{-ct}dt=\frac{1}{c^2}e^{-cx}(1+cx),\qquad c>0,
\end{gather}
to verify this estimate.) As a consequence, we obtain
\begin{gather}
\frac{A^2}{a_+a_-}g_R(A,A;z+\lambda s)\nonumber\\
\qquad{} \sim \frac{\lambda}{2\pi\sqrt{2}}a_-^{\lambda\sqrt{2}-1}\ln\left(\frac{1}{a_-}\right)\exp(-2\pi\sqrt{2}z/a_+),\qquad z\in\C, \qquad a_-\downarrow 0,
\end{gather}
uniformly for $z$ in $\C$-compacts. Thus (a stronger version of) our claim follows. This concludes the proof of~\eqref{gRb}, and so \eqref{gLRLimit} follows as an obvious corollary.


																																												 

\section{Fourier transform formulas}\label{appendixC}

We f\/ix complex numbers $\mu,\nu$ satisfying
\begin{gather}\label{munu}
-a<\im \mu<\im \nu <a,\qquad a=(a_{+}+a_{-})/2.
\end{gather}
Hence we have
\begin{gather}
\im (\nu-\mu)\in(0,a_{+}+a_{-}),
\end{gather}
and the function
\begin{gather}\label{defI}
I(\mu,\nu;x)\equiv G(x-\nu)/G(x-\mu),
\end{gather}
is pole-free in the strip
\begin{gather}\label{xstrip}
\im \nu -a<\im x<\im \mu +a.
\end{gather}
Also, from the $G$-asymptotics~\eqref{Gas} we deduce
\begin{gather}
I(\mu,\nu;x)=O(\exp(\mp \alpha\im (\nu -\mu)\re x/2)),\qquad \re x\to\pm\infty.
\end{gather}
Therefore, for real $y$ the function
\begin{gather}\label{defFy}
F(\mu,\nu;y)\equiv \int_{\R}dx \exp(i\alpha xy)\frac{G(x-\nu)}{G(x-\mu)},
\end{gather}
is well def\/ined, and analytic in $\mu$ and $\nu$ in the region~\eqref{munu}. Moreover, we retain exponential decay of the integrand when we let $y$ vary over the strip
\begin{gather}\label{ystrip}
|\im y|<\im(\nu-\mu)/2,
\end{gather}
so $F(\mu,\nu;y)$ is analytic in $y$ in this strip and given by~\eqref{defFy}.

Our f\/irst aim is to obtain the Fourier transform~\eqref{defFy} in explicit form. Since the func\-tion~$I$~\eqref{defI} has no poles in the strip~\eqref{xstrip}, we can make a contour shift
\begin{gather}
x \to z +(\mu+\nu)/2,\qquad x,z\in\R,
\end{gather}
to deduce
\begin{gather}\label{Fkap}
F(\mu,\nu;y)=\exp(i\alpha y(\mu+\nu)/2)\int_{\R}dz \exp(i\alpha zy)G(\pm z-\kappa),
\end{gather}
where we have set
\begin{gather}\label{kappa}
\kappa \equiv (\nu-\mu)/2,\qquad \im \kappa \in(0,a).
\end{gather}
Hence calculating $F$ amounts to f\/inding the cosine transform of $G(\pm z-\kappa)$ for $\kappa$ in the strip $\im \kappa\in(0,a)$. The reasoning in the proof of the following proposition, however, hinges on staying at f\/irst with~\eqref{defFy}.

A special case of the following result amounts to a hyperbolic analog of a trigonometric beta integral in Ramanujan's lost notebook. More precisely, \eqref{FFt} with $y=(\mu+\nu)/2$ amounts to equation~(1.8) with $\tau<0$ in Stokman's paper~\cite{stok}. Formula~\eqref{FFt} is also obtained in~Chapter~5 of van de Bult's Ph.D. thesis~\cite{bult2} (cf.~Theorem~5.6.8 with $n=1$), as a result of specializing more general formulas. In slightly dif\/ferent guises it occurred previously in various papers (the earliest ones being~\cite{faka,kash,pote}), but a complete proof cannot be found there.

\begin{prop}\label{propositionC.1}
The Fourier transform~\eqref{defFy} admits the explicit evaluation
\begin{gather}\label{FFt}
F(\mu,\nu;y)=(a_{+}a_{-})^{1/2}\exp(i\alpha y(\mu+\nu)/2)G(ia-2\kappa)G(\pm y-ia+\kappa),
\end{gather}
where $\mu,\nu$ satisfy~\eqref{munu}, $y$ satisfies \eqref{ystrip}, and $\kappa$ is given by~\eqref{kappa}.
\end{prop}
\begin{proof}
We f\/irst take $y$ real and choose $\mu,\nu$ such that
\begin{gather}\label{ass}
\im(\nu-\mu)\in  (a_l,a_{+}+a_{-}),\qquad a_l\equiv \max(a_{+},a_{-}).
\end{gather}
This entails that we can shift $\nu$ down by $ia_{\de}$ or $\mu$ up by $ia_{\de}$ without leaving the region~\eqref{munu}. Hence we may shift under the integral sign and use the $G$-A$\De$Es~\eqref{Gades} to obtain
\begin{gather}
F(\mu,\nu-ia_{\de};y)=2\int_{\R}dxe^{i\alpha xy}c_{-\de}(x-\nu+ia_{\de}/2)G(x-\nu)/G(x-\mu),\qquad y\in\R.
\end{gather}
Recalling that $F(\mu,\nu;y)$ is analytic in $y$ when~\eqref{ystrip} holds, this can be rewritten as
\begin{gather}
F(\mu,\nu-ia_{\de};y)=e_{-\de}(-\nu+ia_{\de}/2)F(\mu,\nu;y-ia_{\de}/2)\nonumber\\
\phantom{F(\mu,\nu-ia_{\de};y)=}{}
+e_{-\de}(\nu-ia_{\de}/2)F(\mu,\nu;y+ia_{\de}/2).\label{Fs1}
\end{gather}
Likewise, we obtain
\begin{gather}
F(\mu+ia_{\de},\nu;y)=e_{-\de}(-\mu-ia_{\de}/2)F(\mu,\nu;y-ia_{\de}/2)\nonumber\\
\phantom{F(\mu+ia_{\de},\nu;y)=}{}
+e_{-\de}(\mu+ia_{\de}/2)F(\mu,\nu;y+ia_{\de}/2).\label{Fs2}
\end{gather}

Next, consider the contour shift $x\to x+is$, where
\begin{gather}
s\in (\im\nu-a,\im \mu+a),
\end{gather}
cf.~\eqref{xstrip}. It implies
\begin{gather}
F(\mu,\nu;y)     =     \int_{\R}dx\exp(i\alpha (x+is)y)G(x+is-\nu)/G(x+is-\mu)
\nonumber \\
\phantom{F(\mu,\nu;y)     =}{}  =    \exp(-\alpha sy)F(\mu-is,\nu-is;y),\qquad y\in\R.\label{xshi}
\end{gather}
Since $F(\mu,\nu;y)$ is analytic in $y$ in the strip~\eqref{ystrip}, it now follows from~\eqref{ass} that we have
\begin{gather}\label{Fkey}
F(\mu-is,\nu-is;y)=e^{\alpha sy}F(\mu,\nu;y),\qquad |\im y|\le a_l/2.
\end{gather}

Now in~\eqref{Fs2} we may take $\mu,\nu\to \mu-is,\nu-is$ with $s\in(0,\im\mu +a)$. Then we are entitled to invoke~\eqref{Fkey} to get
\begin{gather}
F(\mu-is+ia_{\de},\nu-is;y)     =     e_{-\de}(-\mu+is-3ia_{\de}/2)e^{\alpha sy}F(\mu,\nu;y-ia_{\de}/2)
\nonumber \\
\phantom{F(\mu-is+ia_{\de},\nu-is;y)     =}{}   +e_{-\de}(\mu-is+3ia_{\de}/2)e^{\alpha sy}F(\mu,\nu;y+ia_{\de}/2).
\end{gather}
In this equation we can let $s$ converge to $a_{\de}$, which yields
\begin{gather}
F(\mu,\nu-ia_{\de};y)     =     e_{-\de}(2y-\mu-ia_{\de}/2)F(\mu,\nu;y-ia_{\de}/2)
\nonumber \\
\phantom{F(\mu,\nu-ia_{\de};y)     =}{}  +e_{-\de}(2y+\mu+ia_{\de}/2)F(\mu,\nu;y+ia_{\de}/2).\label{Fs3}
\end{gather}

Comparing~\eqref{Fs3} and~\eqref{Fs1}, we see that the dif\/ference yields a linear relation between $F(\mu,\nu;y+ia_{\de}/2)$ and $F(\mu,\nu;y-ia_{\de}/2)$. After some simplif\/ication, this relation can be rewritten as
\begin{gather}\label{Fade}
\frac{F(\mu,\nu;y+ia_{\de}/2)}{F(\mu,\nu;y-ia_{\de}/2)}=e_{-\de}(-\mu-\nu)\frac{c_{-\de}(y-ia+\kappa)}{c_{-\de}(y+ia-\kappa)}.
\end{gather}
Introducing
\begin{gather}\label{Fr}
F_r(\mu,\nu;y)\equiv \exp(i\alpha (\mu+\nu)y/2)G(\pm y-ia+\kappa),
\end{gather}
it is easy to verify that $F_r$ also satisf\/ies this $y$-A$\De$E. Since we can choose $\de=+,-$ and $a_{+}/a_{-}\notin\Q$, it readily follows that we must have
\begin{gather}\label{Fform}
F(\mu,\nu;y)=C(\mu,\nu)F_r(\mu,\nu;y),
\end{gather}
with $C$ independent of $y$.

The upshot of our reasoning thus far is that when $\mu$ and $\nu$ are restricted by~\eqref{ass}, then the Fourier transform is of the form~\eqref{Fform}, with $F_r$ given by~\eqref{Fr}. By analyticity in $\mu$ and $\nu$, this relation now extends to the whole region~\eqref{munu} and then, by analyticity in $y$, to the strip~\eqref{ystrip}. In view of~\eqref{Fkap}, we therefore have shown
\begin{gather}
\int_{\R}dz\exp(i\alpha zy)G(\pm z-\kappa)=C(\mu,\nu)G(\pm y-ia+\kappa),\\
 \im \kappa \in (0,a),\qquad |\im y|<\im \kappa.\nonumber
\end{gather}
It follows from this that we have
\begin{gather}
C(\mu,\nu)=D(\mu-\nu),
\end{gather}
so it remains to prove
\begin{gather}
D(\mu-\nu)=(a_{+}a_{-})^{1/2} G(ia+\mu-\nu).
\end{gather}

To this end we substitute
\begin{gather}
F(\mu,\nu;y)=D(\mu-\nu)\exp(i\alpha (\mu+\nu)y/2)G(\pm y-ia+\kappa)
\end{gather}
in~\eqref{Fs1}, and divide the result by
\begin{gather}
D(\mu-\nu)\exp(i\alpha (\mu+\nu)y/2)G(\pm y-ia_{\de}/2-ia+\kappa).
\end{gather}
Using the $G$-A$\De$Es, a straightforward calculation then yields
\begin{gather}\label{Dade}
\frac{D(\mu-\nu+ia_{\de})}{D(\mu-\nu)}=2is_{-\de}(\mu-\nu+ia_{\de}).
\end{gather}
Setting
\begin{gather}
D_r(\mu-\nu)\equiv G(\mu-\nu+ia),
\end{gather}
we see that $D_r$ also satisf\/ies the A$\De$E~\eqref{Dade}. Thus we deduce as before
\begin{gather}
D(\mu-\nu)=\eta D_r(\mu-\nu),
\end{gather}
where $\eta$ can only depend on $a_{+}$ and $a_{-}$.

As a result, we have now proved
\begin{gather}\label{Ftr}
\int_{\R}dx\exp(i\alpha xy)G(\pm x-\kappa)=\eta G(ia-2\kappa)G(\pm y-ia+\kappa),\\
 \im \kappa \in (0,a),\qquad |\im y|<\im \kappa.\nonumber
\end{gather}
To calculate $\eta$, we choose $\kappa$ equal to $ia_{-}/2$. Using the $G$-A$\De$Es, this yields
\begin{gather}
\int_{\R}dx\exp(i\alpha xy)\frac{1}{c_{+}(x)}=\eta G(ia_{+}/2-ia_{-}/2)\frac{1}{c_{-}(y)}.
\end{gather}
Now the integral is elementary, yielding $a_{+}/c_{-}(y)$. Using~\eqref{Geval} we then infer $\eta$ equals $(a_{+}a_{-})^{1/2}$, hence completing the proof.
\end{proof}

In view of~\eqref{Ftr}, the cosine transform of $G(\pm x -\kappa)$ is proportional to $G(\pm y-i\hat{\kappa})$, with
\begin{gather}
\hat{\kappa}\equiv ia-\kappa.
\end{gather}
To be specif\/ic, we have
\begin{gather}\label{ckk}
\left(\frac{2\alpha}{\pi}\right)^{1/2}\int_{0}^{\infty}\cos(\alpha xy)G(\pm x-\kappa)dx=G(\hat{\kappa}-\kappa)G(\pm y -\hat{\kappa}),\\
 \im\kappa\in(0,a),\qquad y\in\R.\nonumber
\end{gather}
Note that the proportionality constant can now be checked by taking the cosine transform of~\eqref{ckk}.

We proceed by deriving a corollary of the above result, which we need in~Subsection~\ref{section5.1}. First, we f\/ix $\nu$ in the strip $\im\nu\in(0,a)$, so that $G(x-\nu)$ has exponential decay for $x\to\pm\infty$, and choose $\mu$ real, so that $G(x-\mu)$ is a phase for real $x$. Consider now the integral
\begin{gather}\label{intau}
\left(\frac{\alpha}{2\pi}\right)^{1/2}\int_{\R}dx \exp(i\alpha xy)G(x-\nu)[\exp(-i\chi -i\alpha \mu^2/4+i\alpha\mu x/2)G(-x+\mu)],\\
 y\in\R. \nonumber
\end{gather}
In virtue of the $G$-asymptotics, the phase factor in square brackets converges to $\exp(i\alpha x^2/4)$ for $\mu\to -\infty$. By dominated convergence, the integral therefore converges to
\begin{gather}\label{mul1}
\left(\frac{\alpha}{2\pi}\right)^{1/2}\int_{\R}dx \exp(i\alpha xy)G(x-\nu)\exp(i\alpha x^2/4),\qquad  y\in\R ,\qquad \im\nu\in(0,a).
\end{gather}
On the other hand, by~\eqref{FFt} the integral~\eqref{intau} equals
\begin{gather}
\exp(-i\chi+i\alpha [y(\mu+\nu)/2+\mu\nu/4])\frac{G(\mu +ia-\nu)}{G(\mu+y+ia-\nu/2)}G(y-ia+\nu/2).
\end{gather}
Using once more the $G$-asymptotics, this has a $\mu\to -\infty$ limit that can be written as
\begin{gather}\label{mul2}
e^{-i\pi/4-4i\chi}G(y-ia+\nu/2)\exp\left( -\frac{i\alpha}{4}\big[(y-ia+\nu/2)^2+4(y+\nu/2)(ia-\nu)+\nu^2\big]\right).\!\!\!
\end{gather}

The resulting equality of~\eqref{mul1} and~\eqref{mul2} can be rewritten in a more illuminating way by using the $G$-cousins $G_R$~\eqref{GRDef} and $G_L$~\eqref{GLDef}. Indeed, let us set
\begin{gather}\label{defz}
z\equiv y+\nu/2,\qquad \im \nu\in(0,a),\qquad y\in\R,
\end{gather}
so that~\eqref{mul1} and~\eqref{mul2} can be written as
\begin{gather}
\exp(-i\alpha\nu^2/4)\left(\frac{\alpha}{2\pi}\right)^{1/2}\int_{\R}dx \exp(i\alpha xz)G(x-\nu)\exp\big(i\alpha (x-\nu)^2/4\big),
\end{gather}
and
\begin{gather}
\exp(-i\alpha\nu^2/4)e^{-i\pi/4-3i\chi}G_L(z-ia)\exp( -i\alpha z(ia-\nu)).
\end{gather}
Hence we have
\begin{gather}
\left(\frac{\alpha}{2\pi}\right)^{1/2}\int_{\R}dx \exp(i\alpha (x+ia-\nu)z)G(x-\nu)\exp\big(i\alpha (x-\nu)^2/4\big)\nonumber\\
\qquad{}
=e^{-i\pi/4-3i\chi}G_L(z-ia),
\end{gather}
with $z$ given by~\eqref{defz}. Choosing now $\nu =2is$,
we obtain the following corollary.

\begin{cor}\label{corollaryC2}
Letting
\begin{gather}
w=x+ia-2is,\qquad z=y+is,\qquad s\in(0,a/2),\qquad x,y\in\R,
\end{gather}
we have
\begin{gather}\label{cor}
\left(\frac{\alpha}{2\pi}\right)^{1/2}\int_{\R} \exp(i\alpha wz)G_R(w-ia)d\re w
=\exp(-i\pi/4-2i\chi)G_L(z-ia),
\end{gather}
where $\chi$ is given by~\eqref{chi}.
\end{cor}

Thus we have recovered the Fourier transform (A.21) in~\cite{Hels}.

\subsection*{Acknowledgments}

We would like to thank M.~Halln\"as for his interest and useful comments. We also thank several referees for pointing out additional references.

\LastPageEnding

\end{document}